\documentclass[a4paper, 9pt, twoside]{article}

\usepackage{amsmath,amssymb}
\usepackage{stmaryrd}
\usepackage{enumerate}
\usepackage{algorithm}
\usepackage[noend]{algpseudocode}
\usepackage{calrsfs}

\usepackage{graphicx}
\usepackage{amssymb}
\usepackage{array}
\usepackage{amsmath}
\usepackage{epstopdf}
\usepackage{epsfig}
\usepackage{subfigure}
\usepackage{pstricks}
\usepackage{fancyheadings}

\usepackage{stackengine}

\usepackage[margin=0.7in]{geometry}
\DeclareMathAlphabet{\pazocal}{OMS}{zplm}{m}{n}

\usepackage{graphicx}

\usepackage{stackengine}

\makeatletter
\def\BState{\State\hskip-\ALG@thistlm}
\makeatother

\newcommand{\lJump}{[\![}
\newcommand{\rJump}{]\!]}

\newtheorem{theorem}{Theorem}
\newtheorem{lemma}{Lemma}

\begin{document}
%%%%%%%%%%%%%%%%%%%%%%
%\journal{Applied Numerical Mathematics}
%\begin{frontmatter}
\title{On  energy stable discontinuous Galerkin spectral element approximations of the perfectly matched layer for  the wave equation}
\author{Kenneth Duru \thanks{Department of Geophysics, Ludwig-Maximilian University, Munich, Germany} , Alice-Agnes Gabriel \thanks{Department of Geophysics, Ludwig-Maximilian University, Munich, Germany} , Gunilla Kreiss \thanks{Division of Scientific Computing, Department of Information Technology, Uppsala University, Sweden}}
%
%\author[label1,label3]{Kenneth Duru}
%\author[label1]{Alice-Agnes Gabriel}
%\author[label2]{Gunilla Kreiss}
% \address[label1]{Department of Geophysics, Ludwig-Maximilian University, Munich, Germany}
%% \address[label2]{Department of Applied Mathematics, Naval Postgraduate School, Monterey, CA.}
% \address[label2]{Division of Scientific Computing, Department of Information Technology, Uppsala University, Sweden.}
% \address[label3]{Corresponding author: kenneth.duru@geophysik.uni-muenchen.de}
%  \author{Kenneth Duru \thanks{Department of Geophysics Stanford University, Stanford, CA.}, Gunilla Kreiss \thanks{Division of Scientific Computing Uppsala University, Uppsala, Sweden.}}
%%%%%%%%%%%%%%%%%%%%%%
%\author{}
\date{October 20, 2017}
%\address{}
%%%%%%%%%%%%%%%%%%%%%%
\maketitle

\pagenumbering{arabic}
%%%%%%%%%%%%%%%%%%%%%%%%%%%%%%

\begin{abstract}
%%%%
%%%%
In this paper, we develop a provably energy stable  discontinuous Galerkin spectral element method (DGSEM) approximation of the perfectly matched layer (PML) for the three and two space dimensional (3D and 2D) linear  acoustic wave equations, in first order form, subject to well-posed linear boundary conditions.  
%%%%
%%%%
First,  using the well-known complex coordinate stretching, we derive an efficient un-split modal PML for the 3D acoustic wave equation,  truncating a cuboidal computational domain.
%%%%
%%%% 
Second, we prove asymptotic stability of the continuous PML by deriving energy  estimates in the Laplace space, for the  3D PML in a heterogeneous acoustic medium,   assuming  piece-wise constant PML damping.
%%%%
%%%%
Third, we develop a DGSEM  for the  wave equation  using physically motivated numerical flux, with penalty weights, which are compatible with all well-posed, internal and external, boundary conditions. When the PML damping vanishes, by construction, our choice of penalty parameters yield an upwind scheme and a discrete energy estimate analogous to the continuous energy estimate. 
%%%%
%%%%
Fourth,  to ensure numerical stability of the discretization when PML damping is present, it is necessary to systematically extend the numerical numerical fluxes, and the  inter-element and boundary  procedures, to the PML auxiliary differential equations. This is critical  for deriving discrete energy estimates analogous to the continuous energy estimates.   
%%%%
Finally, we propose a procedure to compute  PML damping coefficients such that the PML error converges  to zero, at the  optimal convergence rate of the underlying numerical method. 
%%%%
Numerical solutions are evolved in time using the high order Taylor-type  time stepping  scheme of the same order of accuracy of the spatial discretization.  By combining the DGSEM spatial approximation with the high order Taylor-type  time stepping  scheme and the accuracy of the PML we obtain an arbitrarily  accurate wave propagation solver in the time domain.  Numerical experiments are presented in 2D and 3D corroborating the theoretical results.  
%%%%
%%%%
\end{abstract}
%%%%
%%%%
%\begin{keyword}
%%% keywords here, in the form: keyword \sep keyword
%acoustic waves \sep first order systems \sep perfectly matched layer \sep Laplace transforms \sep boundary conditions \sep interface conditions \sep stability \sep high order accuracy \sep discontinuous Galerkin spectral element method \sep penalty method.
%%% MSC codes here, in the form: \MSC code \sep code
%%% or \MSC[2008] code \sep code (2000 is the default)
%\end{keyword}

%\end{frontmatter}
%%%%%%%%%%%%%%%%%%%
%%%%%%%%%%%%%%%%%%%
\section{Introduction}\label{sec:s1}
%%%%%%%%%%%%%%%%%%%
%%%%%%%%%%%%%%%%%%%
The discontinuous Galerkin spectral element method (DGSEM) \cite{HesthavenWarburton2002, ReedHill1973}  is increasingly becoming  attractive as a method of choice for computing approximate solutions of partial differential equations (PDEs).  
%%%%%%%%%%%%%%%%%%%
%%%%%%%%%%%%%%%%%%%
%%%%%%%%%%%%%%%%%%%
%%%%%%%%%%%%%%%%%%%
%The DGSEM combines ideas from high order finite element methods with traditional finite volume methods, yielding local discrete operators with spectral accuracy. 
%%%%%%%%%%%%%%%%%%%
%%%%%%%%%%%%%%%%%%%%%%%
%%%%%%%%%%%%%%%%%%%%%%%
Computational procedures based on DGSEM can be flexible, high order accurate,  provably stable,  and well-suited  for complex large scale wave propagation problems \cite{delaPuenteAmpueroKaser2009, HesthavenWarburton2002, KoprivaGassner2014, Warburton2013, DumbserKaser2006, DumbserPeshkovRomenski}.  
%%%%%%%%%%%%%%%%%%%%%%%
%%%%%%%%%%%%%%%%%%%%%%%
%%%%%%%%%%%%%%%%%%%
However, real world wave propagation problems are often formulated in  large or unbounded domains.  
%%%%%%%%%%%%%%%%%%%
In numerical simulations, large domains must be replaced by smaller computational domains by introducing artificial boundaries. 
%%%%%%%%%%%%%%%%%%%
%%%%%%%%%%%%%%%%%%%
Thus, in order to retain the efficiency and high order accuracy of  DGSEM for wave propagation problems in unbounded domains, artificial boundaries introduced to limit the computational domain must be closed with reliable and accurate boundary conditions such that waves traveling out of the domain disappear without reflections. 
%%%%%%%%%%%%%%%%%%%
%%%%%%%%%%%%%%%%%%%
Otherwise, waves traveling out of the domain generate spurious reflections at artificial boundaries which will travel into the domain and pollute the solution everywhere. 
%%%%%%%%%%%%%%%%%%%

%%%%%%%%%%%%%%%%%%%
%{Chew1994}{DuKrSIAM}{AppeloKreiss2006} {Be_etAl}{Berenger1994} {DuruKozdonKreiss2016}{KDuru2016}
There are now two standard and competing approaches for effective artificial boundary closures, namely: high order local non-reflecting boundary condition (NRBC) \cite{Collino1993, Hagstrom2009, Engquist77, Givoli2004}, and absorbing layers \cite{Tago_etal, AppCol, Chew1994, Berenger1994}. A NRBC is a boundary condition defined at an artificial boundary such that little or no spurious reflections occur as a wave impinges the boundary.  All absorbing layers are constructed by modifying the underlying equations, in a layer surrounding the domain, such that waves  decay rapidly as they propagate in the layer. For this method to be effective, it is important that all waves traveling into the layer, independent of angles of incidence and frequency, be absorbed without reflections. This approach is analogous to the physical treatment of the walls of  anechoic chambers. Absorbing layers with these desirable attributes are called perfectly matched layers (PML) \cite{Berenger1994, Chew1994, KDuru2016}. In this paper, we will focus on the PML, because it is efficient, robust, and straightforward in the treatment of corners.
%%%%%%%%%%%%%%%%%%%
%%%%%%%%%%%%%%%%%%%

%%%%%%%%%%%%%%%%%%%
%%%%%%%%%%%%%%%%%%%
For the PML, or any  artificial boundary closure,  to be useful, it must be well-posed and stable. Otherwise, growth generated by the artificial boundary closure will propagate into the simulation domain and corrupt the solutions everywhere. The well-posedness and temporal stability analysis of the PML has been considered extensively in the literature \cite{DuKrSIAM, AppeloKreiss2006, Be_etAl, HalpernPetit-BergezRauch2011}.
For general systems, there is no guarantee that all solutions decay with time. In \cite{Be_etAl}, however, the geometric stability condition was introduced to characterize the temporal stability of initial value problems for PMLs. If this condition is not satisfied, then there are modes of high spatial frequencies with temporally growing amplitudes.  This result have been extended to PML initial boundary value problems (IBVPs) \cite{KDuru2016, DuKrSIAM}.  Even when the geometric stability condition is satisfied, however, numerical experiments have also shown that the PML can be unstable \cite{KaltenbacherKaltenbacherSim2013, AbeniusEdelvikJohansson2005}. For models that satisfy the geometric stability condition, like the acoustic wave equation, recent results \cite{KDuru2016, DuruKozdonKreiss2016} have revealed the impact of numerical boundary procedures on the  stability of discrete PMLs, using high order  summation-by-parts (SBP) finite difference method. 
%%%%%%%%%%%%%%%%%%%
%%%%%%%%%%%%%%%%%%%

%%%%%%%%%%%%%%%%%%%
%%%%%%%%%%%%%%%%%%%
By the results in \cite{DuKrSIAM, DuKr12, HalpernPetit-BergezRauch2011} the PML for the acoustic wave equation can be proven  well-posed and asymptotically stable. However, the PML and NRBC, involve auxiliary variables and equations that are often not covered by standard DGSEMs.
There are, though, a few exceptions \cite{ModaveAtleChanandWarburton2017, ModaveLambrechtsGeuzaine2016}. A straightforward DGSEM or finite element approximation for the PML for acoustic wave equation can result in an unstable scheme, \cite{KaltenbacherKaltenbacherSim2013, XieKomatitschMartinMatzen, AbeniusEdelvikJohansson2005}. We also refer the reader to section 5 of this paper for more elaborate numerical experiments and discussions. Many previous attempts to effectively include the PML in modern numerical methods, such the DGSEM for acoustic waves, suffered from numerical instabilities, particularly at long times \cite{KaltenbacherKaltenbacherSim2013, AbeniusEdelvikJohansson2005, Tago_etal, HalpernPetit-BergezRauch2011}. A major difficulty is that, in general, the PML (for symmetric systems) is asymmetric. Therefore, deriving energy estimates for the PML that are useful in designing stable and accurate DGSEM or finite element approximations  can be challenging.  Exponential and/or linear growth is often seen in numerical simulations.   Ad hoc procedures,  such as artificial numerical dissipation, were employed to tame numerical instabilities \cite{KaltenbacherKaltenbacherSim2013, XieKomatitschMartinMatzen, HalpernPetit-BergezRauch2011}. However, the stabilization procedures  can also introduce some undesirable effects, and destroy the fidelity of DGSEM solutions.
%%%%%%%%%%%%%%%%%%%
%%%%%%%%%%%%%%%%%%%

%%%%
%%%%

In this paper, we consider the PML for  3D and 2D linear  acoustic wave equations, in first order form, subject to well-posed linear boundary conditions.
%%%%
%%%%
Using the well-known complex coordinate stretching \cite{Chew1994}, we derive an efficient un-split modal PML for the 3D acoustic wave equation,  truncating a cuboidal computational domain.  In the Laplace space our PML is equivalent to the standard PML model. In the time-domain, however, we judiciously choose auxiliary variables, and invert the Laplace transforms. Our choice of auxiliary variables reduces dramatically the number of auxiliary variables and equations needed in the layer. For the 3D acoustic wave equation, only two auxiliary variables  are needed to surround a cuboidal domain with the PML. The  acoustic wave equation in 2D require only one  auxiliary variable, to surround a rectangular domain with the PML. 
%%%%

%%%%
Our primary objective, in the present study, is the development  of energy stable DGSEM approximations for the   PML in a bounded domain, subject to well-posed boundary conditions. 
%%%%
%%%%
 A main result in this study is the derivation of a weak form and an energy estimate for the continuous 3D PML  IBVP in the Laplace space,  assuming piecewise constant PML damping coefficients in a heterogenous acoustic medium. Thus,  proving asymptotic stability of the continuous PML.   The energy estimate is also valid for the 2D PML, and for the strip and corner regions of the PML. 
 %%%%
%%%%
Next, we develop a DGSEM approximation for the  PML  using physically motivated numerical flux, with penalty weights, which are compatible with all well-posed, internal and external, boundary conditions. When the PML damping vanishes, by construction, our choice of penalty parameters yield an upwind scheme and a discrete energy estimate analogous to the continuous energy estimate. 
%%%%
Another main result is the derivation of a discrete energy estimate for the  DGSEM discrete PML, in a 3D heterogenous acoustic medium, assuming element-wise constant PML damping.
%%%%
To do this, it is necessary to systematically extend the numerical numerical fluxes, and the  inter-element and boundary  procedures, to the PML auxiliary differential equations. This is critical  for deriving discrete energy estimates analogous to the continuous energy estimates, in the Laplace space. 
%%%%
%%%%
The semi-discrete energy estimate guarantees the stability and  accuracy of the DGSEM PML solution  at any future time.
%%%%
For polynomial approximations of degree $P$, the discrete energy estimate is valid for all quadrature rules that are exact for polynomial integrand of degree $2P - 1$.  Therefore, Gauss-Legendre-Lobatto (GLL),  Gauss-Legendre (GL) and Gauss-Legendre-Radau  (GLR)  quadrature rules  are qualified candidates.
%%%%

%%%%
%%%%
As soon the DGSEM approximation is introduced, the discrete PML is no longer a perfectly absorbing medium. The PML parameters must be tuned to achieve optimal performance.
%%%%
%%%%
Another objective of this paper is  the derivation of effective damping parameters, for DGSEM approximations of a finite width PML. 
%%%%
%%%%
We propose a procedure to compute  PML damping coefficients, for the DGSEM approximation of the PML,  such that the PML error converges  to zero, at the  optimal convergence rate of the underlying numerical method. 
%%%%
%%%%
Time integration is performed using  the high order Taylor-type  time stepping  scheme \cite{DumbserKaser2006, DumbserPeshkovRomenski} of the same order of accuracy as the spatial approximations.   
We present numerical experiments, in 2D and 3D, using Lagrange basis with GLL, GL and GLR quadrature nodes, separately.
The numerical experiments verify  accuracy,  stability and optimal convergence of PML errors.   
%%%%
%%%%

%%%%
%%%%
The remainder of the paper will proceed as follows. In section 2 we introduce the linear acoustic wave equation in 3D, and the corresponding PML equations. Continuous energy estimates are derived in section 3.  In section 4 we present numerical approximations and prove numerical stability. Numerical experiments are presented in section 5 verifying the analysis of previous sections. In section 6, we draw conclusions and suggest directions for future work.
%%%%%%%%%%%%%%%%%%%
%%%%%%%%%%%%%%%%%%%

\section{Equations}
%%%%%%%%%%%%%%%%%%%
%%%%%%%%%%%%%%%%%%%
In this section, we present the   3D linear acoustic wave equation, in a heterogeneous bounded domain. We introduce a general linear well-posed boundary conditions and derive energy estimates.  We conclude  the section by deriving the  modal PML via complex coordinate transformations. 
%%%%%%%%%%%%%%%%%%%
%%%%%%%%%%%%%%%%%%%
\subsection{The acoustic wave equation}
%%%%%%%%%%%%%%
Consider the  3D cuboidal domain
{
\small
\begin{equation}\label{eq:physical_domain}
\Omega = \{(x,y): -x_b \le x \le x_b, \quad -y_b \le y\le y_b, \quad -z_b \le z\le z_b; \quad x_b, y_b, z_b > 0 \},
\end{equation}
%%%%%%%%%%%%%%%%%%%%%%
}
with the linear  acoustic wave equation in first order form
%%%%%%%%%%%%%%%%%%%%%%
\begin{equation}\label{eq:linear_wave}
%%%%%%%%%%%%%%%%%%%%%%
\begin{split}
%%%%%%%%%%%%%%%%%%%
%%%%%%%%%%%%%%%%%%%
\frac{1}{\kappa} \frac{\partial p}{\partial t}  + \frac{\partial u}{\partial x} +  \frac{\partial v}{\partial y} +  \frac{\partial w}{\partial z} &= 0, \\
\rho\frac{\partial  u}{\partial t} + \frac{\partial p}{\partial x} &= 0, \\
\rho \frac{\partial v}{\partial t} +\frac{\partial p}{\partial y} &= 0, \\
\rho \frac{\partial w}{\partial t} + \frac{\partial p}{\partial z} &= 0,
%\mathbf{P}\frac{\partial{\mathbf{u}}}{\partial t} = \mathbf{A}\frac{\partial{\mathbf{u}}}{\partial x} + \mathbf{B}\frac{\partial{\mathbf{u}}}{\partial y},\quad  t> 0, \quad (x,y) \in \Omega, \quad \mathbf{u}(x,y,0) = \mathbf{u}^0(x,y) \in H^{s}, \quad s \ge 0,
%%%%%%%%%
%%%%%%%%%%%%%%%%%%%
%%%%%%%%%%%%%%%%%%%
  \end{split}
  \end{equation}
%%%%%%%%%%%%%%%%%%%
where the unknowns
%%%%%%%%%%%%%%%%%%%
are the acoustic pressure $p(x,y, z, t)$ and the velocity vector $[u(x,y, z, t), v(x,y, z,t), w(x,y, z,t)]^T$. 
%%%%%%%%%%%%%%%%%%%
Here  $\kappa(x,y, z) > 0$ is the bulk modulus, $\rho(x,y, z) > 0$ is the  density of the fluid. 
The acoustic wave speed is defined by $c = \sqrt{\kappa/\rho}$. 
%%%%%%%%%%%%%%%%%%%
%%%%%%%%%%%%%%%%%%%
%%%%%%%%%%%%%%%%%%%
Define the physical energy density,
\[
dE(x,y,z, t) = \frac{1}{2}\left[\frac{1}{\kappa}|p|^2 + {\rho}\left(|u|^2 + |v|^2 + |w|^2\right)\right] > 0,
\]
%%%%%%%%%%%%%%%%%%%
%%%%%%%%%%%%%%%%%%%
 and the mechanical energy $E(t)$, associated with the system  \eqref{eq:linear_wave}
{
\small
\begin{equation}\label{eq:physical_energy}
E(t) = \int_{\Omega}dE(x,y,z, t) dxdydz > 0.
\end{equation}
%%%%%%%%%%%%%%%%%%%
}
%%%%%%%%%%%%%%%%%%%
To show that the  Cauchy problem,  \eqref{eq:linear_wave} with the decay condition $| p|, |u|, |v|, |w| \to 0$ at $x_b \to \infty$, $y_b \to \infty$, $z_b \to \infty$,  is well-posed we multiply   \eqref{eq:linear_wave}  with $\boldsymbol{\phi}^T(x,y,z) $ from the left, where $\boldsymbol{\phi}(x,y,z) = \left({\phi}_p(x,y,z), {\phi}_u(x,y,z), {\phi}_v(x,y,z), {\phi}_w(x,y,z)\right) \in L^2\left(\Omega\right)$ is an arbitrary test function, and integrate over the whole spatial domain, $\Omega$,  we have
%%%%%%%%%%%%%%%%%%%
{
\small
\begin{equation}\label{eq:weak_wave}
%%%%%%%%%%%%%%%%%%%%%%
\begin{split}
%%%%%%%%%%%%%%%%%%%
%%%%%%%%%%%%%%%%%%%
\int_{\Omega}{\phi}_p\left(\frac{1}{\kappa} \frac{\partial p}{\partial t}  + \frac{\partial u}{\partial x} +  \frac{\partial v}{\partial y}\right)dxdydz = 0, \\
\int_{\Omega}{\phi}_u\left(\rho\frac{\partial u}{\partial t}  +\frac{\partial p}{\partial x}\right) dxdydz = 0, \\
\int_{\Omega}{\phi}_v\left(\rho \frac{\partial v}{\partial t} + \frac{\partial p}{\partial y}\right) dxdydz= 0, \\
\int_{\Omega}{\phi}_w\left(\rho \frac{\partial w}{\partial t} + \frac{\partial p}{\partial z}\right) dxdydz= 0. 
%%%%%%%%%%%%%%%%%%%
%%%%%%%%%%%%%%%%%%%
  \end{split}
  \end{equation}
%%%%%%%%%%%%%%%%%%%
}
Summing equation \eqref{eq:weak_wave} together we have
{
\small
%%%%%%%%%%%%%%%%%%%
\begin{align}\label{eq:product_1}
\int_{\Omega}\left(\frac{1}{\kappa} {\phi}_p \frac{\partial p}{\partial t}  +\rho\left( {\phi}_u\frac{\partial u}{\partial t}  + {\phi}_v\frac{\partial v}{\partial t} + {\phi}_w\frac{\partial w}{\partial t}\right)\right)dxdydz &= -\int_{\Omega} \left({\phi}_p\frac{\partial u}{\partial x}  + {\phi}_p \frac{\partial v}{\partial y} +  {\phi}_p \frac{\partial w}{\partial z}  + {\phi}_u \frac{\partial p}{\partial x} + {\phi}_v \frac{\partial p}{\partial y}  + {\phi}_w \frac{\partial p}{\partial z}  \right) dxdydz.
\end{align}
%%%%%%%%%%%%%%%%%%%
}
%%%%%%%%%%%%%%%%%%%
Integrating-by-parts the right hand side of \eqref{eq:product_1} gives
%%%%%%%%%%%%%%%%%%%
{
\small
%%%%%%%%%%%%%%%%%%%
\begin{equation}\label{eq:product_2}
\begin{split}
&\int_{\Omega}\left(\frac{1}{\kappa} {\phi}_p \frac{\partial p}{\partial t}  +\rho\left( {\phi}_u\frac{\partial u}{\partial t}  + {\phi}_v\frac{\partial v}{\partial t} + {\phi}_w\frac{\partial w}{\partial t} \right)\right)dxdydz = 
-\frac{1}{2}\int_{\Omega} \left({\phi}_u \frac{\partial p}{\partial x}- p\frac{\partial {\phi}_u}{\partial x}  + {\phi}_p\frac{\partial u}{\partial x} - u\frac{\partial {\phi}_p}{\partial x} \right) dxdydz\\
&-\frac{1}{2}\int_{\Omega} \left( {\phi}_v \frac{\partial p}{\partial y} - p\frac{\partial {\phi}_v}{\partial y}  + {\phi}_p \frac{\partial v}{\partial y}- v\frac{\partial {\phi}_p}{\partial y}  \right) dxdydz-\frac{1}{2}\int_{\Omega} \left({\phi}_w \frac{\partial p}{\partial z}- p\frac{\partial {\phi}_w}{\partial z}  + {\phi}_p \frac{\partial w}{\partial z}- w\frac{\partial {\phi}_p}{\partial z}\right) dxdydz\\
&-\frac{1}{2}\int_{y = -y_b}^{y = y_b}\int_{z = -z_b}^{z = z_b}\left(\left({\phi}_u p +{\phi}_p u\right)\Big|_{x=-x_b}^{x = x_b}\right)dydz -\frac{1}{2}\int_{x = -x_b}^{x = x_b}\int_{z = -z_b}^{z = z_b}\left(\left({\phi}_v p +{\phi}_p v\right)\Big|_{y=-y_b}^{y = y_b}\right)dxdz\\
& -\frac{1}{2}\int_{x = -x_b}^{x = x_b}\int_{y = -y_b}^{y = y_b}\left(\left({\phi}_w p +{\phi}_p w\right)\Big|_{z=-z_b}^{z = z_b}\right)dxdy\\
  \end{split}
\end{equation}
%%%%%%%%%%%%%%%%%%%
}
%%%%%%%%%%%%%%%%%%%
Replacing $\left(\phi_p, \phi_u, \phi_v, \phi_w\right)$ with $\left(p, u, v, w\right)$ in \eqref{eq:product_2}, in the right hand side  the volume terms vanish, having
%%%%%%%%%%%%%%%%%%%
{
\small
%%%%%%%%%%%%%%%%%%%
\begin{align}\label{eq:product_3}
\frac{d}{dt}\int_{\Omega}\frac{1}{2}\left[\frac{1}{\kappa}|p|^2 + \rho\left(|u|^2 + |v|^2 + |w|^2\right)\right]  dxdydz &= -\int_{-y_b}^{y_b}\int_{z = -z_b}^{z = z_b}\left[u(x_b, y, z, t) p(x_b, y, z, t) -  u(-x_b, y, z, t) p(-x_b, y, z,t) \right]dydz  \nonumber \\ 
&- \int_{-x_b}^{ x_b}\int_{z = -z_b}^{z = z_b} \left[v(x, y_b, z, t) p(x, y_b, z, t) -v(x, -y_b, z, t) p(x, -y_b, z, t)\right]dxdz\\ \nonumber
&- \int_{-x_b}^{ x_b}\int_{y = -y_b}^{y = y_b} \left[w(x, y, z_b, t) p(x, y, z_b, t) -w(x,  y, -z_b, t) p(x, y, -z_b, t)\right]dxdy.
\end{align}
%%%%%%%%%%%%%%%%%%%
}
%%%%%%%%%%%%%%%%%%%
 The decay condition, $| p|, |u|, |v|, |w| \to 0$ at $x_b \to \infty$, $y_b \to \infty$ and $z_b \to \infty$, yields the energy equation
%%%%%%%%%%%%%%%%%%%
{
\small
%%%%%%%%%%%%%%%%%%%
\begin{align}
\frac{d}{dt}E(t) = 0.
\end{align}
%%%%%%%%%%%%%%%%%%%
}
%%%%%%%%%%%%%%%%%%%
The energy is conserved, $E(t) = E(0)$ for all $t\ge 0$. 

\subsection{Boundary conditions}
%%%%%%%%%%%%%%%%%%%%%%
We will now consider a bounded domain.  Well-posed boundary conditions are needed to close the rectangular surfaces  of the boundaries of the cuboidal domain.
  Boundary conditions are enforced by modifying the amplitude of the incoming characteristics.  Thus, the number of boundary conditions must be equal to the number of incoming characteristics on the boundary.  In general,  boundary data for the incoming characteristics  can be expressed as a linear combination of the outgoing characteristics \cite{GustafssonKreissOliger1995, DuruGabrielIgel2017}. 
 
%%%%%
 To begin, we introduce the acoustic wave impedance
$
Z = \rho c,
$
where 
$
 c = \sqrt{\kappa/\rho} 
$
is the speed of sound.
%%%%%%%%%%%%%%%%%%%%%%
In the $x$-direction,  $y$-direction or $z$-direction, there are two characteristics, $\chi^{\left(\pm i\right)}, \quad i = x, y,z$, $\chi^{\left(- i\right)}$ propagating to the negative direction and the other $\chi^{\left(+ i \right)}$  propagating to the positive direction,
%%%%%%%%%%%%%%%%%%%%%%
\begin{align}\label{eq:characteristics}
\chi^{(\pm x)}: = \frac{1}{2}\left(Zu \mp p\right), \quad \chi^{(\pm y)}: = \frac{1}{2}\left(Zv \mp p\right), \quad \chi^{(\pm z)}: = \frac{1}{2}\left(Zw \mp p\right).
\end{align}
%%%%%%%%%%%%%%%%%%%%%%
Therefore,  at any boundary surface there is one incoming characteristic and one outgoing characteristic.
%%%%%%%%%%%%%%%%%%%%%%
We  pose the linear boundary condition, $\chi^{\left(\pm i\right)} = r_i \chi^{\left(\mp i\right)}$, at $\pm i_b$ where $r_i$ is the reflection coefficient. That is
%%%%%%%%%%%%%%%%%%%%%%
{
\small
%%%%%%%%%%%%%%%%%%%%%%
\begin{subequations}\label{eq:boundary_condition_acoustic}
\begin{alignat}{2}
%\begin{align}
\label{eq:bcx_acoustic}
\frac{1-r_x}{2}Zu \mp \frac{1+r_x}{2}p = 0, \quad \text{at} \quad x = \pm{x}_b,\\
%\end{align}
%%%%%%%%%%%%%%%%%%%
%%%%%%%%%%%%%%%%%%%
%\begin{align}
\label{eq:bcy_acoustic}
 \frac{1-r_y}{2} Zv \mp \frac{1+r_y}{2}p= 0, \quad \text{at} \quad y = \pm{y}_b,\\
%\label{eq:characteristics}
\label{eq:bcz_acoustic}
 \frac{1-r_z}{2} Zw \mp \frac{1+r_z}{2}p= 0, \quad \text{at} \quad z = \pm{z}_b.
%\label{eq:characteristics}
  \end{alignat}
 \end{subequations}
%%%%%%%%%%%
}

In \eqref{eq:boundary_condition_acoustic}, the non-dimensional real numbers $r_i$   with $i = x, y,z$,   and $|r_i| \le 1$ are reflection coefficients. The reflection coefficients model different physical situations. For example, we have $r_i = -1$: soft wall,  $r_i = 1$: hard wall, $r_i = 0$: absorbing,  boundary condition.
Thus from the energy equation \eqref{eq:product_3} we obtain
%%%%%%%%%%%%%%%%%%%
{
\small
%%%%%%%%%%%%%%%%%%%
\begin{align}\label{eq:product_4}
\frac{d}{dt}E(t) = -\int_{y = -y_b}^{y=y_b}\int_{z = -z_b}^{z=z_b}\mathbf{BT}^{(x)} dydz - \int_{x = -x_b}^{x=x_b}\int_{z = -z_b}^{z=z_b}\mathbf{BT}^{(y)} dxdz - \int_{x = -x_b}^{x=x_b}\int_{y = -y_b}^{y=y_b}\mathbf{BT}^{(z)} dxdy,
\end{align}
%%%%%%%%%%%%%%%%%%%
}
where
{
\small
\begin{align}\label{eq:boundary_terms}
\nonumber
&\mathbf{BT}^{(x)} = -\left(u(x_b, y, z, t) p(x_b, y, z,t) -  u(-x_b, y, z, t) p(-x_b, y, z, t)\right), \\
&  \mathbf{BT}^{(y)} = -\left(v(x, y_b, z, t) p(x, y_b, z, t) -v(x, -y_b, z, t) p(x, -y_b, z, t)\right), \\
&\mathbf{BT}^{(z)} = -\left(w(x, y, z_b, t) p(x, y, z_b, t) - w(x, y, -z_b, t) p(x, y, -z_b, t)\right).
\end{align}
}
It is easy to show, see \cite{DuruGabrielIgel2017}, that with the boundary conditions \eqref{eq:boundary_condition_acoustic} we have
%%%%%%%%%%%%%%%%%%%
%%%%%%%%%%%%%%%%%%%
\[
\mathbf{BT}^{(x)} = \left(\frac{1-|r_x(-x_b)|^2}{Z(-x_b,y, z)}|\chi^{(-x)}(-x_b, y, z)|^2 +\frac{1-|r_x(x_b)|^2}{Z(x_b, y, z)}|\chi^{(+x)}(x_b, y, z) |^2\right),
\]
\[
\mathbf{BT}^{(y)} = \left(\frac{1-|r_y(-y_b)|^2}{Z(x,-y_b, z)}|\chi^{(-y)}(x, -y_b, z)|^2 +\frac{1-|r_y(y_b)|^2}{Z(x, y_b, z)}|\chi^{(+y)}(x, y_b, z) |^2\right) ,
\]
\[
\mathbf{BT}^{(z)} = \left(\frac{1-|r_z(-z_b)|^2}{Z(x,y, -z_b)}|\chi^{(-z)}(x, y, -z_b)|^2 +\frac{1-|r_z(z_b)|^2}{Z(x, y, z_b)}|\chi^{(+z)}(x, y, z_b) |^2 \right),
\]
%%%%%%%%%%%%%%%%%%%
%%%%%%%%%%%%%%%%%%%
where $ \chi^{(\pm i)} $,  $i = x, y, z$,  are the characteristic variables defined in \eqref{eq:characteristics}.
%%%%%%%%%%%%%%%%%%%
%%%%%%%%%%%%%%%%%%%
 Thus if $|r_i| \le 1$ then  $\mathbf{BT}^{(i)} \ge 0$.
%%%%%%%%%%%%%%%%%%%
%%%%%%%%%%%%%%%%%%%
In particular, if  $|r_i^{(i)}| = 1$, then   from the boundary conditions \eqref{eq:boundary_condition_acoustic} we have $\mathbf{BT}^{(i)} = 0$.
%%%%%%%%%%%%%%%%%%%
%%%%%%%%%%%%%%%%%%%
Thus, the boundary terms are positive semi-definite, that is $\mathbf{BT}^{(i)} \ge 0$,
%$\mathbf{BT}^{(x)}(\mathbf{u}) \le 0, \quad \mathbf{BT}^{(y)}(\mathbf{u}) \le 0$.
%%%%%%%%%%%%%%%%%%%
%%%%%%%%%%%%%%%%%%%
%Enforcing the boundary conditions \eqref{eq:boundary_condition_acoustic}--\eqref{eq:boundary_condition_elasticity}   in the right hand side  of   \eqref{eq:product_4} and using  the identities defined \eqref{eq:identiy_1}--\eqref{eq:identiy_2} gives
%%%%%%%%%%%%%%%%%%%
%{
%\small
%\begin{align}\label{eq:BTs}
%\mathbf{BT}^{(x)} \le 0, \quad \mathbf{BT}^{(y)} \le 0,
%\end{align}
%%%%%%%%%%%%%%%%%%%%
%}
and
%%%%%
{
\small
\begin{align}\label{energy_conservation}
\frac{d}{dt}E(t) \le 0 \iff \mathrm{E}\left(t\right) \le \mathrm{E}\left(0\right), \quad \forall t \ge 0,
\end{align}
%%%%%%%%%%%%%%%%%%%
}
%%%%%%%%%%%%%%%%%%%
where the energy $\mathrm{E}(t)$ is defined by  \eqref{eq:physical_energy}.
%%%%%%%%%%%%%%%%%%%
%%%%%%%%%%%%%%%%%%%
\subsection{Interface conditions}
%%%%%%%%%%%%%%%%%%%
We formulate physical interface conditions that will be used to patch DGSEM elements together. We consider  planar interfaces and focus on  the  $x$-direction. The conditions can be easily extended to the $y$- and $z$-direction, and to curvilinear elements.
To begin we consider an interface at $x = 0$, and decompose the  domain into two subdomains having $\Omega = \Omega_{-x} \cup \Omega_{+x}$, with $\Omega_{-x} = [-x_b, 0]\times[-y_b, y_b]\times[-z_b, z_b]$ and $\Omega_{+x} = [0, x_b]\times[-y_b, y_b]\times[-z_b, z_b]$. We denote fields in the negative subdomain, $x < 0$, by superscript ${(-x)}$, and  fields in the positive subdomain, $x > 0$,  are denoted by superscript  ${(+x)}$. As before, the number of interface conditions must be equal to the number of  the number of characteristics going in and out of an interface.  We define jumps in the normal velocity across the interface, $\lJump {u}^{(x)}  \rJump = {u}^{(+x)}-{u}^{(-x)}$.
%%%%%%%%%%%%%%%%%%% 
%%%%%%%%%%%%%%%%%%%

%%%%%%%%%%%%%%%%%%%
%%%%%%%%%%%%%%%%%%%
 Since there are two characteristics going in and out of an interface, we formulate the interface condition,
%%%%%%%%%%%%%%%%%%%
%For later use, we summarize the interface condition:
%%%%%%%%%%%%%%%%%%%%%%%
{
\small
%%%%%%%%%%%%%%%%%%%%%%%
\begin{align}\label{eq:physical_interface_acoustic}
p^{(-x)} = p^{(+x)} = p^{(x)},  \quad \lJump  {u}^{(x)}  \rJump = 0, \quad x = 0.
\end{align}
%%%%%%%%%%%%%%%%%%%%%%%
}
The interface condition \eqref{eq:physical_interface_acoustic} enforces the continuity of the pressure field and the normal velocity across an interface.
%%%%%%%%%%%%%%%%%%%%%%%
%%%%%%%%%%%%%%%%%%%%%%%
 It is possible to generalize \eqref{eq:physical_interface_acoustic} to impedance conditions.
%%%%%%%%%%%%%%%%%%%%%%%
%\begin{align}\label{eq:physical_interface_acoustic_cont}
% P\lJump  {v}_x^{(x)}  \rJump = 0, \quad x = 0, \quad P\lJump  {v}_y^{(y)}  \rJump = 0, \quad y = 0,
%\end{align}
%%%%%%%%%%%%%%%%%%%%%%%%
%%%%%%%%%%%%%%%%%%%%%%%%
%\begin{align}\label{eq:physical_interface_maxwell_cont}
% E_z \lJump  {H}_y^{(x)}  \rJump = 0, \quad x = 0, \quad  E_z\lJump  {H}_x^{(y)}  \rJump = 0, \quad y = 0,
%\end{align}
%%%%%%%%%%%%%%%%%%%%
%%%%%%%%%%%%%%%%%%%%
%\begin{align}\label{eq:physical_interface_elastic_cont}
%  &\sigma_{xx} \lJump  {v}_x^{(x)}  \rJump = 0, \quad x = 0, \quad \sigma_{yy} \lJump  {v}_y^{(y)}  \rJump = 0, \quad y = 0,\\ \nonumber
%  &\sigma_{xy} \lJump  {v}_y^{(x)}  \rJump = 0, \quad x = 0, \quad \sigma_{xy} \lJump  {v}_x^{(y)}  \rJump = 0, \quad y = 0,
%\end{align}
%%%%%%%%%%%%%%%%%%%
%%%%%%%%%%%%%%%%%%%
The interface conditions \eqref{eq:physical_interface_acoustic} are both physically and and mathematically consistent. Using the energy method we have
{
\small
\begin{align}\label{eq:product_5}
\frac{d}{dt}E(t) = \int_{-y_b}^{y_b}\int_{-z_b}^{z_b}\mathbf{IT}^{(x)} dydz = 0,
\end{align}
%%%%%%%%%%%%%%%%%%%
}
%%%%%%%%%%%%%%%%%%%
where
{
\small
\begin{align}\label{eq:interface_term_x}
 \mathbf{IT}^{(x)}\left({p}^{(x)}, \lJump {u}^{(x)}  \rJump \right) = p^{(x)}\lJump  {u}^{(x)}  \rJump = 0.
\end{align}
 %%%%%%%%%%%%%%%%%%%
%%%%%%%%%%%%%%%%%%%
}
The interface terms vanish identically, and the energy is conserved
{
\small
\begin{align}\label{eq:product_6}
\frac{d}{dt}E(t) = 0.
\end{align}
}

In \cite{DuruGabrielIgel2017}, using interface conditions such as \eqref{eq:physical_interface_acoustic}, we developed new physics based numerical fluxes suitable for patching DGSEM elements together.   A fundamental step in the construction of the physically motivated numerical fluxes is to reformulate the boundary condition  \eqref{eq:boundary_condition_acoustic} and  interface condition \eqref{eq:physical_interface_acoustic} by introducing transformed (hat-) variables, $\left(\widehat{p}, \widehat{u}, \widehat{v}, \widehat{w}\right)$,  so that we can simultaneously construct (numerical) boundary/interface data for the particle velocities and the pressure.   The hat-variables encode the solution of the IBVP on the boundary/interface. To be more specific, the hat-variables are solutions of the Riemann problem constrained against physical boundary/interface conditions  \eqref{eq:boundary_condition_acoustic}  and \eqref{eq:physical_interface_acoustic}. Since the hat-variables are constructed to satisfy the boundary/interface conditions,  \eqref{eq:boundary_condition_acoustic}  and \eqref{eq:interface_term_x} exactly, we must have
%%%%%%
{
\small
\begin{align}\label{eq:identity_bc}
\mathbf{BT}^{(x)}(\widehat{p}, \widehat{u}) \ge 0, \quad \mathbf{BT}^{(y)}(\widehat{p}, \widehat{v}) \ge 0,  \quad \mathbf{BT}^{(z)}(\widehat{p}, \widehat{w}) \ge 0,
\end{align}
%%%%%%
\begin{align}\label{eq:identity_interface}
 \mathbf{IT}^{(x)}\left(\widehat{p}^{(x)}, \lJump \widehat{u}^{(x)}  \rJump \right)  = 0, \quad \mathbf{IT}^{(y)}\left(\widehat{p}^{(y)}, \lJump \widehat{v}^{(y)}  \rJump \right)  =0, \quad \mathbf{IT}^{(z)}\left(\widehat{p}^{(z)}, \lJump \widehat{w}^{(z)}  \rJump \right)  =0.
\end{align}
%%%%
}
 The indentities \eqref{eq:identity_bc}--\eqref{eq:identity_interface}, will be used in proving numerical stability. We refer the reader to \cite{DuruGabrielIgel2017} for more elaborate discussions.

\subsection{The PML for the wave equation in  first order form}
%%%%%%%%%
Here, will use the well known complex coordinate stretching technique \cite{Chew1994}, to construct a modal PML,  \cite{KDuru2016, DuruKozdonKreiss2016} for the  system \eqref{eq:linear_wave}. As above, we consider the homogeneous rectangular domain \eqref{eq:physical_domain}.
 %%%%%%%%%%%%%%%%%%
To begin with, let the Laplace transform, in time, of $\mathbf{u}\left(x,y,z,t\right)$ be defined by
{
\small
\begin{equation}
\widetilde{\mathbf{u}}(x,y, z, s)  = \int_0^{\infty}e^{-st}{\mathbf{u}}\left(x,y,z,t\right)\text{dt},  \quad s = a + ib, \quad \Re{s} = a > 0, \quad i = \sqrt{-1}.
\end{equation}
}
%%%%%%%%%%%%%%%%%%%%%%%%%%%%%%%%
We consider  a setup where the PML is included in  all space directions,  $x-$,  $y-$ and $z-$direction.  Take the Laplace transform,  in time,  of    equation  \eqref{eq:linear_wave}. 
The PML can be constructed direction-by-direction using 
$\partial/\partial{x} \to 1/S_x \partial/\partial{x} $, $\partial/\partial{y} \to 1/S_y \partial/\partial{y} $, $\partial/\partial{z} \to 1/S_z \partial/\partial{z} $, respectively. 
%%%%%%%%%%%%%%%%%%%
%%%%%%%%%%%%%%%%%%%
Here
%%%%%%%%%%%%%%%%%%%
{
\small
%%%%%%%%%%%%%%%%%%%
\begin{align}\label{eq:PML_metric}
S_x = 1 +\frac{d_x(x)}{s}, \quad  S_y = 1 + \frac{d_y(y)}{s},  \quad  S_z = 1 + \frac{d_z(z)}{s},
\end{align}
%%%%%%%%%%%%%%%%%%%
}
%%%%%%%%%%%%%%%%%%%
 are the complex PML metrics, with $s$ denoting the Laplace dual time variable,  and $d_x(x),d_y(y), d_z(z) \ge 0$  are  the damping functions. In the Laplace domain, $\partial/\partial t \to s$, the PML for \eqref{eq:linear_wave} truncating the computational cuboidal domain in all directions is 
 %%%%%%%%%%%%%%%%%%%
 \begin{equation}\label{eq:acoustic_pml_3D_Laplace}
\begin{split}
%%%%%%%%%%%%%%%%%%%
 %%%%%%%%%%%%%%%%%%%
\frac{1}{\kappa} s\widetilde{p} + \frac{1}{S_x}\frac{\partial \widetilde{u}}{\partial x} + \frac{1}{S_y}\frac{\partial \widetilde{v}}{\partial y} +\frac{1}{ S_z}\frac{\partial \widetilde{w}}{\partial z} & = 0,  \\
s\rho\widetilde{u} + \frac{1}{S_x}\frac{\partial \widetilde{p}}{\partial x} & = 0, \\
s\rho\widetilde{v} + \frac{1}{ S_y}\frac{\partial \widetilde{p}}{\partial y}& = 0,\\
s\rho\widetilde{w} + \frac{1}{ S_z}\frac{\partial \widetilde{p}}{\partial z}& = 0.\\
  \end{split}
  \end{equation}
  %%%%%%%%%%%%%%%%%%%
  %
  %%%%%%%%%%%%%%%%%%%
 The difficulty lies in transforming \eqref{eq:acoustic_pml_3D_Laplace} from the Laplace domain back to the time domain, without introducing
too many auxiliary variables. Following the technic in \cite{DuruKreiss2013}, we will judiciously choose auxiliary variables, and then invert the Laplace transforms. 
  %%%%%%%%%%%%%%%%%%%
  
  %%%%%%%%%%%%%%%%%%%
The time-dependent PML  will be  obtained by the following steps: Multiply the first and second equations in \eqref{eq:acoustic_pml_3D_Laplace}  with $S_x$,  and multiply the third and fourth equations in \eqref{eq:acoustic_pml_3D_Laplace} with $S_y$ and $S_z$, respectively. Choosing the auxiliary variables
  %%%%%%%%%%%%%%%%%%%
 \[
 s\widetilde{\sigma} = \left(d_x - d_y \frac{S_x}{S_y}\right)\frac{\partial \widetilde{v}}{\partial y}  , \quad  s\widetilde{\psi} = \left(d_x - d_z \frac{S_x}{S_z}\right)\frac{\partial \widetilde{w}}{\partial z},
 \]
 %%%%%%%%%%%%%%%%%%%
 and inverting the Laplace transforms we obtain
 %%%%%%%%%%%%%%%%%%%
 the time-dependent modal PML 
 {
\small
\begin{equation}\label{eq:acoustic_pml_3D_timedomain}
\begin{split}
%%%%%%%%%%%%%%%%%%%
%%%%%%%%%%%%%%%%%%%
\frac{1}{\kappa} \left(\frac{\partial p}{\partial t} + d_x p \right) +\frac{\partial u}{\partial x} +  \frac{\partial v}{\partial y} +\frac{\partial w}{\partial z} - \sigma - \psi & = 0, \\
\rho\left(\frac{\partial u}{\partial t} + d_x u \right) +\frac{\partial p}{\partial x} & = 0, \\
\rho \left(\frac{\partial v}{\partial t} + d_y v \right) +\frac{\partial p}{\partial y} & = 0,\\
\rho \left(\frac{\partial w}{\partial t} + d_z w \right) +\frac{\partial p}{\partial z} & = 0,\\
 \left(\frac{\partial \sigma}{\partial t} + d_y \sigma \right) + \left(d_y - d_x\right)\frac{\partial v}{\partial y} & = 0,\\
 \left(\frac{\partial \psi}{\partial t} + d_z \psi \right) + \left(d_z - d_x\right)\frac{\partial w}{\partial z}& = 0.
%\mathbf{P}\frac{\partial{\mathbf{u}}}{\partial t} = \mathbf{A}\frac{\partial{\mathbf{u}}}{\partial x} + \mathbf{B}\frac{\partial{\mathbf{u}}}{\partial y},\quad  t> 0, \quad (x,y) \in \Omega, \quad \mathbf{u}(x,y,0) = \mathbf{u}^0(x,y) \in H^{s}, \quad s \ge 0,
%%%%%%%%%
%%%%%%%%%%%%%%%%%%%
%%%%%%%%%%%%%%%%%%%
  \end{split}
  \end{equation}
%%%%%%%%%%%%%%%%%%%
%%%%%%%%%%%%%%%%%%%
}
%%%%%%%%%%%%%%%%%%%
%%%%%%%%%%%%%%%%%%%
%%%%%%%%%%%%%%%%%%%
There are certainly other ways to choose auxiliary variables; see, for example, \cite{KDuru2016}. However, all resulting PML models can be shown to be linearly equivalent to
\eqref{eq:acoustic_pml_3D_timedomain}. Since the PML model \eqref{eq:acoustic_pml_3D_timedomain} corresponds to \eqref{eq:acoustic_pml_3D_Laplace}, it follows that the equations
are perfectly matched to the acoustic wave equation \eqref{eq:linear_wave} by construction, see \cite{KDuru2016}.

 We will initialize all fields in  the PML with zero initial data and  terminate the PML  \eqref{eq:acoustic_pml_3D_timedomain} with the boundary conditions  \eqref{eq:boundary_condition_acoustic}.
%%%%%%%%%%
%%%%%%%%%%

Note  that the damping functions and auxiliary  functions vanish almost everywhere except in the layers defining the PML.
%%%%%%%%%%
%%%%%%%%%%
For example,  in the  $x$-dependent PML,  damping function $d_x(x)$  is nonzero, $d_x(x) > 0, d_y(y) =0, d_z(z) =0$,  only in the vertical PML layers truncating the left and right edges of the computational domain. In the multi-dimensional PML, there are edges and corner regions where two or all  damping functions are simultaneously nonzero, $d_x(x) > 0, d_y(y) >0, d_z(z) >0$.

 Without damping, $d_x(x)\equiv 0$,  $d_y(y)\equiv 0$,  $d_z(z)\equiv 0$, we recover the  wave equation  \eqref{eq:linear_wave},  which satisfies the energy estimate  \eqref{energy_conservation}.
%%%%%%%%%%
%%%%%%%%%%
However, the energy estimate is not  applicable to the PML, \eqref{eq:acoustic_pml_3D_timedomain},  when any of the damping functions, $d_x(x) > 0$,  $d_y(y) > 0$, or $d_z(z) > 0$, is non-zero.  Other technics are needed to investigate the well-posedness and stability properties of the PML.
%%%%%%%%% 
%%%%%%%%%
%%%%%%%%%

For the Cauchy PML problem, Fourier transform in space and Laplace transform in time yields no growing solutions \cite{DuKr12, AppeloKreiss2006, Be_etAl, HalpernPetit-BergezRauch2011}.  The analysis has been extended to PML IBVPs \cite{KDuru2016, DuKrSIAM}. However, these results  do not yield any type of energy estimate, and give no insight in  the construction of stable  DGSEM or finite element approximations of the PML in bounded computational domains. In the next  section, we will  derive energy estimates for the PML, suitable for developing energy stable  DGSEM or finite element approximations of the PML.
%%%%%%%%%

\section{Weak formulation and energy identity for the PML in the Laplace space}
%%%%%%%%%
%%%%%%%%%
Here, we formulate the corresponding weak formulation for the PML and derive energy estimates in the Laplace space subject to the boundary condition \eqref{eq:boundary_condition_acoustic}. The energy estimate establishes the well-posedness and the asymptotic stability of the continuous PML, and will be critical in developing an energy stable DGSEM approximation of the PML IBVP, \eqref{eq:acoustic_pml_3D_timedomain} , \eqref{eq:boundary_condition_acoustic}.

 Multiply the PML equation \eqref{eq:acoustic_pml_3D_timedomain} by the test function $\boldsymbol{\phi}^{T}(x,y,z)$, with
 \[
 \boldsymbol{\phi}(x,y,z) = \left({\phi}_p(x,y,z), {\phi}_u(x,y,z), {\phi}_v(x,y,z), {\phi}_w(x,y,z), {\phi}_{\sigma}(x,y,z), {\phi}_{\psi}(x,y,z)\right),
 \]
%%%%%%%%%
belonging to the mixed space
%%%%%%%%%
 \[
 \left({\phi}_p(x,y,z),  {\phi}_{\sigma}(x,y,z), {\phi}_{\psi}(x,y,z)\right) \in H^{1}\left(\Omega\right), \quad \left({\phi}_u(x,y,z),  {\phi}_{v}(x,y,z), {\phi}_{w}(x,y,z)\right) \in L^{2}\left(\Omega\right),
 \]
  and integrate over the whole domain, having
%%%%%%%%%
%%%%%%%%%
% 
\begin{equation}\label{eq:weak_form_pml_3D}
\begin{split}
\int_{\Omega}{\phi}_p \left( \frac{1}{\kappa} \left(\frac{\partial p}{\partial t} + d_x p \right)  +  \frac{\partial {u}}{\partial x} +  \frac{\partial {v}}{\partial y} +   \frac{\partial {w}}{\partial z}  - \sigma - \psi \right) dxdydz&= 0, \\
\int_{\Omega} {\phi}_u \left(\rho\left(\frac{\partial u}{\partial t} + d_x u \right)  + \frac{\partial {p}}{\partial x}  \right) dxdydz &= 0, \\
\int_{\Omega}{\phi}_v \left(\rho \left(\frac{\partial v}{\partial t} + d_y v \right)   + \frac{\partial {p}}{\partial y} \right)dxdydz  &= 0,\\
\int_{\Omega}{\phi}_w \left(\rho \left(\frac{\partial w}{\partial t} + d_z w \right) + \frac{\partial {p}}{\partial z} \right)dxdydz  &= 0,\\
 \int_{\Omega}{\phi}_{\sigma} \left(\frac{\partial \sigma}{\partial t} + d_y \sigma  + \left(d_y - d_x\right)\frac{\partial {v}}{\partial y}\right)dxdydz &= 0,\\
  \int_{\Omega}{\phi}_{\psi}  \left(\frac{\partial \psi}{\partial t} + d_y \psi  +  \left(d_z - d_x\right)\frac{\partial {w}}{\partial z}\right)dxdydz &= 0,
 %%
%\frac{1}{\kappa} s \widetilde{p}  &= -\frac{1}{S_x}\frac{\partial \widetilde{u}}{\partial x} -  \frac{1}{S_y}\frac{\partial \widetilde{v} }{\partial y} + \frac{1}{\kappa}f_p(x,y), \\
%\rho{ s\widetilde{u}} &= -\frac{1}{S_x}\frac{\partial \widetilde{p}}{\partial x} + \rho f_u(x,y), \\
%\rho s{\widetilde{v}} &= -\frac{1}{S_y}\frac{\partial \widetilde{p}}{\partial y} + \rho f_v(x,y),
%\mathbf{P} s\widetilde{\mathbf{u}} = \frac{1}{S_x} \mathbf{A}\frac{\partial{\widetilde{\mathbf{u}} }}{\partial x} + \frac{1}{S_y} \mathbf{B}\frac{\partial{\widetilde{\mathbf{u}} }}{\partial y} + \mathbf{P}{\mathbf{u}}^0, \quad \Re{s} > 0,
  \end{split}
  \end{equation}
   subject the boundary conditions \eqref{eq:boundary_condition_acoustic}.

Generalizing the $L_2$-norm to complex valued functions is necessary when dealing with the Laplace transformed systems. Let $a > 0$, be a positive real number and denote the complex variables, $u, v \in \mathbb{C}$. Define the weighted $L_2$ inner product and the corresponding norm
\begin{align}\label{eq:inner_product_complex}
\left(u, v\right)_{a} = \int_{\Omega}{v^*a u dxdydz}, \quad \|u\|_{a}^2 = \left(u, u\right)_{a}.
\end{align}
 Here $v^*$ denotes the complex conjugate of $v$.
 
 Let $\left({f}_p(x,y,z),  f_u(x,y,z),  f_v(x,y,z),   f_w(x,y,z), f_{\sigma}(x,y,z),   f_{\psi}(x,y,z)\right)^T$ denote the initial condition. Take the Laplace transform in time  of the PML equations \eqref{eq:acoustic_pml_3D_timedomain} and the boundary conditions \eqref{eq:boundary_condition_acoustic}, having
 \begin{equation}\label{eq:acoustic_pml_3D_Laplace0}
\begin{split}
\frac{s}{\kappa} S_x\widetilde{p}  + \frac{\partial \widetilde{u}}{\partial x} +  \frac{\partial \widetilde{v}}{\partial y} + \frac{\partial \widetilde{w}}{\partial z}  + \widetilde{\sigma}  + \widetilde{\psi}-  \frac{1}{\kappa} f_p(x,y,z) &=0, \\
\rho s{S_x} \widetilde{u} +  \frac{\partial \widetilde{p}}{\partial x} - \rho f_u(x, y, z)  &= 0, \\
 \rho s {S_y} \widetilde{v}  + \frac{\partial \widetilde{p}}{\partial y} - \rho f_v(x,y,z)  &= 0,\\
 \rho s {S_z} \widetilde{w}  + \frac{\partial \widetilde{p}}{\partial z} - \rho f_w(x,y,z) &= 0,\\
s{S_y}\widetilde{\sigma}  + \left(d_y - d_x\right)\frac{\partial \widetilde{v}}{\partial y} - f_{\sigma}(x,y,z)&= 0,\\
s{S_z}\widetilde{\psi}  + \left(d_z - d_x\right)\frac{\partial \widetilde{w}}{\partial z} - f_{\psi}(x,y,z)&= 0,\\
   \end{split}
  \end{equation}
 with 
  {
\small
%%%%%%%%%%%%%%%%%%%%%%
\begin{subequations}\label{eq:boundary_condition_acoustic_laplace}
\begin{alignat}{2}
%\begin{align}
\label{eq:bcx_acoustic}
\frac{1-r_x}{2}Z\widetilde{u} \mp \frac{1+r_x}{2}\widetilde{p} = 0, \quad \text{at} \quad x = \pm{x}_b,\\
%\end{align}
%%%%%%%%%%%%%%%%%%%
%%%%%%%%%%%%%%%%%%%
%\begin{align}
\label{eq:bcy_acoustic}
 \frac{1-r_y}{2} Z\widetilde{v} \mp \frac{1+r_y}{2}\widetilde{p}= 0, \quad \text{at} \quad y = \pm{y}_b,\\
%\label{eq:characteristics}
\label{eq:bcz_acoustic}
 \frac{1-r_z}{2} Z\widetilde{w} \mp \frac{1+r_z}{2}\widetilde{p}= 0, \quad \text{at} \quad z = \pm{z}_b.
  \end{alignat}
 \end{subequations}
%%%%%%%%%%%
}
 Here $s$, with $\Re{s} \ge 0$, is the dual time variable.
%%%%%%%%%%%
Note that 
%%%%%%%%%%%
\[
\frac{1-r_x}{2}Z \widetilde{p}^* \widetilde{u} = - \frac{1+r_x}{2}\widetilde{p}^{*}\widetilde{p} \le 0, \quad x = -x_b, \quad  
\frac{1-r_x}{2}Z\widetilde{p}^* \widetilde{u} =  \frac{1+r_x}{2}\widetilde{p}^{*}\widetilde{p} \ge 0, \quad x = x_b,
\]
%%%%%%%%%%%
%%%%%%%%%%%
\[
\frac{1-r_y}{2}Z \widetilde{p}^* \widetilde{v} = - \frac{1+r_y}{2}\widetilde{p}^{*}\widetilde{p} \le 0, \quad y = -y_b, \quad  
\frac{1-r_y}{2}Z \widetilde{p}^* \widetilde{v} =  \frac{1+r_y}{2}\widetilde{p}^{*}\widetilde{p} \ge 0, \quad y = y_b,
\]
%%%%%%%%%%%
%%%%%%%%%%%
\[
\frac{1-r_z}{2}Z \widetilde{p}^* \widetilde{w} = - \frac{1+r_z}{2}\widetilde{p}^{*}\widetilde{p} \le 0, \quad z = -z_b, \quad  
\frac{1-r_z}{2}Z \widetilde{p}^* \widetilde{w} =  \frac{1+r_z}{2}\widetilde{p}^{*}\widetilde{p} \ge 0, \quad z = z_b.
\]
%%%%%%%%%%%
Mutiply \eqref{eq:acoustic_pml_3D_Laplace0} by $\boldsymbol{\phi}^{H}(x,y,z)$ and integrate over the whole domain, we have
\begin{equation}\label{eq:acoustic_pml_3D_Laplace1}
\begin{split}
\int_{\Omega}{\phi}_p^*\left(\frac{s}{\kappa} S_x\widetilde{p}  + \frac{\partial \widetilde{u}}{\partial x} +  \frac{\partial \widetilde{v}}{\partial y} + \frac{\partial \widetilde{w}}{\partial z}  + \widetilde{\sigma}  + \widetilde{\psi}-  \frac{1}{\kappa} f_p(x,y,z) \right)dxdydz&=0, \\
\int_{\Omega}{\phi}_u^* \left(\rho s{S_x} \widetilde{u} +  \frac{\partial \widetilde{p}}{\partial x} - \rho f_u(x, y, z) \right)dxdydz &= 0, \\
\int_{\Omega}{\phi}_v^*\left( \rho s {S_y} \widetilde{v}  + \frac{\partial \widetilde{p}}{\partial y} - \rho f_v(x,y,z) \right)dxdydz &= 0,\\
\int_{\Omega}{\phi}_w^*\left( \rho s {S_z} \widetilde{w}  + \frac{\partial \widetilde{p}}{\partial z} - \rho f_w(x,y,z)\right)dxdydz &= 0,\\
 \int_{\Omega}{\phi}_{\sigma}^*\left(s{S_y}\widetilde{\sigma}  + \left(d_y - d_x\right)\frac{\partial \widetilde{v}}{\partial y} - f_{\sigma}(x,y,z)\right)dxdydz &= 0,\\
  \int_{\Omega}{\phi}_{\psi}^*\left(s{S_z}\widetilde{\psi}  + \left(d_z - d_x\right)\frac{\partial \widetilde{w}}{\partial z} - f_{\psi}(x,y,z)\right)dxdydz &= 0.\\
  \end{split}
  \end{equation}
  %%%%%%%%%%%%%%
%%%%%%%%%%%
Let $ j = x, y, z$, and introduce the complex numbers 
\begin{align}\label{eq:scaled_variables_2}
s = a + i b, \quad S_j =  \frac{s + d_j}{ s} = \frac{a+ d_j + ib}{ a + ib }, \quad |s| = \sqrt{a^2 + b^2}, \quad \Re{s}=a \ge 0.
\end{align} 
%%%%%%%%%%%
and the  real numbers
%%%%%%%%%%%
%\begin{align}\label{eq:scaled_variables_1}
% a^{\prime} = \frac{a}{|s|}, \quad b^{\prime} = \frac{b}{|s|}, \quad \epsilon_j = \frac{d_j}{|s|}, \quad \sqrt{{a^{\prime} }^2 + {b^{\prime} }^2} = 1, \quad 0 \le a^{\prime}  \le 1, 
%\end{align} 
%%%%%%%%%%%
%%%%%%%%%%%
%%%%%%%%%%%
\begin{align}\label{eq:scaled_variables_1}
0\le \Re\left(\frac{1}{S_j}\right) =  \frac{|s|^2 + d_j a}{|s|^2 + \left(d_j + 2a\right) d_j } \le 1, \quad 0 \le \gamma_j\left(s, d_j\right) := \Re\left(\frac{\left(sS_j\right)^*}{S_j}\right) =  \frac{a|s|^2 + d_j\left(2|s|^2 + d_j a \right)}{|s|^2 + d_j \left(d_j + 2a\right) }  \le |s|.
\end{align} 
When the PML damping vanishes $d_j = 0$, we have
$
\Re\left(\frac{1}{S_j}\right)  = 1, \quad \gamma_j\left(s, 0\right) = a.
$
%%%%
%%%%
If $ s = 0$, then
$
\Re\left(\frac{1}{S_j}\right)  = 0, \quad \gamma_j\left(0, d_j\right) = 0, 
$
for all $d_j \ge 0$. However, if $s \ne 0$ with $a \ge 0$ and $d_j> 0$, then $\gamma_j\left(s, d_j\right) > 0$.
%
%\[
%0< \Re\left(\frac{1}{S_j}\right) =  \frac{1 + \epsilon_j a^{\prime}}{1 + \left(\epsilon_j + 2a^{\prime}\right) \epsilon_j } \le 1, \quad \frac{1}{|s|} \Re\left(\frac{\left(sS_j\right)^*}{S_j}\right) =  \gamma_j(a^{\prime}, \epsilon_j)  > 0,
%\]
%with
%%%%%%%%%%%%
%\begin{align}\label{eq:scaled_variables_2}
%\gamma_j\left(a^{\prime}, \epsilon_j\right)  = \frac{a^{\prime} + \epsilon_j\left(2 + \epsilon_j a^{\prime} \right)}{1 + \epsilon_j \left(\epsilon_j + 2a^{\prime}\right) }, \quad 0 \le \gamma_j\left(a^{\prime}, \epsilon_j\right)   \le 1, \quad \forall a^{\prime},  \epsilon_j :  \quad 0 \le a^{\prime}  \le 1, \quad \epsilon_j  \ge 0.
%\end{align} 
%%%%%%%%%%%%
%When the PML damping vanishes, that is $\epsilon_j = 0$, we have  $\gamma_j(a^{\prime}, 0) =  a^{\prime}$. 
%%%%%%%%%%%%
%%%%%%%%%%%%
%If $a^{\prime} = 0$, then 
%$\gamma_j\left(0, \epsilon_j\right)  = {2 \epsilon_j }/{\left(1 + \epsilon_j^2\right)}$.
%%%%%%%%%%%%
%%%%%%%%%%%%
%Thus,  $\gamma_j\left(0, 0\right) = 0$, and if $a^{\prime} > 0$, or $ \epsilon_j > 0$, then $\gamma_j(a^{\prime}, \epsilon_j) > 0 $.
%%%%%%%%%%%
%
Introduce the complex vector variables
\[ 
\widetilde{\mathbf{U}}\left(s\right)   = \left[\widetilde{u}_1(s), \widetilde{u}_2(s), \widetilde{u}_3(s), \widetilde{u}_4(s)\right]^T, \quad \widetilde{\mathbf{V}}\left(s\right)  = \left[\widetilde{v}_1(s), \widetilde{v}_2(s), \widetilde{v}_3(s), \widetilde{v}_4(s) \right]^T,
\]
%%%%%%%%%%%
and define the weighted inner product
%%%%%%%%%%%
\begin{equation}\label{eq:elastic_pml_2D_20}
\begin{split}
&\Big\langle \widetilde{\mathbf{U}}\left(s\right) ,   \widetilde{\mathbf{V}}\left(s\right) \Big\rangle:= \left(\widetilde{u}_1,  \widetilde{v}_1\right)_{a/\kappa}+ \left(\widetilde{u}_2, \widetilde{v}_2\right)_{\gamma_x/\rho} + \left(\widetilde{u}_3, \widetilde{v}_3\right)_{\gamma_y/\rho}  + \left(\widetilde{u}_4, \widetilde{v}_4\right)_{\gamma_z/\rho},
  \end{split}
  \end{equation}
  with the corresponding energy-norm
 %%%%%%%%%%%
 %%%%%%%%%%%
\begin{align}\label{eq:energy-norm}
  \widetilde{E}\left(\widetilde{\mathbf{U}}\left(s\right) \right): = \sqrt{\Big\langle \widetilde{\mathbf{U}}\left(s\right) ,   \widetilde{\mathbf{U}}\left(s\right) \Big\rangle} \ge 0.
\end{align}
 %%%%%%%%%%%
 Note that if $a >0, d_j > 0$, then $ \widetilde{E}\left(\widetilde{\mathbf{U}}\left(s\right)  \right) > 0$ for all $\widetilde{\mathbf{U}}\left(s\right) \ne 0$.
 %%%%%%%%%%%
 We formulate a result peculiar to $ {a} \to 0$.
 %%%%%%%%% %%
\begin{lemma}\label{eq:remark1}
%%%%%%%%% %%%%%%%%%
Let 
%%%%%%%%% %%%%%%%%%
%%%%%%%%% %%%%%%%%%
$ 
\widetilde{\mathbf{U}}\left(s\right)
$
be a complex valued function, with 
%%%%%%%%% %%%%%%%%%
%%%%%%%%% %%%%%%%%%
%%%%%%%%% %%%%%%%%%
$ s = a + ib$,  $ {a} \ge 0$,   $b \in \mathbb{R}$,  and $d_j \ge 0$. 
%%%%%%%%% %%%%%%%%%
%%%%%%%%% %%%%%%%%%
Denote
%%%%%%%%% %%%%%%%%%
%%%%%%%%% %%%%%%%%%
$
\widetilde{E}^2\left(\widetilde{\mathbf{U}}\left(s\right)  \right) = \widetilde{E}^2_1\left(s\right)  +  \widetilde{E}^2_2\left(s\right), 
$
%%%%%%%%% %%%%%%%%%
with
%%%%%%%%% %%%%%%%%%
$
\widetilde{E}^2_1\left(s\right) =  \left(\widetilde{u}_1,  \widetilde{u}_1\right)_{a/\kappa}, \quad 
\widetilde{E}^2_2\left(s\right)=   \left(\widetilde{u}_2, \widetilde{u}_2\right)_{\gamma_x/\rho} + \left(\widetilde{u}_3, \widetilde{u}_3\right)_{\gamma_y/\rho}  + \left(\widetilde{u}_4, \widetilde{u}_4\right)_{\gamma_z/\rho},
$
 the norm induced by the inner product \eqref{eq:elastic_pml_2D_20}.
 We have
\begin{itemize}
%%%%%%%%% %%%%%%%%%
%%%%%%%%% %%%%%%%%%
\item[1)]  $0 \le \gamma_j(s, d_j) \le  |s| $ for all $a\ge 0$, $b \in \mathbb{R}$,  and $d_j \ge 0$. 
%%%%%%%%% %%%%%%%%%
\item[2)]  $s \to 0 \iff \gamma_j(s, d_j) \to 0   $ for all  $d_j \ge 0$. 
\item[3)]  $\widetilde{E}^2_1\left(s\right) \to 0$ for $ {a} \to 0$ and $\widetilde{E}^2_2\left(s\right) > 0$, for all $ {a} \ge 0$, $|s| > 0$, and $d_j >0$.
%%%%%%%%% %%%%%%%%%
%\item[3)] $\widetilde{E}\left(\widetilde{\mathbf{U}}\left(s\right)  \right)  \to 0$ if and only if  $\frac{\partial \widetilde{p} }{\partial x} \to 0$, $\frac{\partial \widetilde{p} }{\partial y} \to 0$, $\frac{\partial \widetilde{p} }{\partial z} \to 0$. 
%%%%%%%%%%%%%%%%%%
\item[4)] If $ {a} \to 0$ then $\widetilde{E}^2\left(\widetilde{\mathbf{U}}\left(s\right)  \right) \ge 0$ is a semi-norm, for all $|s| > 0$, and $d_j >0$.
\item[5)] If $ {a} > 0$ then $\widetilde{E}^2\left(\widetilde{\mathbf{U}}\left(s\right)  \right) > 0$ for all $\widetilde{\mathbf{U}}\left(s\right) \ne 0$,   thus defines a norm for all  $d_j >0$.
\end{itemize}
%%%%%%%%% %%%%%%%%%
\end{lemma}
%%%%%%%%% %%%%%%%%%
Thus, if $ {a} \to 0$ then the energy $\widetilde{E}\left(\widetilde{\mathbf{U}}\left(s\right)  \right) $ defined in \eqref{eq:energy-norm} is a semi-norm. This is not surprising since, with $ {a} \to 0$,  the time-derivative of the pressure field vanishes,  ${a} |s\widetilde{p}|^2 \to 0$. We note,  however, if 
 $ 
\widetilde{\mathbf{U}}\left(s\right)  = \left(s\widetilde{p}, \frac{1}{S_x} \frac{\partial \widetilde{p} }{\partial x}, \frac{1}{S_y} \frac{\partial \widetilde{p} }{\partial y}, \frac{1}{S_z} \frac{\partial \widetilde{p} }{\partial z}\right)^T,
$
%%%%%%%%% %%%%%%%%%
 then for $a\ge 0$, $b\ne 0$, $d_j > 0$,  we have $\widetilde{E}\left(\widetilde{\mathbf{U}}\left(s\right)  \right)  \to 0$ if and only if $ {a} \to 0$ and $\widetilde{p}(s, x, y,z)  \to \mathrm{const}.$
%%%%%%%%% %%%%%%%%%
 We have 
%%%%%%%%%%%%%%%%%%%
\begin{theorem}\label{Theo:Stability_PML_Laplace}
%%%%%%%%% %%%%%%%%%
Consider  the PML equation  in the Laplace space \eqref{eq:acoustic_pml_3D_Laplace0}  subject to the boundary conditions \eqref{eq:boundary_condition_acoustic_laplace} with $\Re{s}  = a \ge 0$ and piece wise constant damping  $d_x  \ge 0, d_y \ge 0,  d_z\ge 0$. Let 
%%%%%%%%% %%%%%%%%%
%%%%%%%%% %%%%%%%%%
\[ 
\widetilde{\mathbf{U}}\left(s\right)  = \left(s\widetilde{p}, \frac{1}{S_x} \frac{\partial \widetilde{p} }{\partial x}, \frac{1}{S_y} \frac{\partial \widetilde{p} }{\partial y}, \frac{1}{S_z} \frac{\partial \widetilde{p} }{\partial z}\right)^T, \quad 
%\]
%%%%%%%%% %%%%%%%%%
%\[ 
\widetilde{\mathbf{F}}\left(s\right)  = \left(\widetilde{F}_p(x,y,z), \frac{\rho}{S_x}   f_u(x,y,z),  \frac{\rho}{S_y}   f_v(x,y,z), \frac{\rho}{S_z}   f_w(x,y,z)\right)^T,
\]
%%%%%%%%% %%%%%%%%%
%%%%%%%%% %%%%%%%%%
with
%%%%%%%%% %%%%%%%%%
%%%%%%%%% %%%%%%%%%
\[
\widetilde{F}_p(x,y,z) = \frac{1}{S_x}f_{p}(x,y,z) - \frac{\kappa}{sS_yS_x}f_{\sigma}(x,y,z)  -  \frac{\kappa}{sS_zS_x}f_{\psi}(x,y,z).
\]
%%%%%%%%% %%%%%%%%%
%%%%%%%%% %%%%%%%%%
 We have
%%%%%%%%% %%%%%%%%%
%%%%%%%%% %%%%%%%%%
\begin{equation}\label{eq:energy_estimate_pml_laplace_corner_cont}
 \widetilde{E}^2\left( \widetilde{\mathbf{U}}\left(s\right)  \right) + \mathrm{BT}(s)    \le    \widetilde{E}\left( \widetilde{\mathbf{U}}\left(s\right)  \right)\widetilde{E}\left(\widetilde{\mathbf{F}}\left(s\right) \right), \quad \mathrm{BT}(s)   = {\Re\left(\frac{1}{S_x}\right)\widetilde{BT}^{(x)} + \Re\left(\frac{1}{S_y}\right)\widetilde{BT}^{(y)}  + \Re\left(\frac{1}{S_z}\right)\widetilde{BT}^{(z)} } \ge 0,
\end{equation}
%%%%%%%%% %%%%%%%%%
where 
{
\small
\[
\widetilde{BT}^{(x)}  = |s| \int_{-y_b}^{y_b} \int_{-z_b}^{z_b}{\widetilde{u}}{\widetilde{p}^*}\Big|_{-x_b}^{x_b} dydz \ge 0,
%\]
%%%%%%%%%%%%%%%%%%%
\quad
%%%%%%%%%%%%%%%%%%%
%\[
\widetilde{BT}^{(y)}  =|s| \int_{-x_b}^{x_b} \int_{-z_b}^{z_b}{\widetilde{v}}{\widetilde{p}^*}\Big|_{-y_b}^{y_b} dxdz \ge 0,
%\]
%%%%%%%%%%%%%%%%%%%
\quad
%%%%%%%%%%%%%%%%%%%
%\[
\widetilde{BT}^{(z)}  =|s| \int_{-x_b}^{x_b} \int_{-y_b}^{y_b}{\widetilde{w}}{\widetilde{p}^*}\Big|_{-z_b}^{z_b} dxdy \ge 0.
\]
}
\end{theorem}
%%%%%%%%% %%%%%%%%%
%%%%%%%%% %%%%%%%%%
The proof of Theorem \ref{Theo:Stability_PML_Laplace} can be found  in \ref{sec:appendix1}. 
%%%%%%%%% %%%%%%%%%

%%%%%%%%% %%%%%%%%%
The importance of Theorem \ref{Theo:Stability_PML_Laplace} is twofold: One, it establishes the well-posedness and asymptotic stability of the the PML IBVP, \eqref{eq:acoustic_pml_3D_timedomain}, \eqref{eq:boundary_condition_acoustic}, in a heterogeneous acoustic medium with piecewise constant damping function, $d_x  \ge 0, d_y \ge 0,  d_z\ge 0$. Two, it will be useful in designing provably  stable DGSEM approximations of the PML \eqref{eq:acoustic_pml_3D_timedomain} subject to the well-posed boundary conditions \eqref{eq:boundary_condition_acoustic}.
%%%%%%%%%%%%%%%%%%
%%%%%%%%%%%%%%%%%%%
%%%%%%%%%%%%%%%%%%%
%%%%%%%%%%%%%%%%%%%
\section{The discontinuous Galerkin spectral element method}
%%%%%%%%%%%%%%%%%%
In this section, we present the DGSEM approximations for the undamped system \eqref{eq:linear_wave} and the PML  \eqref{eq:acoustic_pml_3D_timedomain}, subject to the boundary conditions  \eqref{eq:boundary_condition_acoustic}. We will  use the physically motivated numerical fluxes develop in \cite{DuruGabrielIgel2017} to patch DGSEM elements into the global domain. The physically motivated numerical flux is upwind by construction and gives an energy estimate analogous to \eqref{energy_conservation}.  The boundary and inter-element procedure will begin with the integral form \eqref{eq:product_1}.
As we will see later, the procedure and analysis carries over when numerical approximations are introduced.
%%%%%%%%%%%%%%%%%%%
%%%%%%%%%%%%%%%%%%%

%%%%%%%%%%%%%%%%%%%
\subsection{Weak boundary and inter-element procedures, and the energy identity}
%%%%%%%%%%%%%%%%%%%
We begin by discretizing the domain $(x,y,z)  \in \Omega = [-x_b, x_b]\times[-y_b, y_b]\times[-z_b, z_b]$ into $L\times M\times N$ elements denoting the $lmn$-th element by $\Omega_{lmn} = [x_l, x_{l+1}]\times [y_m, y_{m+1}]\times [z_n, z_{n+1}]$, where $l = 1, 2, \dots, L$, $m= 1, 2, \dots, M$, $n = 1, 2, \dots, N$, with $x_1 = -x_b$,  $y_1 = -y_b$, $z_1 = -z_b$  and $x_{L+1} = x_b$, $y_{M+1} = y_b$, $z_{N+1} = z_b$. The volume integral yields
%%%%%%%%%%%%%%%%%%%%%%%
\[
\int_{\Omega} f(x,y,z) dxdydz = \sum_{k=1}^{K}\sum_{l=1}^{L}\sum_{m=1}^{M}\int_{\Omega_{klm}} f(x,y,z) dxdydz.
\]
Thus, the elemental weak form corresponding to \eqref{eq:weak_form_pml_3D} is
%%%%%%%%%%%%%%%%%%%
%%%%%%%%%%%%%%%%%%%
%\begin{equation}\label{eq:weak_form_pml_2D}
%\begin{split}
%\int_{\Omega}{\phi}_p \frac{1}{\kappa} \left(\frac{\partial p}{\partial t} + d_x p \right) dxdy&= - \int_{\Omega}{\phi}_p\left(\frac{\partial {u}}{\partial x} +  \frac{\partial \widetilde{v}}{\partial y} + \widetilde{w}  \right) dxdy, \\
%\int_{\Omega} {\phi}_u \rho\left(\frac{\partial u}{\partial t} + d_x u \right) dxdy &= -\int_{\Omega}{{\phi}_u}\frac{\partial {p}}{\partial x}  dxdy, \\
%\int_{\Omega}{\phi}_v \rho \left(\frac{\partial v}{\partial t} + d_y v \right) dxdy &= -\int_{\Omega} {{\phi_v}}\frac{\partial {p}}{\partial y}dxdy,\\
% \int_{\Omega}{\phi}_w  \left(\frac{\partial w}{\partial t} + d_y w \right)dxdy &= -\int_{\Omega}{\phi}_w\left(d_y - d_x\right)\frac{\partial {v}}{\partial y} dxdy,\\
% %%
%%\frac{1}{\kappa} s \widetilde{p}  &= -\frac{1}{S_x}\frac{\partial \widetilde{u}}{\partial x} -  \frac{1}{S_y}\frac{\partial \widetilde{v} }{\partial y} + \frac{1}{\kappa}f_p(x,y), \\
%%\rho{ s\widetilde{u}} &= -\frac{1}{S_x}\frac{\partial \widetilde{p}}{\partial x} + \rho f_u(x,y), \\
%%\rho s{\widetilde{v}} &= -\frac{1}{S_y}\frac{\partial \widetilde{p}}{\partial y} + \rho f_v(x,y),
%%\mathbf{P} s\widetilde{\mathbf{u}} = \frac{1}{S_x} \mathbf{A}\frac{\partial{\widetilde{\mathbf{u}} }}{\partial x} + \frac{1}{S_y} \mathbf{B}\frac{\partial{\widetilde{\mathbf{u}} }}{\partial y} + \mathbf{P}{\mathbf{u}}^0, \quad \Re{s} > 0,
%  \end{split}
%  \end{equation}
 %%%%%%%%%%%%%%%%%%%%%%%
 %%%%%%%%%%%%%%%%%%%%%%%
\begin{align}\label{eq:weak_0}
\int_{\Omega_{klm}} {\phi}_p \left(\frac{1}{\kappa} \left(\frac{\partial p}{\partial t} + d_x p \right)+ \frac{\partial {u}}{\partial x} +  \frac{\partial {v}}{\partial y} + \frac{\partial {w}}{\partial z}  + \sigma + \psi \right) dxdydz&= 0,\\
%%%%%%%
\int_{\Omega_{klm}} {\phi}_u\left( \rho\left(\frac{\partial u}{\partial t} + d_x u \right) + \frac{\partial {p}}{\partial x} \right) dxdydz &=  0, \\
%%%%%
\int_{\Omega_{klm}}{\phi}_v\left( \rho \left(\frac{\partial v}{\partial t} + d_y v \right) + \frac{\partial {p}}{\partial y} \right) dxdydz &= 0,\\
%%%%%
\int_{\Omega_{klm}}{\phi}_w \left(\rho \left(\frac{\partial w}{\partial t} + d_z w \right) + \frac{\partial {p}}{\partial z}\right)dxdydz &= 0,\\
%%%%%
\int_{\Omega_{klm}}{\phi}_{\sigma}  \left(\frac{\partial \sigma }{\partial t} + d_y \sigma  + \left(d_y - d_x\right)\frac{\partial {v}}{\partial y}\right) dxdydz &= 0,\\
%%%%%
\int_{\Omega_{klm}}{\phi}_{\psi}  \left(\frac{\partial \psi }{\partial t} + d_z \psi + \left(d_z - d_x\right)\frac{\partial {w}}{\partial z}\right)dxdydz &= 0.
%%%%%
\end{align}
%%%%%%%%%%%%%%%%%%%%%%%
%%%%%%%%%%%%%%%%%%%%%%%
Next   we consider the element boundaries, $x = x_l, x_{l+1}$, $y = y_m, y_{m+1}$, $z = z_n, z_{n+1}$, and generate boundary and interface data $\widehat{p}(x, y, z, t)$, $ \widehat{u} (x, y,z, t)$, at $x = x_l, x_{l+1}$, $\widehat{p}(x, y, z, t)$, $ \widehat{v} (x, y, z, t)$, at $y = y_m, y_{m+1}$ and $\widehat{p}(x, y, z, t)$, $ \widehat{w} (x, y,z, t)$, at $z = z_n, z_{n+1}$. The hat-variables encode the solutions of the IBVP at the element boundaries. Please see \cite{DuruGabrielIgel2017} for more elaborate discussions.
%%%%%%%%%%%%%%%%%%%%%%%
%%%%%%%%%%%%%%%%%%%%%%%
Note that, by both physical and mathematical considerations, the only way information can be propagated into an element is through the incoming characteristics,  $ \chi^{\pm i} $,  $i = x, y, z$ defined in \eqref{eq:characteristics}, on the boundaries. We construct flux fluctuations by penalizing data against incoming characteristics only,
%%%%%%%%%%%%%%%%%%%%%%%
{
\small
\begin{equation}
\begin{split}
 F^x(x_k, y, z, t):=\chi^{(-x)} -\widehat{\chi}^{(-x)} &=\frac{Z_{s}(x_k, y, z)}{2} \left(u(x_k, y, z, t)-\widehat{u}(x_k, y, z, t) \right) + \frac{1}{2}\left(p(x_k, y, z,t)- \widehat{p} (x_k, y, z,t)\right),\\
 F^y(x, y_l, z, t):=\chi^{(-y)} -\widehat{\chi}^{(-y)} &=\frac{Z_{s}(x, y_l, z)}{2} \left(v(x, y_l, z, t)-\widehat{v}(x, y_l, z,t) \right) + \frac{1}{2}\left(p(x, y_l, z,t)- \widehat{p} (x, y_l, z,t)\right),\\
 F^{z}(x, y, z_m, t):=\chi^{(-z)} -\widehat{\chi}^{(-z)} &=\frac{Z_{s}(x, y, z_m)}{2} \left(w(x, y, z_m, t)-\widehat{w}(x, y, z_m ,t) \right) + \frac{1}{2}\left(p(x, y, z_m, t)- \widehat{p} (x, y, z_m, t)\right),
 \end{split}
\end{equation}
%%%%%%%%%%%%%%%%%%%%%%%
%%%%%%%%%%%%%%%%%%%%%%%
%%%%%%%%%%%%%%%%%%%%%%%
\begin{equation}
\begin{split}
 G^x(x_{k+1}, y, z, t):=\chi^{(+x)} -\widehat{\chi}^{(+x)} &=\frac{Z_{s}(x_{k+1}, y, z)}{2} \left(u(x_{k+1}, y, z, t)-\widehat{u}(x_{k+1}, y, z, t) \right) - \frac{1}{2}\left(p(x_{k+1}, y, z,t)- \widehat{p} (x_{k+1}, y, z,t)\right),\\
 G^y(x, y_{l+1}, z, t):=\chi^{(+y)} -\widehat{\chi}^{(+y)} &=\frac{Z_{s}(x, y_{l+1}, z)}{2} \left(v(x, y_{l+1}, z, t)-\widehat{v}(x, y_{l+1}, z,t) \right) - \frac{1}{2}\left(p(x, y_{l+1}, z,t)- \widehat{p} (x, y_{l+1}, z,t)\right),\\
 G^{z}(x, y, z_{m+1}, t):=\chi^{(+z)} -\widehat{\chi}^{(+z)} &=\frac{Z_{s}(x, y, z_{m+1})}{2} \left(w(x, y, z_{m+1}, t)-\widehat{w}(x, y, z_{m+1} ,t) \right) - \frac{1}{2}\left(p(x, y, z_{m+1}, t)- \widehat{p} (x, y, z_{m+1}, t)\right).
 \end{split}
\end{equation}
%%%%%%%%%%%%%%%%%%%%%%%
}
The, fluctuations, $F^{i}$, $G^{i}$, penalize data against the incoming characteristic at the element boundaries. Since we have not introduced any approximation yet,  the hat-variables are exact and $F^{i}\equiv 0$, $G^{i}\equiv 0$.  In particular, at the external boundaries we have
%%%%%%%%%%%%%%%%%%%%
%%%%%%%%%%%%%%%%%%%%
{
\small
%%%%%%%%%%%%%%%%%%%%%%
\begin{subequations}\label{eq:boundary_condition_acoustic_hat}
\begin{alignat}{2}
%\begin{align}
\label{eq:bcx_acoustic_hat}
F^{x}:=\frac{1-r_x}{2}Zu + \frac{1+r_x}{2}p = 0, \quad \text{at} \quad x = -{x}_b,
\quad 
G^{x}:=\frac{1-r_x}{2}Zu - \frac{1+r_x}{2}p = 0, \quad \text{at} \quad x = +{x}_b, \\
%\end{align}
%%%%%%%%%%%%%%%%%%%
%%%%%%%%%%%%%%%%%%%
%\begin{align}
\label{eq:bcy_acoustic_hat}
 F^{y}:=\frac{1-r_y}{2}Zv + \frac{1+r_y}{2}p = 0, \quad \text{at} \quad y = -{y}_b,
\quad 
G^{y}:=\frac{1-r_y}{2}Zv - \frac{1+r_y}{2}p = 0, \quad \text{at} \quad y = +{y}_b, \\
%\label{eq:characteristics}
\label{eq:bcz_acoustic_hat}
 F^{z}:=\frac{1-r_z}{2}Zw + \frac{1+r_z}{2}p = 0, \quad \text{at} \quad z = -{z}_b,
\quad 
G^{z}:=\frac{1-r_z}{2}Zw - \frac{1+r_z}{2}p = 0, \quad \text{at} \quad z = +{z}_b.
%\label{eq:characteristics}
  \end{alignat}
 \end{subequations}
 %%%%%%%%%%%%%%%%%%%
%%%%%%%%%%%%%%%%%%%
%%%%%%%%%%%
}

Consequently we append the flux fluctuations ${F}^{i} = 0$, ${G}^{i} = 0$ to the weak form \eqref{eq:weak_0}, and we have
%%%%%%%%%%%%%%%%%%%
{
\small
\begin{align}\label{eq:weak_1_eq_1}
&\int_{\Omega_{lmn}} {\phi}_p \left(\frac{1}{\kappa} \left(\frac{\partial p}{\partial t} + d_x p \right)+ \frac{\partial {u}}{\partial x} +  \frac{\partial {v}}{\partial y} + \frac{\partial {w}}{\partial z}  + \sigma + \psi \right) dxdydz= \\ \nonumber
&- \int_{y_m}^{y_{m+1}}\int_{z_n}^{z_{n+1}}{\left(\frac{\phi_p(x_l, y, z)}{Z_s(x_{l}, y,z)} F^{x}(x_l,y, z, t) -  \frac{\phi_p(x_{l+1}, y, z)}{Z_s(x_{l+1}, y, z)}G^{x}(x_{l+1}, y, z, t)\right)}dydz\\ \nonumber
%%%%%%
& - \int_{x_l}^{x_{l+1}}\int_{z_n}^{z_{n+1}}{\left( \frac{\phi_p(x, y_m, z)}{Z_s(x, y_m, z)}F^{y}(x, y_m, z, t) -  \frac{\phi_p(x, y_{m+1}, z)}{Z_s(x, y_{m+1}, z)}G^{y}(x, y_{m+1}, z, t)\right)}dxdz\\ \nonumber
%%%%%
& - \int_{x_l}^{x_{l+1}}\int_{y_m}^{y_{m+1}}{\left( \frac{\phi_p(x, y, z_n)}{Z_s(x, y, z_n)}F^{z}(x, y, z_n, t) -  \frac{\phi_p(x, y, z_{n+1})}{Z_s(x, y, z_{n+1})}G^{y}(x, y, z_{n+1}, t)\right)}dxdy,
%%%%%%%
\end{align}
\begin{align}
\label{eq:weak_1_eq_2}
&\int_{\Omega_{lmn}} {\phi}_u\left( \rho\left(\frac{\partial u}{\partial t} + d_x u \right) + \frac{\partial {p}}{\partial x} \right) dxdydz =  \\ \nonumber
%%%%%
& - \int_{y_m}^{y_{m+1}}\int_{z_n}^{z_{n+1}}{\left({\phi_u(x_l, y, z)} F^{x}(x_l,y, z, t) +  {\phi_u(x_{l+1}, y, z)}G^{x}(x_{l+1}, y, z, t)\right)}dydz, 
%%%%%
\end{align}
\begin{align}
\label{eq:weak_1_eq_3}
&\int_{\Omega_{lmn}}{\phi}_v\left( \rho \left(\frac{\partial v}{\partial t} + d_y v \right) + \frac{\partial {p}}{\partial y} \right) dxdydz = \\ \nonumber
%%%%%
& - \int_{x_l}^{x_{l+1}}\int_{z_n}^{z_{n+1}}{\left({\phi_v(x, y_m, z)}F^{y}(x, y_m, z, t) +  {\phi_v(x, y_{m+1}, z)}G^{y}(x, y_{m+1}, z, t)\right)}dxdz,
%
%%%%%
\end{align}
\begin{align}
\label{eq:weak_1_eq_4}
&\int_{\Omega_{lmn}}{\phi}_w \left(\rho \left(\frac{\partial w}{\partial t} + d_z w \right) + \frac{\partial {p}}{\partial z}\right)dxdydz = \\ \nonumber
%%%%%
& - \int_{x_l}^{x_{l+1}}\int_{y_m}^{y_{m+1}}{\left({\phi_w(x, y, z_n)}F^{z}(x, y, z_n, t) +  {\phi_w(x, y, z_{n+1})}G^{y}(x, y, z_{n+1}, t)\right)}dxdy,
%%%%%
\end{align}
\begin{align}
\label{eq:weak_1_eq_5}
& \int_{\Omega_{lmn}}{\phi}_{\sigma}  \left(\frac{\partial \sigma }{\partial t} + d_y \sigma  + \left(d_y - d_x\right)\frac{\partial {v}}{\partial y}\right) dxdydz = \\ \nonumber
& - \underbrace{\int_{x_l}^{x_{l+1}}\int_{z_n}^{z_{n+1}}{\omega_y \left(d_y - d_x\right)\left( \frac{\phi_{\sigma}(x, y_m, z)}{Z_s(x, y_m, z)}F^{y}(x, y_m, z, t) -  \frac{\phi_{\sigma}(x, y_{m+1}, z)}{Z_s(x, y_{m+1}, z)}G^{y}(x, y_{m+1}, z, t)\right)}dxdz}_{\text{PML stabilizing flux fluctuation} = 0},
%%%%%
\end{align}
\begin{align}
\label{eq:weak_1_eq_6}
& \int_{\Omega_{lmn}}{\phi}_{\psi}  \left(\frac{\partial \psi }{\partial t} + d_z \psi + \left(d_z - d_x\right)\frac{\partial {w}}{\partial z}\right)dxdydz = \\ \nonumber
& -  \underbrace{\int_{x_l}^{x_{l+1}}\int_{y_m}^{y_{m+1}}{\omega_z\left(d_z - d_x\right)\left( \frac{\phi_{\psi}(x, y, z_n)}{Z_s(x, y, z_n)}F^{z}(x, y, z_n, t) -  \frac{\phi_{\psi}(x, y, z_{n+1})}{Z_s(x, y, z_{n+1})}G^{y}(x, y, z_{n+1}, t)\right)}dxdy}_{\text{PML stabilizing flux fluctuation} = 0}.
%%%%%%
\end{align}
%%%
}
%%%%%%
Here, $\omega_y$ and $\omega_z$ are PML stabilizing parameters to be determined by requiring stability of the discrete PML.
%%%%%%
%%%%%%
Note that since ${F}^{i} = 0$, ${G}^{i} = 0$, equation \eqref{eq:weak_1_eq_1}--\eqref{eq:weak_1_eq_6} is completely equivalent to \eqref{eq:weak_0}. We will now derive an energy estimate analogous to \eqref{energy_conservation}. To do this we set the damping functions to zero, $d_x = d_y = d_z = 0$, the auxiliary variables terms vanish, $\sigma = \psi = 0$, and equations \eqref{eq:weak_1_eq_5}--\eqref{eq:weak_1_eq_6} drop out.  We follow exactly the same steps as in \cite{DuruGabrielIgel2017}. That is, we integrate-by-parts the volume term in the right hand side of \eqref{eq:weak_1_eq_1}--\eqref{eq:weak_1_eq_4}, obtaining
%%%

\begin{theorem}\label{theorem:main_1}
Consider the weak form \eqref{eq:weak_1_eq_1}--\eqref{eq:weak_1_eq_6} of the PML equation. When the damping functions  vanish, $d_x = d_y = d_z =0$, We have the energy identity
\begin{equation}\label{eq:weak_energy_elements}
\begin{split}
&\frac{d}{dt} {E}(t) = -\sum_{l=1}^{L}\sum_{m=1}^{M}\sum_{n=1}^{N}\int_{y_m}^{y_{m+1}}\int_{z_n}^{z_{n+1}}\left(\frac{1}{Z(x_l,y,z)}|F^{x}(x_l, y, z,t)|^2 + \frac{1}{Z(x_{l+1}, y,z)}|G^{x}(x_{l+1},y, z, t)|^2 \right)dydz  \\
&-\sum_{l=1}^{L}\sum_{m=1}^{M}\sum_{n=1}^{N}\int_{x_l}^{x_{l+1}}\int_{z_n}^{z_{n+1}}\left(\frac{1}{Z(x, y_m, z)}|F^{y}(x, y_m, z, t)|^2 + \frac{1}{Z(x, y_{m+1}, z)}|G^{y}(x,y_{m+1}, z, t)|^2 \right)dxdz   \\
&-\sum_{l=1}^{L}\sum_{m=1}^{M}\sum_{n=1}^{N}\int_{x_l}^{x_{l+1}}\int_{y_m}^{y_{m+1}}\left(\frac{1}{Z(x, y, z_n)}|F^{y}(x, y, z_n, t)|^2 + \frac{1}{Z(x, y, z_{n+1})}|G^{y}(x,y, z_{n+1}, t)|^2 \right)dxdy   \\
& - \sum_{m = 1}^{M}\sum_{n=1}^{N}\int_{y_m}^{y_{m+1}}\int_{z_n}^{z_{n+1}}\mathbf{BT}^{(x)}dydz - \sum_{l = 1}^{L}\sum_{n=1}^{N}\int_{x_l}^{x_{l+1}}\int_{z_m}^{z_{m+1}}\mathbf{BT}^{(y)}dxdz - \sum_{l = 1}^{L}\sum_{m=1}^{M}\int_{x_l}^{x_{l+1}}\int_{y_m}^{y_{m+1}}\mathbf{BT}^{(z)}dydz,
%%%%%%%%%%%%%%%%%%%%%%%
\end{split}
\end{equation}
%%%%%%%%%%%%%%%%%%%%%%%
where
\[
\mathbf{BT}^{(x)} = \frac{1-|r_x(-x_b)|^2}{Z(-x_b,y, z)}|\chi^{(-x)}(-x_b, y, z, t)|^2 +\frac{1-|r_x(x_b)|^2}{Z(x_b, y, z)}|\chi^{(+x)}(x_b, y, z, t) |^2 \ge 0,
\]
\[
\mathbf{BT}^{(y)} = \frac{1-|r_y(-y_b)|^2}{Z(x,-y_b, z)}|\chi^{(-y)}(x, -y_b, z, t)|^2 +\frac{1-|r_y(y_b)|^2}{Z(x, y_b, z)}|\chi^{(+y)}(x, y_b, z, t) |^2 \ge 0,
\]
\[
\mathbf{BT}^{(z)} = \frac{1-|r_z(-z_b)|^2}{Z(x,y, -z_b)}|\chi^{(-z)}(x, y, -z_b, t)|^2 +\frac{1-|r_z(z_b)|^2}{Z(x, y, z_b)}|\chi^{(+z)}(x, y, z_b, t) |^2 \ge 0,
\]
%%%%%%%%%%%%%%%%%%%%%%%
with $ \chi^{\pm i} $,  $i = x, y, z$,  being the characteristic variables defined in \eqref{eq:characteristics}.
%%%%%%%%%%%%%%%%%%%%%%%
\end{theorem}
%%%%
Note that since we have not introduce any approximation yet the flux fluctuations vanish identically, ${F}^{i} \equiv 0$, ${G}^{i} \equiv 0$. Thus, \eqref{eq:weak_energy_elements} in  Theorem \ref{theorem:main_1}  is identical to \eqref{eq:product_4}, \eqref{energy_conservation}. When numerical approximations are introduced numerical fluctuations may not vanish, it will then contribute to artificial dissipation which will vanish in the limit of mesh resolution.
%%%%%%%%%%%%%%%%%%%%%%%
%%%%%%%%%%%%%%%%%%%%%%%

%%%%%%%%%%%%%%%%%%%%%%%
%%%%%%%%%%%%%%%%%%%%%%%
It is particularly noteworthy that Theorem \ref{theorem:main_1}  is not applicable when the PML is active, that is $d_x \ge 0$, $d_y \ge 0$, $d_z \ge 0$. We will need to take Laplace transform of \eqref{eq:weak_1_eq_1}--\eqref{eq:weak_1_eq_6} in time, and derive a result analogous to Theorem \ref{Theo:Stability_PML_Laplace}. We can prove
%%%%%%%%%%%%%%%%%%%
\begin{theorem}\label{Theo:Stability_PML_Laplace_Mutiple_Elements}
%%%%%%%%% %%%%%%%%%
Consider the   weak form of the PML equation  in the Laplace space \eqref{eq:acoustic_pml_3D_Laplace1}  subject to boundary conditions \eqref{eq:boundary_condition_acoustic_laplace} with $\Re{s}  = a \ge 0$ and constant damping  $d_x  \ge 0, d_y \ge 0,  d_z\ge 0$. Let 
%%%%%%%%% %%%%%%%%%
%%%%%%%%% %%%%%%%%%
\[ 
\widetilde{\mathbf{U}}\left(s\right)  = \left(s\widetilde{p}, \frac{1}{S_x} \frac{\partial \widetilde{p} }{\partial x}, \frac{1}{S_y} \frac{\partial \widetilde{p} }{\partial y}, \frac{1}{S_z} \frac{\partial \widetilde{p} }{\partial z}\right)^T, \quad 
%\]
%%%%%%%%% %%%%%%%%%
%\[ 
\widetilde{\mathbf{F}}\left(s\right)  = \left(\widetilde{F}_p(x,y,z), \frac{\rho}{S_x}   f_u(x,y,z),  \frac{\rho}{S_y}   f_v(x,y,z), \frac{\rho}{S_z}   f_w(x,y,z)\right)^T,
\]
with
\[
\widetilde{F}_p(x,y,z) = \frac{1}{S_x}f_{p}(x,y,z) - \frac{\kappa}{sS_yS_x}f_{\sigma}(x,y,z)  -  \frac{\kappa}{sS_zS_x}f_{\psi}(x,y,z).
\]
%%%%%%%%% %%%%%%%%%
We have
\begin{equation}\label{eq:energy_estimate_pml_laplace_corner_dg}
\widetilde{E}^2\left( \widetilde{\mathbf{U}}\left(s\right)  \right) + \mathrm{BT}(s)   \le \widetilde{E}\left( \widetilde{\mathbf{U}}\left(s\right) \right) \widetilde{E}\left(\widetilde{\mathbf{F}}\left(s\right) \right),
\end{equation}
%%%%%%%%% %%%%%%%%%
with
{
\small
\[
\mathrm{BT}(s)  = {\Re\left(\frac{1}{S_x}\right)\widetilde{BT}^{(x)} + \Re\left(\frac{1}{S_y}\right)\widetilde{BT}^{(y)}  + \Re\left(\frac{1}{S_z}\right)\widetilde{BT}^{(z)} } +  {\Re\left(\frac{1}{S_x}\right)\widetilde{IT}^{(x)} + \Re\left(\frac{1}{S_y}\right)\widetilde{IT}^{(y)}  + \Re\left(\frac{1}{S_z}\right)\widetilde{IT}^{(z)} }  \ge 0,
\]
}
where 
{
\small
%%%%%%%%%%%%%%%%%%%
\[
\widetilde{BT}^{(x)}  = |s|\sum_{m=1}^{M}\sum_{n=1}^{N}\int_{y_m}^{y_{m+1}}\int_{z_n}^{z_{n+1}}\left(\frac{1-|r_x(-x_b)|^2}{Z(-x_b,y, z)}|\widetilde{\chi}^{(-x)}(-x_b, y, z, s)|^2 +\frac{1-|r_x(x_b)|^2}{Z(x_b, y, z)}|\widetilde{\chi}^{(+x)}(x_b, y, z, s) |^2\right)dydz  ,
\]
%%%%%%%%%%%%%%%%%%%
%%%%%%%%%%%%%%%%%%%
\[
\widetilde{BT}^{(y)}  = 
|s| \sum_{l=1}^{L}\sum_{n=1}^{N}\int_{x_l}^{x_{l+1}}\int_{z_n}^{z_{n+1}}\left(\frac{1-|r_y(-y_b)|^2}{Z(x,-y_b, z)}|\widetilde{\chi}^{(-y)}(x, -y_b, z, s)|^2 +\frac{1-|r_y(y_b)|^2}{Z(x, y_b, z)}|\widetilde{\chi}^{(+y)}(x, y_b, z, s) |^2 \right)dxdz ,  
\]
%%%%%%%%%%%%%%%%%%%
%%%%%%%%%%%%%%%%%%%
\[
\widetilde{BT}^{(z)}  =
|s|\sum_{l=1}^{L}\sum_{m=1}^{M}\int_{x_l}^{x_{l+1}}\int_{y_m}^{y_{m+1}}\left(\frac{1-|r_z(-z_b)|^2}{Z(x,y, -z_b)}|\widetilde{\chi}^{(-z)}(x, y, -z_b, s)|^2 +\frac{1-|r_z(z_b)|^2}{Z(x, y, z_b)}|\widetilde{\chi}^{(+z)}(x, y, z_b, s) |^2\right)dxdy   ,
\]

%%%%%%%%%%%%%%%%%%%
\[
\widetilde{IT}^{(x)}  = |s|\sum_{l=1}^{L}\sum_{m=1}^{M}\sum_{n=1}^{N}\int_{y_m}^{y_{m+1}}\int_{z_n}^{z_{n+1}}\left(\frac{1}{Z(x_l,y,z)}|\widetilde{F}^{x}(x_l, y, z,s)|^2 + \frac{1}{Z(x_{l+1}, y,z)}|\widetilde{G}^{x}(x_{l+1},y, z, s)|^2 \right)dydz  ,
\]
%%%%%%%%%%%%%%%%%%%
%%%%%%%%%%%%%%%%%%%
\[
\widetilde{IT}^{(y)}  = 
|s| \sum_{l=1}^{L}\sum_{m=1}^{M}\sum_{n=1}^{N}\int_{x_l}^{x_{l+1}}\int_{z_n}^{z_{n+1}}\left(\frac{1}{Z(x, y_m, z)}|\widetilde{F}^{y}(x, y_m, z, s)|^2 + \frac{1}{Z(x, y_{m+1}, z)}|\widetilde{G}^{y}(x,y_{m+1}, z, s)|^2 \right)dxdz   ,
\]
%%%%%%%%%%%%%%%%%%%
%%%%%%%%%%%%%%%%%%%
\[
\widetilde{IT}^{(z)}  =
|s|\sum_{l=1}^{L}\sum_{m=1}^{M}\sum_{n=1}^{N}\int_{x_l}^{x_{l+1}}\int_{y_m}^{y_{m+1}}\left(\frac{1}{Z(x, y, z_n)}|\widetilde{F}^{y}(x, y, z_n, s)|^2 + \frac{1}{Z(x, y, z_{n+1})}|\widetilde{G}^{y}(x,y, z_{n+1}, s)|^2 \right)dxdy   .
\]
%%%%%%%%%%%%%%%%%%%
%%%%%%%%%%%%%%%%%%%
%& - \sum_{m = 1}^{M}\sum_{n=1}^{N}\int_{y_m}^{y_{m+1}}\int_{z_n}^{z_{n+1}}\mathbf{BT}^{(x)}dydz - \sum_{l = 1}^{L}\sum_{n=1}^{N}\int_{x_l}^{x_{l+1}}\int_{z_m}^{z_{m+1}}\mathbf{BT}^{(y)}dxdz - \sum_{l = 1}^{L}\sum_{m=1}^{M}\int_{x_l}^{x_{l+1}}\int_{y_m}^{y_{m+1}}\mathbf{BT}^{(z)}dydz
}
\end{theorem}
%%%%%%%%% %%%%%%%%%
Note that since we have not introduce any approximation yet the interface terms in   \eqref{eq:energy_estimate_pml_laplace_corner_dg} vanish completely, $\widetilde{IT}^{(i)} \equiv 0$, $i = x, y, z$. The energy estimate \eqref{eq:energy_estimate_pml_laplace_corner_dg} in Theorem \ref{Theo:Stability_PML_Laplace_Mutiple_Elements}   is identical to the energy estimate \eqref{eq:energy_estimate_pml_laplace_corner_cont} in Theorem \ref{theorem:main_1}. As we will see below, when numerical approximation is introduced we will see that Theorem \ref{Theo:Stability_PML_Laplace_Mutiple_Elements} can be extended to the DGSEM approximation of the PML.
%%%%%%%%%%%%%%%%%%

\subsection{The Galerkin approximation}
%%%%%%%%%%%%%%%%%%%%%%%
%%%%%%%%%%%%%%%%%%%%%%%
To begin, we map the element $\Omega_{lmn} = [x_l, x_{l+1}]\times [y_m, y_{m+1}]\times [z_n, z_{n+1}]$ to a reference element $(\xi, \eta, \theta) \in [-1, 1]^{3}$ by the linear transformation 
%%%%%%%%%%%%%%%%%%%%%%%
{
\small
\begin{align}\label{eq:transf}
x = x_l + \frac{\Delta{x}_l}{2}\left(1 + \xi \right),  \quad
y = y_m + \frac{\Delta{y}_m}{2}\left(1 + \eta \right), \quad 
z = z_n + \frac{\Delta{z}_n}{2}\left(1 + \theta \right),  
\end{align}
}
with
{
\small
\[
\quad \Delta{x}_l = x_{l+1} - x_l, \quad  \Delta{y}_m = y_{m+1} - y_m, \quad  \Delta{z}_n = y_{n+1} - y_n.
\]
}
%%%%%%%%%%%%%%%%%%%%%%%
Applying the transformation \eqref{eq:transf} to the weak problem \eqref{eq:weak_1_eq_1}--\eqref{eq:weak_1_eq_6}  yields
%%%%%%%%%%%%%%%%%%%%%%%
%%%%%%%%%%%%%%%%%%%%%%%
{
\small
\begin{align}
\label{eq:pde_1} 
&\int_{-1}^{1}\int_{-1}^{1}\int_{-1}^{1}  {\phi}_p\left(\xi, \eta, \theta\right) \left(\frac{1}{\kappa} \left(\frac{\partial p}{\partial t} + d_x p \right)  + \frac{2}{\Delta{x}_l}\frac{\partial u}{\partial \xi} +  \frac{2}{\Delta{y}_m}\frac{\partial v}{\partial \eta} + \frac{2}{\Delta{z}_n}\frac{\partial {w}}{\partial \theta}  + \sigma + \psi  \right)d\xi d\eta d\theta = \\ \nonumber
& - \frac{2}{\Delta{x}_l}\int_{-1}^{1}\int_{-1}^{1}{\left(\frac{\phi_p(-1, \eta, \theta)}{Z(-1, \eta, \theta)} F^{x}(-1, \eta, \theta, t) -  \frac{\phi_p(1, \eta, \theta)}{Z_s(1, \eta, \theta)}G^{x}(1, \eta, \theta, t)\right)}d\eta d\theta \\ \nonumber
%%%%%%
& - \frac{2}{\Delta{y}_m}\int_{-1}^{1}\int_{-1}^{1}{\left( \frac{\phi_p(\xi, -1, \theta)}{Z(\xi, -1, \theta)}F^{y}(\xi, -1, \theta, t) -  \frac{\phi_p(\xi, 1, \theta)}{Z(\xi, 1, \theta)}G^{y}(\xi, 1, \theta, t)\right)}d\xi  d\theta \\ \nonumber
& - \frac{2}{\Delta{z}_n}\int_{-1}^{1}\int_{-1}^{1}{\left( \frac{\phi_p(\xi, \eta, -1)}{Z(\xi, \eta, -1)}F^{z}(\xi, \eta, -1, t) -  \frac{\phi_p(\xi, \eta, 1)}{Z(\xi, \eta, 1)}G^{y}(\xi, \eta, 1, t)\right)}d\xi d\eta,  \ \\
%%%%%
%%%%%
\label{eq:pde_2} 
&\int_{-1}^{1}\int_{-1}^{1}\int_{-1}^{1}  {\phi}_u\left(\xi, \eta, \theta\right)  \left(\rho\left(\frac{\partial u}{\partial t} + d_x u \right) \frac{2}{\Delta{x}_k}\frac{\partial p}{\partial \xi} \right) d\xi d\eta d\theta =  \\ \nonumber
%%%%%
& - \frac{2}{\Delta{x}_l}\int_{-1}^{1}\int_{-1}^{1}{\left({\phi_u(-1, \eta, \theta)} F^{x}(-1, \eta, \theta, t) + {\phi_u(1, \eta, \theta)}G^{x}(1, \eta, \theta, t)\right)}d\eta d\theta,  \\ 
%%%%%
%%%%%
\label{eq:pde_3}
&\int_{-1}^{1}\int_{-1}^{1}\int_{-1}^{1}  {\phi}_v\left(\xi, \eta, \theta\right) \left(\rho \left(\frac{\partial v}{\partial t} + d_y v \right)  +  \frac{2}{\Delta{y}_l}\frac{\partial p}{\partial \eta} \right)d\xi d\eta d\theta = \\ \nonumber
%%%%%
& - \frac{2}{\Delta{y}_m}\int_{-1}^{1}\int_{-1}^{1}{\left({\phi_v(\xi, -1, \theta)}F^{y}(\xi, -1, \theta, t) +  {\phi_v(\xi, 1, \theta)}G^{y}(\xi, 1, \theta, t)\right)}d\xi  d\theta, \\
%%%%%
\label{eq:pde_4} 
%%%%%
&\int_{-1}^{1}\int_{-1}^{1}\int_{-1}^{1}  {\phi}_w\left(\xi, \eta, \theta\right) \left(\rho \left(\frac{\partial w}{\partial t} + d_z w \right)  + \frac{2}{\Delta{z}_m}\frac{\partial {p}}{\partial \theta} \right)d\xi d\eta d\theta = \\ \nonumber
%%%%%
& - \frac{2}{\Delta{z}_n}\int_{-1}^{1}\int_{-1}^{1}{\left({{\phi}_w\left(\xi, \eta, -1\right)}F^{z}(\xi, \eta, -1, t) +  {\phi_w(\xi, \eta, 1)}G^{y}(\xi, \eta, 1, t)\right)}d\xi d\eta,  
%%%%%%
%%%%%%
\end{align}
%%%%%

%%%%%
\begin{align}
%%%%%%
 \label{eq:aux_1} 
%%%%%
&\int_{-1}^{1}\int_{-1}^{1}\int_{-1}^{1}  {\phi}_{\sigma}\left(\xi, \eta, \theta\right) \left( \left(\frac{\partial \sigma }{\partial t} + d_y \sigma \right) +  \frac{2}{\Delta{y}_l} \left(d_y - d_x\right)\frac{\partial v}{\partial \eta} \right) d\xi d\eta d\theta = \\  \nonumber
%%%%%
& - \underbrace{\frac{2}{\Delta{y}_m}\int_{-1}^{1}\int_{-1}^{1}{ { \omega_y} \left(d_y - d_x\right) \left( \frac{\phi_{\sigma}(\xi, -1, \theta)}{Z(\xi, -1, \theta)}F^{y}(\xi, -1, \theta, t) -  \frac{\phi_{\sigma}(\xi, 1, \theta)}{Z(\xi, 1, \theta)}G^{y}(\xi, 1, \theta, t)\right)}d\xi  d\theta}_{\text{PML stabilizing flux fluctuation = 0}},
\\ 
%%%%%
 \label{eq:aux_2}
 %%%%%
&\int_{-1}^{1}\int_{-1}^{1}\int_{-1}^{1}  {\phi}_{\psi}\left(\xi, \eta, \theta\right) \left( \left(\frac{\partial \psi }{\partial t} + d_z \psi \right)   + \frac{2}{\Delta{z}_m} \left(d_z - d_x\right)\frac{\partial {w}}{\partial \theta} \right)d\xi d\eta d\theta = \\ \nonumber
& - \underbrace{\frac{2}{\Delta{z}_n}\int_{-1}^{1}\int_{-1}^{1}{{ \omega_z} \left(d_z - d_x\right)\left( \frac{\phi_{\psi}(\xi, \eta, -1)}{Z(\xi, \eta, -1)}F^{z}(\xi, \eta, -1, t) -  \frac{\phi_{\psi}(\xi, \eta, 1)}{Z(\xi, \eta, 1)}G^{y}(\xi, \eta, 1, t)\right)}d\xi d\eta}_{\text{PML stabilizing flux fluctuation = 0}} .
%%%%%
\end{align}
%%%
}
Inside the transformed  element  $(\xi, \eta, \theta) \in [-1, 1]^3$, approximate the elemental solution by a polynomial interpolant ${u}^{lmn}(\xi,\eta,\theta, t)$,  and write 
\begin{equation}\label{eq:variables_elemental}
{u}^{lmn}(\xi,\eta,\theta, t) =  \sum_{i = 1}^{P+1} \sum_{j = 1}^{P+1}\sum_{k = 1}^{P+1}{u}_{ijk}^{lmn}(t) \phi_{ijk}(\xi, \eta,\theta),
\end{equation}
%%%%%%%%%%%%%%%%%%%
where ${u}_{ijk}^{lmn}(t)$, are the elemental degrees of freedom to be determined, 
 and $ \phi_{ijk}(\xi, \eta,\theta)$ are the $ijk$-th interpolating polynomial. We consider tensor products of  nodal basis with $ \phi_{ijk}(\xi, \eta,\theta) = \mathcal{L}_i(\xi)\mathcal{L}_j(\eta)\mathcal{L}_k(\theta)$,  where $\mathcal{L}_i(\xi)$, $\mathcal{L}_j(\eta)$, $\mathcal{L}_k(\theta)$, are one dimensional nodal interpolating Lagrange polynomials of degree $P$, with $ \mathcal{L}_i(\xi_m) = \delta_{im}$. 
The interpolating nodes $\xi_m$, $m = 1, 2, \dots, P+1$, are the nodes of a Gauss quadrature with
%%%%%%%%%%%%%%%%%
%%%%%%%%%%%%%%%%%
\begin{equation}\label{eq:quad_rule_3D}
 \sum_{i = 1}^{P+1} \sum_{j = 1}^{P+1}  \sum_{k = 1}^{P+1}f(\xi_i, \eta_j, \theta_k)h_ih_jh_k \approx \int_{-1}^{1}\int_{-1}^{1}\int_{-1}^{1}f(\xi, \eta,\theta) d\xi d\eta d\theta,
\end{equation}
%%%%%%%%%%%%%%%%%
where $h_i > 0$, $h_j>0$, $h_k>0$, are the quadrature weights.
We will only use quadrature rules such  that for all polynomial integrand $f(\xi)$ of degree $\le 2P-1$, the corresponding one dimensional rule is exact, 
%\begin{equation}
 $\sum_{m = 1}^{P+1} f(\xi_m)h_m = \int_{-1}^{1}f(\xi) d\xi.$
%\end{equation}
  Admissible  candidates can be Gauss-Legendre-Lobatto quadrature rule with GLL nodes, Gauss-Legendre quadrature rule with GL nodes and Gauss-Legendre-Radau quadrature rule with GLR nodes. While both endpoints, $\xi = -1, 1$, are part of   GLL quadrature nodes,  the   GLR quadrature contains only the first endpoint $\xi = -1$ as a node. Lastly,  for the GL quadrature, both endpoints, $\xi = -1, 1$, are not quadrature nodes. Note that when an endpoint is not a quadrature node, $\xi_1 \ne -1$ or $\xi_{P+1} \ne 1$, extrapolation is needed to compute numerical fluxes at the element boundary, $\xi = -1$ or $\xi = 1$. We also remark that the GLL quadrature rule is exact for polynomial integrand of degree $2P-1$, GLR quadrature rule is exact for polynomial integrand of degree $2P$, and GL quadrature rule is exact for polynomial integrand of  degree $2P+1$.

 %%%%%%%%%%%%%%%%%%%
 
 We now make a classical  Galerkin approximation by choosing test functions  in the same space as the basis functions, so that the residual is orthogonal to the space of test functions. By rearranging the elemental degrees of freedom $[{u}_{ijk}^{lmn}(t) ]$ row-wise as a vector, $\mathbf{u}^{lmn}(t) $,  of length $(P+1)^d$ where $d = 3$ is the number of space dimensions, we have the evolution equation

{ 
 \small
\begin{equation}\label{eq:disc_elemental_pde1_pml_semi_disc}
\begin{split}
 &\left({\boldsymbol{\kappa}^{lmn}}\right)^{-1} \left(\frac{d \boldsymbol{p}^{lmn}( t)}{ d t} + \mathbf{d}_x \boldsymbol{p}^{lmn}( t) \right) + \mathbf{D}_x \boldsymbol{u}^{lmn}( t)  +  \mathbf{D}_y \boldsymbol{v}^{lmn}( t)  +  \mathbf{D}_y \boldsymbol{w}^{lmn}( t) +  \boldsymbol{\sigma}^{lmn}( t) + \boldsymbol{\psi}^{lmn}( t)  = \\
 %%%%%%
- &\mathbf{H}_x^{-1}\left(\frac{\mathbf{e}_x(-1)}{\boldsymbol{Z}}\mathbf{F}^{xlmn}(-1, \boldsymbol{\eta}, \boldsymbol{\theta}, t)- \frac{\mathbf{e}_x(1)} {\boldsymbol{Z}} \mathbf{G}^{xlmn}(1, \boldsymbol{\eta}, \boldsymbol{\theta}, t)  \right)
 %%%%%
 -\mathbf{H}_y^{-1}\left(\frac{\mathbf{e}_y(-1)}{\boldsymbol{Z}}\mathbf{F}^{ylmn}(\boldsymbol{\xi}, 1, \boldsymbol{\theta}, t)- \frac{\mathbf{e}_y(1)} {\boldsymbol{Z}} \mathbf{G}^{ylmn}(\boldsymbol{\xi}, 1, \boldsymbol{\theta}, t)  \right)\\
%%%%%
-&\mathbf{H}_z^{-1}\left(\frac{\mathbf{e}_z(-1)}{\boldsymbol{Z}}\mathbf{F}^{zlmn}(\boldsymbol{\xi}, \boldsymbol{\eta}, -1, t)- \frac{\mathbf{e}_z(1)} {\boldsymbol{Z}} \mathbf{G}^{zlmn}(\boldsymbol{\xi}, \boldsymbol{\eta}, 1, t)  \right),
\end{split}
\end{equation}
%%%%%%%
%%%%%%%
 \begin{equation}\label{eq:disc_elemental_pde2_pml_semi_disc}
\begin{split}
& {\boldsymbol{\rho}^{lmn}} \left(\frac{d \boldsymbol{u}^{lmn}( t)}{ d t} + \mathbf{d}_x \boldsymbol{u}^{lmn}( t) \right)  + \mathbf{D}_x \boldsymbol{p}^{lmn}( t)  = 
- \mathbf{H}_x^{-1}\left({\mathbf{e}_x(-1)}\mathbf{F}^{xlmn}(-1, \boldsymbol{\eta}, \boldsymbol{\theta}, t)+ {\mathbf{e}_x(1)}  \mathbf{G}^{xlmn}(1, \boldsymbol{\eta}, \boldsymbol{\theta}, t)  \right),
%&-\left(H^{-1}_x\otimes I\right) \left( \left({\Phi}(-1)\otimes I\right)\left(\mathbf{F}^{(x)}\right)^{kl}(-1, \eta, t) + \left({\Phi}(1)\otimes I \right)\left(\mathbf{G}^{(x)}\right)^{kl}(1, \eta, t)\right)\\
\end{split}
\end{equation}
%%%%%%
%%%%%%
 \begin{equation}\label{eq:disc_elemental_pde3_pml_semi_disc}
\begin{split}
& {\boldsymbol{\rho}^{lmn}} \left( \frac{d \boldsymbol{v}^{lmn}( t)}{ d t} + \mathbf{d}_y \boldsymbol{v}^{lmn}( t) \right)  + \mathbf{D}_y \boldsymbol{p}^{lmn}( t) =  
 -\mathbf{H}_y^{-1}\left({\mathbf{e}_y(-1)}\mathbf{F}^{ylmn}(\boldsymbol{\xi}, -1, \boldsymbol{\theta}, t)+ {\mathbf{e}_y(1)} \mathbf{G}^{ylmn}(\boldsymbol{\xi}, 1, \boldsymbol{\theta}, t)  \right),
\end{split}
\end{equation}
%%%%%%
%%%%%%
\begin{equation}\label{eq:disc_elemental_pde4_pml_semi_disc}
\begin{split}
& {\boldsymbol{\rho}^{lmn}} \left( \frac{d \boldsymbol{w}^{lmn}( t)}{ d t} + \mathbf{d}_z \boldsymbol{w}^{lmn}( t) \right)  + \mathbf{D}_z \boldsymbol{p}^{lmn}( t) =  
 -\mathbf{H}_z^{-1}\left({\mathbf{e}_z(-1)}\mathbf{F}^{zlmn}(\boldsymbol{\xi}, \boldsymbol{\eta}, -1, t) + {\mathbf{e}_z(1)}  \mathbf{G}^{zlmn}(\boldsymbol{\xi}, \boldsymbol{\eta}, 1, t)  \right),
\end{split}
\end{equation}
%%%%%%
%%%%%%
%%%%%%
\begin{equation}\label{eq:disc_elemental_pde5_pml_semi_disc}
\begin{split}
& \left( \frac{d \boldsymbol{\sigma}^{lmn}( t)}{ d t} + \mathbf{d}_y \boldsymbol{\sigma}^{lmn}( t) \right)  + { \left(\mathbf{d}_y - \mathbf{d}_x\right)}\mathbf{D}_y \boldsymbol{v}^{lmn}( t) =  
  -\underbrace{\omega_y\mathbf{H}_y^{-1} { \left(\mathbf{d}_y - \mathbf{d}_x\right)}\left( \frac{\mathbf{e}_y(-1)}{\boldsymbol{Z}}\mathbf{F}^{ylmn}(\boldsymbol{\xi}, -1, \boldsymbol{\theta}, t)- \frac{\mathbf{e}_y(1)} {\boldsymbol{Z}} \mathbf{G}^{ylmn}(\boldsymbol{\xi}, 1, \boldsymbol{\theta}, t)\right)}_{\text{PML stabilizing flux fluctuation} \to 0},
\end{split}
\end{equation}
%%%%%%%
%%%%%%%
\begin{equation}\label{eq:disc_elemental_pde6_pml_semi_disc}
\begin{split}
& \left( \frac{d \boldsymbol{\psi}^{lmn}( t)}{ d t} + \mathbf{d}_y \boldsymbol{\psi}^{lmn}( t) \right)  + { \left(\mathbf{d}_z - \mathbf{d}_x\right)}\mathbf{D}_z \boldsymbol{w}^{lmn}( t) =  
  -\underbrace{\omega_z\mathbf{H}_z^{-1} { \left(\mathbf{d}_z - \mathbf{d}_x\right)}\left( \frac{\mathbf{e}_z(-1)}{\boldsymbol{Z}}\mathbf{F}^{zlmn}(\boldsymbol{\xi}, \boldsymbol{\eta}, -1, t)- \frac{\mathbf{e}_z(1)} {\boldsymbol{Z}} \mathbf{G}^{zlmn}(\boldsymbol{\xi}, \boldsymbol{\eta}, 1, t)  \right)}_{\text{PML stabilizing flux fluctuation} \to 0}.
\end{split}
\end{equation}
}

%%%%%%%%%%%%%%%%%%
Here, the spatial operators are denoted
\[
\mathbf{D}_x = \frac{2}{\Delta{x}}\left(D\otimes I\otimes I\right), \quad \mathbf{D}_y = \frac{2}{\Delta{y}}\left(I\otimes D\otimes I\right), \quad \mathbf{D}_z = \frac{2}{\Delta{z}}\left(I\otimes I\otimes D\right),
\]
%%%%%%%%%%%%%%%%%%
%%%%%%%%%%%%%%%%%%
\[
\mathbf{H}_x = \frac{\Delta{x}}{2}\left(H\otimes I\otimes I\right), \quad \mathbf{H}_y = \frac{\Delta{y}}{2}\left(I\otimes H\otimes I\right), \quad \mathbf{H}_z = \frac{\Delta{z}}{2}\left(I\otimes I\otimes H\right),
\]
%%%%%%%%%%%%%%%%%%
\[
\mathbf{B}_x(\xi, \eta) = \left(\boldsymbol{\Phi}(\xi, \eta)\otimes I\otimes I\right), \quad \mathbf{B}_y(\xi, \eta) = \left(I\otimes \boldsymbol{\Phi}(\xi, \eta)\otimes I\right), \quad \mathbf{B}_z(\xi, \eta) = \left(I\otimes I\otimes \boldsymbol{\Phi}(\xi,\eta)\right),
\]
%%%%%%%%%%%%%%%%%%
\[
\mathbf{e}_x(\xi) = \left(\boldsymbol{e}(\xi)\otimes I\otimes I\right), \quad \mathbf{e}_y(\xi) = \left(I\otimes \boldsymbol{e}(\xi)\otimes I\right), \quad \mathbf{e}_z(\xi) = \left(I\otimes I\otimes \boldsymbol{e}(\xi)\right),
\]
%%%%%%%%%%%%%%%%%%
\[
\mathbf{H}_x\mathbf{H}_y = \frac{\Delta{x}}{2} \frac{\Delta{y}}{2} \left(H\otimes H\otimes I\right), \quad \mathbf{H}_x\mathbf{H}_z =  \frac{\Delta{x}}{2} \frac{\Delta{z}}{2} \left(I\otimes H\otimes H\right), \quad  \mathbf{H}_y\mathbf{H}_z =  \frac{\Delta{y}}{2} \frac{\Delta{z}}{2}\left(I\otimes H\otimes H\right),
\]
%%%%%%%%%%%%%%%%%%
%%%%%%%%%%%%%%%%%%
\[
 \mathbf{H} =  \mathbf{H}_x\mathbf{H}_y\mathbf{H}_z  = \frac{\Delta{x}}{2} \frac{\Delta{y}}{2} \frac{\Delta{z}}{2}\left(H\otimes H\otimes H\right).
\]
%%%%%%%%%%%%%%%%%%
%%%%%%%%%%%%%%%%%%
where
%%%%%%%%%%%%%%%%%%
%%%%%%%%%%%%%%%%%%
\[
\boldsymbol{\Phi}(\xi, \eta) = \boldsymbol{e}(\xi) \boldsymbol{e}^T(\eta), \quad  \boldsymbol{e}(\eta) = [\mathcal{L}_i(\eta), \mathcal{L}_i(\eta), \cdots, \mathcal{L}_{P+1}(\eta)]^T.
%\]
%%%%%%%%%%%%%%%%%
%%%%%%%%%%%%%%%%%
%\[
%\boldsymbol{e}_{N+1}^{(y)} = \left(H^N \otimes  {\Phi}^N(1)\right), \quad \boldsymbol{e}_{1}^{(y)} = \left(H^N \otimes  {\Phi}^N(-1)\right).
\]
The one dimensional  matrices $ H$, $Q$  are defined by
%%%%%%%%%%%%%%%%%
%%%%%%%%%%%%%%%%%
{\small
\begin{equation}
 H = \mathrm{diag}[h_1, h_2, \cdots, h_{N+1}], \quad Q_{ij} = \sum_{m = 1}^{N+1} h_m \mathcal{L}_i(\xi_m)  {\mathcal{L}_j^{\prime}(\xi_m)} = \int_{-1}^{1}\mathcal{L}_i(\xi)  {\mathcal{L}_j^{\prime}(\xi)} d\xi,  \quad {\Phi}_{ij}(\xi, \eta) =  \mathcal{L}_i(\xi)  \mathcal{L}_j(\eta).
\end{equation}
}
%%%%%%%%%%%%%%%%%
%%%%%%%%%%%%%%%%%
Note that the matrix
%%%%%%%%%%%%%%%%%
%%%%%%%%%%%%%%%%%
\begin{equation}
D = H^{-1} Q \approx \frac{\partial}{\partial \xi},
\end{equation}
%%%%%%%%%%%%%%%%%
is a one space dimensional spectral difference approximation of the first derivative.
%%%%%%%%%%%%%%%%%

Using the fact that the quadrature rule is exact  for all polynomial  integrand of degree $\le 2P-1$ 
%%%%%%%%%
  implies that
%%%%%%%%%
\begin{equation}\label{eq:sbp_property_a}
Q + Q^T =  B,
\quad
B = \boldsymbol{\Phi}(1, 1)-\boldsymbol{\Phi}(-1, -1) =  \boldsymbol{e}(1) \boldsymbol{e}^T(1) - \boldsymbol{e}(-1) \boldsymbol{e}^T(-1).
\end{equation}
Equation \eqref{eq:sbp_property_a} is the discrete equivalence of the integration-by-parts property.
If boundary points $\xi = -1,1$ are quadrature nodes and we consider nodal bases  with $ \mathcal{L}_j(\xi_i) = \delta_{ij}$ then we have $B = \text{diag}[-1, 0,0, \dots, 0, 1].$ 
%%%%%%%%%%%%%%%%%
%In the finite difference literature \cite{DelReyFernandezBoomZingg2014, DuruandDunham2016} equation \eqref{eq:sbp_property_a} is analogous to the so-called  summation-by-parts (SBP) property.
%%%%%%%%%%%%%%%%%

We will now derive a semi-discrete energy equation analogous to \eqref{eq:product_4}, \eqref{energy_conservation} and Theorem \ref{theorem:main_1}.
Introduce the elemental semi-discrete energy density
{
\small
\[
dE^{lmn}_{ijk}=\frac{1}{2}\left[\frac{1}{\kappa^{lmn}(\xi_i, \eta_j, \theta_k)}|p^{lmn}(\xi_i, \eta_j, \theta_k, t)|^2 + \rho^{lmn}(\xi_i, \eta_j, \theta_k)\left(|u^{lmn}(\xi_i, \eta_j, \theta_k, t)|^2 + |v^{lmn}(\xi_i, \eta_j, \theta_k, t)|^2 + |w^{lmn}(\xi_i, \eta_j, \theta_k, t)|^2\right)\right]
\]
}
and the corresponding semi-discrete energy
{
\small
\begin{equation}\label{eq:physical_energy_semi_discrete}
\mathcal{E}(t) =   \sum_{l=1}\sum_{m=1}\sum_{n=1}\frac{\Delta{x}_l}{2}\frac{\Delta{y}_m}{2}\frac{\Delta{z}_n}{2}\sum_{i = 1}^{P+1} \sum_{j = 1}^{P+1}  \sum_{k = 1}^{P+1}dE^{lmn}_{ijk}(t) h_ih_jh_k.
\end{equation}
%%%%%%%%%%%%%%%%%%%
}
%%%%%%%%%%%%%%%%%%%
%%%%%%%%%%%%%%%%%%%
Note that the semi-discrete energy $\mathcal{E}(t)$ in \eqref{eq:physical_energy_semi_discrete} is obtained by replacing the volume integral in \eqref{eq:physical_energy} by the quadrature rule \eqref{eq:quad_rule_3D}.
%%%%%%%%%%%%%%%%%%%
%%%%%%%%%%%%%%%%%%%
\begin{theorem}\label{theorem:discrete_stability_no_pml}
Consider the semi-discrete approximation \eqref{eq:disc_elemental_pde1_pml_semi_disc}--\eqref{eq:disc_elemental_pde4_pml_semi_disc}. When the damping vanishes, $d_x = d_y = d_z = 0$, the solution of the semi-discrete approximation satisfies the energy identity
{
\small
\begin{equation}\label{eq:discrete_energy_no_pml}
\begin{split}
&\frac{d}{dt} {\mathcal{E}}(t) = -\sum_{l=1}\sum_{m=1}\sum_{n=1}\frac{\Delta{y}_m}{2}\frac{\Delta{z}_n}{2}\left(\sum_{j=1}^{P+1}\sum_{k=1}^{P+1}\left(\frac{1}{Z(-1,\eta_j, \theta_k)}|F^{xlmn}(-1, \eta_j, \theta_k,t)|^2 + \frac{1}{Z^{}(1, \eta_j, \theta_k)}|G^{xlmn}(1, \eta_j, \theta_k, t)|^2 \right){ h_j}{h_k}\right)  \\
&-\sum_{l=1}\sum_{m=1}\sum_{n=1}\frac{\Delta{x}_l}{2}\frac{\Delta{z}_n}{2}\left(\sum_{i=1}^{P+1}\sum_{k=1}^{P+1}\left(\frac{1}{Z(\xi_i, -1,\theta_k)}|F^{ylmn}(\xi_i, -1,\theta_k, t)|^2 + \frac{1}{Z^{}(\xi_i, 1,\theta_k)}|G^{ylmn}(\xi_i, 1,\theta_k, t)|^2 \right){ h_i}{h_k}\right)  \\
&-\sum_{l=1}\sum_{m=1}\sum_{n=1}\frac{\Delta{x}_l}{2}\frac{\Delta{y}_m}{2}\left(\sum_{i=1}^{P+1}\sum_{j=1}^{P+1}\left(\frac{1}{Z(\xi_i,\eta_j, -1)}|F^{zlmn}(\xi_i,\eta_j, -1, t)|^2 + \frac{1}{Z(\xi_i,\eta_j, 1)}|G^{zlmn}(\xi_i,\eta_j, 1, t)|^2 \right){ h_i}{h_j}\right)  \\
%%%%%%%%%%%%%%%%
& - \sum_{m = 1}\sum_{n=1}\frac{\Delta{y}_m}{2}\frac{\Delta{z}_n}{2}\left(\sum_{j=1}^{P+1}\sum_{k=1}^{P+1}\left(\frac{1-|r_x(-x_b)|^2}{Z(-1,\eta_j, \theta_k)}|\chi^{-x1mn}(-1,\eta_j, \theta_k, t)|^2 +\frac{1-|r_x(x_b)|^2}{Z(1,\eta_j, \theta_k)}|\chi^{+xLmn}(1,\eta_j, \theta_k, t) |^2\right){ h_j}{h_k}\right)\\
%%%%%%%%%%%%%%%%%%%%%%
%%%%%%%%%%%%%%%%
& - \sum_{l = 1}\sum_{n=1}\frac{\Delta{x}_l}{2}\frac{\Delta{z}_n}{2}\left(\sum_{i=1}^{P+1}\sum_{k=1}^{P+1}\left(\frac{1-|r_y(-y_b)|^2}{Z(\xi_i,-1, \theta_k)}|\chi^{-yl1n}(\xi_i,-1, \theta_k, t)|^2 +\frac{1-|r_y(y_b)|^2}{Z(\xi_i, 1, \theta_k)}|\chi^{+ylMn}(\xi_i, 1, \theta_k, t) |^2\right){ h_i}{h_k}\right)\\
%%%%%%%%%%%%%%%%%%%%%%
%%%%%%%%%%%%%%%%
& - \sum_{l = 1}\sum_{m=1}\frac{\Delta{x}_l}{2}\frac{\Delta{y}_m}{2}\left(\sum_{i=1}^{P+1}\sum_{j=1}^{P+1}\left(\frac{1-|r_z(-z_b)|^2}{Z(\xi_i,\eta_j, -1)}|\chi^{-zlm1}(\xi_i,\eta_j, -1, t)|^2 +\frac{1-|r_z(z_b)|^2}{Z(\xi_i,\eta_j, 1)}|\chi^{+zlmN}(\xi_i,\eta_j, 1, t) |^2\right){ h_i}{h_j}\right)\\
%%%%%%%%%%%%%%%%%%%%%%
%%%%%%%%%%%%%%%%%%%%%%%
\end{split}
\end{equation}
}
%%%%%%%%%%%%%%%%%%%%%%%
with  the elemental quantities $ \chi^{\pm i lmn} $,  $i = x, y, z$  being the characteristics defined in \eqref{eq:characteristics}.
%%%%%%%%%%%%%%%%%%%%%%%
\end{theorem}
%%%%
%%%%%%%%%%%%%%%%%%%%%%%
The proof of Theorem \ref{theorem:discrete_stability_no_pml} can be easily adapted from \cite{DuruGabrielIgel2017}.
Therefore we omit it here. When the damping vansishes,  $d_x = d_y = d_z = 0$,  by \eqref{eq:discrete_energy_no_pml} we know that the semi-discrete approximation is asymptotically stable. However, Theorem \ref{theorem:discrete_stability_no_pml} is not valid when the PML is present, $d_j > 0$ for any $j =x,y,z$.
%%%%%%%%%%%%%%%%%%%%%%%
%%%%%%%%%%%%%%%%%%%%%%%

%%%%%%%%%%%%%%%%%%%%%%%
%%%%%%%%%%%%%%%%%%%%%%%
Note that for the continuous PML it is  technically difficult to derive an  energy estimate for the PML in the  time domain. 
Working in the Laplace domain simplifies the analysis. If there is an energy estimate in the Laplace space, in principle, we can invert the Laplace transform obtaining an energy estimate in the time-domain.
%%%%%%%%%%%%%%%%%%%%%%%
%%%%%%%%%%%%%%%%%%%%%%%
%As before, we  can  derive energy estimates in the Laplace domain,  ensuring asymptotic stability, see theorems \ref{Theo:Stability_PML_Laplace} and \ref{Theo:Stability_PML_Laplace_Mutiple_Elements}.
%%%%%%%%%%%%%%%%%%%%%%%
%%%%%%%%%%%%%%%%%%%%%%%
As in the continuous setting, for the semi-discrete PML, we will perform Laplace transform in time and derive an energy equation analogous  to \eqref{eq:energy_estimate_pml_laplace_corner_cont} and \eqref{eq:energy_estimate_pml_laplace_corner_dg}. To simplify the analysis, we will consider a two element DGSEM approximation containing only one inter-element boundary in the $x$-direction, but includes all the six physical boundaries of the cuboidal computational domain. The analysis  can easily be extended to arbitrarily  many-element DGSEM approximation of the PML.
%%%%%%%%%%%%%%%%%%%%%%%
%%%%%%%%%%%%%%%%%%%%%%%

 Now, take the Laplace transform of  \eqref{eq:disc_elemental_pde1_pml_semi_disc}--\eqref{eq:disc_elemental_pde4_pml_semi_disc} in time,  we have
{ 
 \small
\begin{equation}\label{eq:disc_elemental_pde1_pml_laplace}
\begin{split}
 \left({\boldsymbol{\kappa}^{lmn}}\right)^{-1} s{\widetilde{\boldsymbol{p}}^{lmn}(s)}  &=  -\frac{1}{S_x}\left(\mathbf{D}_x{\widetilde{\boldsymbol{u}}^{lmn}(s)}  +\mathbf{H}_x^{-1}\left(\frac{\mathbf{e}_x(-1)}{\boldsymbol{Z}}{\widetilde{\mathbf{F}}}^{xlmn}(-1, \boldsymbol{\eta}, \boldsymbol{\theta}, s)- \frac{\mathbf{e}_x(1)} {\boldsymbol{Z}} {\widetilde{\mathbf{G}}}^{xlmn}(1, \boldsymbol{\eta}, \boldsymbol{\theta}, s)  \right)\right)\\
 %%%%%
& -\frac{1}{S_y}\left(\mathbf{D}_y{\widetilde{\boldsymbol{v}}^{lmn}(s)}  +\mathbf{H}_y^{-1}\left(\frac{\mathbf{e}_y(-1)}{\boldsymbol{Z}}{\widetilde{\mathbf{F}}}^{ylmn}(\boldsymbol{\xi}, 1, \boldsymbol{\theta}, s)- \frac{\mathbf{e}_y(1)} {\boldsymbol{Z}} {\widetilde{\mathbf{G}}}^{ylmn}(\boldsymbol{\xi}, 1, \boldsymbol{\theta}, s)  \right)\right)\\
%%%%%
&-\frac{1}{S_z}\left(\mathbf{D}_z{\widetilde{\boldsymbol{w}}^{lmn}(s)} + \mathbf{H}_z^{-1}\left(\frac{\mathbf{e}_z(-1)}{\boldsymbol{Z}}{\widetilde{\mathbf{F}}}^{zlmn}(\boldsymbol{\xi}, \boldsymbol{\eta}, -1, s)- \frac{\mathbf{e}_z(1)} {\boldsymbol{Z}} {\widetilde{\mathbf{G}}}^{zlmn}(\boldsymbol{\xi}, \boldsymbol{\eta}, 1, s)  \right)\right)\\
& -\frac{\left(1-\omega_y\right)}{S_y}\left(\mathbf{H}_y^{-1}\left(\mathbf{d}_x- \mathbf{d}_y\right)\left(\frac{\mathbf{e}_y(-1)}{\boldsymbol{Z}}{\widetilde{\mathbf{F}}}^{ylmn}(\boldsymbol{\xi}, 1, \boldsymbol{\theta}, s)- \frac{\mathbf{e}_y(1)} {\boldsymbol{Z}} {\widetilde{\mathbf{G}}}^{ylmn}(\boldsymbol{\xi}, 1, \boldsymbol{\theta}, s)  \right)\right)\\
%%%%%
&-\frac{\left(1-\omega_z\right)}{S_z}\left( \mathbf{H}_z^{-1}\left(\mathbf{d}_x- \mathbf{d}_z\right)\left(\frac{\mathbf{e}_z(-1)}{\boldsymbol{Z}}{\widetilde{\mathbf{F}}}^{zlmn}(\boldsymbol{\xi}, \boldsymbol{\eta}, -1, s)- \frac{\mathbf{e}_z(1)} {\boldsymbol{Z}} {\widetilde{\mathbf{G}}}^{zlmn}(\boldsymbol{\xi}, \boldsymbol{\eta}, 1, s)  \right)\right)
\end{split}
\end{equation}
}
%%%%%%%
{
\small
%%%%%%%
 \begin{equation}\label{eq:disc_elemental_pde2_pml_laplace}
\begin{split}
&  {\boldsymbol{\rho}^{lmn}}  s{\widetilde{\boldsymbol{u}}^{lmn}(s)}  = -  \frac{1}{S_z}\left( \mathbf{D}_x  \widetilde{\boldsymbol{p}}^{lmn}(s)   
+ \mathbf{H}_x^{-1}\left({\mathbf{e}_x(-1)}{\widetilde{\mathbf{F}}}^{xlmn}(-1, \boldsymbol{\eta}, \boldsymbol{\theta}, s)+ {\mathbf{e}_x(1)}  {\widetilde{\mathbf{G}}}^{xlmn}(1, \boldsymbol{\eta}, \boldsymbol{\theta}, s)  \right)\right)
%&-\left(H^{-1}_x\otimes I\right) \left( \left({\Phi}(-1)\otimes I\right)\left(\mathbf{F}^{(x)}\right)^{kl}(-1, \eta, t) + \left({\Phi}(1)\otimes I \right)\left(\mathbf{G}^{(x)}\right)^{kl}(1, \eta, t)\right)\\
\end{split}
\end{equation}
%%%%%%
%%%%%%
 \begin{equation}\label{eq:disc_elemental_pde3_pml_laplace}
\begin{split}
&   {\boldsymbol{\rho}^{lmn}}  s{\widetilde{\boldsymbol{v}}^{lmn}(s)}   = - \frac{1}{S_y}\left( \mathbf{D}_y \widetilde{\boldsymbol{p}}^{lmn}(s) 
 + \mathbf{H}_y^{-1}\left({\mathbf{e}_y(-1)}{\widetilde{\mathbf{F}}}^{ylmn}(\boldsymbol{\xi}, 1, \boldsymbol{\theta}, s)+ {\mathbf{e}_y(1)} {\widetilde{\mathbf{G}}}^{ylmn}(\boldsymbol{\xi}, 1, \boldsymbol{\theta}, s)  \right) \right)\\
\end{split}
\end{equation}
%%%%%%
%%%%%%
\begin{equation}\label{eq:disc_elemental_pde4_pml_laplace}
\begin{split}
&  {\boldsymbol{\rho}^{lmn}}  s{\widetilde{\boldsymbol{w}}^{lmn}(s)} =-\frac{1}{S_z}\left( \mathbf{D}_z {\widetilde{\boldsymbol{w}}}^{lmn}(s) 
 + \mathbf{H}_z^{-1}\left({\mathbf{e}_z(-1)}{\widetilde{\mathbf{F}}}^{zlmn}(\boldsymbol{\xi}, \boldsymbol{\eta}, -1, s) + {\mathbf{e}_z(1)}  {\widetilde{\mathbf{G}}}^{zlmn}(\boldsymbol{\xi}, \boldsymbol{\eta}, 1, s)  \right)\right)\\
\end{split}
\end{equation}
%%%%%%
%%%%%%
%%%%%%
}
%%%%%%
%%%%%%
Note that, with  $\omega_y = 1, \omega_z = 1$ the last two terms in the right hand side of \eqref{eq:disc_elemental_pde1_pml_laplace} vanish identically.

Consider the DGSEM approximation \eqref{eq:disc_elemental_pde1_pml_laplace}--\eqref{eq:disc_elemental_pde4_pml_laplace} in two elements, separated at $x = x_2$. Note that $L = 2$, $M = 1$, $N = 1$, so there is only one internal element boundary. Let us denote the solutions in the first element by $\widetilde{p}^{1}(\xi,\eta,\theta,s)$, $\widetilde{u}^{1}(\xi,\eta,\theta, s)$, $\widetilde{v}^{1}(\xi,\eta,\theta, s)$, $\widetilde{w}^{1}(\xi,\eta,\theta, s)$ and the solution the second element denoted by $\widetilde{p}^{2}(\xi,\eta,\theta, s)$, $\widetilde{u}^{2}(\xi,\eta,\theta, s)$, $\widetilde{v}^{2}(\xi,\eta,\theta, s)$, $\widetilde{w}^{2}(\xi,\eta,\theta, s)$.
%%%%%%%%%%%%%%%%%%%%%%%
%%%%%%%%%%%%%%%%%%%%%%%
 Introduce the matrices
%%%%%%%%%%%%%%%%%%%%%%%
%%%%%%%%%%%%%%%%%%%%%%%
{
\small
\begin{align*}
 &\widetilde{\mathbf{D}}_x =  \frac{1}{S_x} \left( \begin{pmatrix} 
\mathbf{D}_z & \mathbf{0} \\
 \mathbf{0}  & \mathbf{D}_x\\
 \end{pmatrix}
 + 
%%%%%
\begin{pmatrix} 
\mathbf{H}_x^{-1} & \mathbf{0} \\
 \mathbf{0}  & \mathbf{H}_x^{-1}\\
 \end{pmatrix} 
 %%%%%%%
 \left(
 \frac{1+r_x}{2}
 \begin{pmatrix} 
\mathbf{B}_x{(-1,-1)}  & \mathbf{0} \\
 \mathbf{0}  & -\mathbf{B}_x{(1,1)}\\
 \end{pmatrix} 
 +
 \frac{1}{2}
%%%%%
 %%%%%
 \begin{pmatrix} 
-\mathbf{B}_x{(1,1)} &  \mathbf{B}_x{(1,-1)} \\
-\mathbf{B}_x^T{(1,-1)}   & \mathbf{B}_x{(-1,-1)} \\
 \end{pmatrix} 
 \right)
 \right),
 \\
  &\widetilde{\mathbf{D}}_y =  \frac{1}{S_y} \left( \begin{pmatrix} 
\mathbf{D}_y & \mathbf{0} \\
 \mathbf{0}  & \mathbf{D}_y\\
 \end{pmatrix}
 + 
%%%%%
\begin{pmatrix} 
\mathbf{H}_y^{-1} & \mathbf{0} \\
 \mathbf{0}  & \mathbf{H}_y^{-1}\\
 \end{pmatrix} 
 %%%%%%%
 \left(
 \frac{1+r_y}{2}
 \begin{pmatrix} 
\mathbf{B}_y{(-1,-1)}  & \mathbf{0} \\
 \mathbf{0}  & -\mathbf{B}_y{(1,1)}\\
 \end{pmatrix} 
 \right)
 \right),
 \\
  &\widetilde{\mathbf{D}}_z =  \frac{1}{S_z} \left( \begin{pmatrix} 
\mathbf{D}_z & \mathbf{0} \\
 \mathbf{0}  & \mathbf{D}_z\\
 \end{pmatrix}
 + 
%%%%%
\begin{pmatrix} 
\mathbf{H}_z^{-1} & \mathbf{0} \\
 \mathbf{0}  & \mathbf{H}_z^{-1}\\
 \end{pmatrix} 
 %%%%%%%
 \left(
 \frac{1+r_z}{2}
 \begin{pmatrix} 
\mathbf{B}_z{(-1,-1)}  & \mathbf{0} \\
 \mathbf{0}  & -\mathbf{B}_z{(1,1)}\\
 \end{pmatrix} 
 \right)
 \right),
\end{align*}
}
%%%%%
%%%%%
{
\small
\begin{align*}
 &\widetilde{\mathbf{H}}_x(s, d_x)  =   \begin{pmatrix} 
\mathbf{H} & \mathbf{0} \\
 \mathbf{0}  & \mathbf{H}\\
 \end{pmatrix}
\left( \mathbf{I} + 
%%%%%
\begin{pmatrix} 
\mathbf{H}_x^{-1} & \mathbf{0} \\
 \mathbf{0}  & \mathbf{H}_x^{-1}\\
 \end{pmatrix} 
 %%%%%%%
 \left(
 \frac{(1- r_x)c}{2 sS_x}
 \begin{pmatrix} 
\mathbf{B}_x{(-1,-1)}  & \mathbf{0} \\
 \mathbf{0}  & \mathbf{B}_x{(1,1)} \\
 \end{pmatrix} 
 +
 \frac{c}{2sS_x}
%%%%%
 %%%%%
 \begin{pmatrix} 
\mathbf{B}_x{(1,1)}  & - \mathbf{B}_x{(1,-1)} \\
-\mathbf{B}_x^T{(1,-1)}   & \mathbf{B}_x{(-1,-1)} \\
 \end{pmatrix} 
 \right)
 \right)^{-1},
 \\
 &\widetilde{\mathbf{H}}_y(s, d_y)  =   \begin{pmatrix} 
\mathbf{H} & \mathbf{0} \\
 \mathbf{0}  & \mathbf{H}\\
 \end{pmatrix}
\left( \mathbf{I} + 
%%%%%
\begin{pmatrix} 
\mathbf{H}_y^{-1} & \mathbf{0} \\
 \mathbf{0}  & \mathbf{H}_y^{-1}\\
 \end{pmatrix} 
 %%%%%%%
 \left(
 \frac{(1- r_y)c}{2 sS_y}
 \begin{pmatrix} 
\mathbf{B}_y{(-1,-1)}  & \mathbf{0} \\
 \mathbf{0}  & \mathbf{B}_y{(1,1)} \\
 \end{pmatrix} 
 \right)
 \right)^{-1},
 \\
 &\widetilde{\mathbf{H}}_z(s, d_z) =   \begin{pmatrix} 
\mathbf{H} & \mathbf{0} \\
 \mathbf{0}  & \mathbf{H}\\
 \end{pmatrix}
\left( \mathbf{I} + 
%%%%%
\begin{pmatrix} 
\mathbf{H}_z^{-1} & \mathbf{0} \\
 \mathbf{0}  & \mathbf{H}_z^{-1}\\
 \end{pmatrix} 
 %%%%%%%
 \left(
 \frac{(1- r_z)c}{2 sS_z}
 \begin{pmatrix} 
\mathbf{B}_z{(-1,-1)}  & \mathbf{0} \\
 \mathbf{0}  & \mathbf{B}_z{(1,1)} \\
 \end{pmatrix} 
 \right)
 \right)^{-1},
 \end{align*}
}
with $c = Z/\rho$. For $\widetilde{\mathbf{H}}$  a complex matrix, let  $\widetilde{\mathbf{H}}^{\dagger}$ denote  the complex conjugate transpose of  $\widetilde{\mathbf{H}}$ and $\Re{\left(\widetilde{\mathbf{H}}\right)} = \frac{1}{2}\left(\widetilde{\mathbf{H}} + \widetilde{\mathbf{H}}^{\dagger}\right)$.  With $\omega_y = \omega_z = 1$, the last two terms in \eqref{eq:disc_elemental_pde1_pml_laplace} vanish.  
After eliminating the velocity fields, the Laplace transformed, the two element DGSEM  approximation of the PML,  \eqref{eq:disc_elemental_pde1_pml_laplace}--\eqref{eq:disc_elemental_pde4_pml_laplace},   can be written as 
%%%%%%
%%%%%%
 { 
 \small
\begin{equation}\label{eq:disc_elemental_pde_secondorder_pml}
\begin{split}
 &
  ss^*
 %%%%%%%
 \begin{pmatrix} 
\mathbf{H} {\boldsymbol{\kappa}}^{-1}& \mathbf{0} \\
 \mathbf{0}  & \mathbf{H} {\boldsymbol{\kappa}}^{-1}\\
 \end{pmatrix}  
 %%%%%%%
 \begin{bmatrix} 
s\widetilde{\boldsymbol{p}}^{1} \\
%%%%%%%
s\widetilde{\boldsymbol{p}}^{2}
 \end{bmatrix} 
%%%%%%%
+ 
\frac{(s^*S_x^*)}{S_x} 
%%%%%%%
\widetilde{\mathbf{D}}_x^{\dagger}  
\begin{pmatrix} 
\boldsymbol{\rho}^{-1} & \mathbf{0} \\
 \mathbf{0}  & \boldsymbol{\rho}^{-1} \\
 \end{pmatrix} 
 %%%%%%%
  \widetilde{\mathbf{H}}_x
 \widetilde{\mathbf{D}}_x
 \begin{bmatrix} 
\widetilde{\boldsymbol{p}}^{1} \\
%%%%%%%
\widetilde{\boldsymbol{p}}^{2}
 \end{bmatrix} 
%%%%%%%
%%%%%%%
+ 
\frac{(sS_y)^*}{S_y} 
%%%%%%%
\widetilde{\mathbf{D}}_y^{\dagger}  
\begin{pmatrix} 
\boldsymbol{\rho}^{-1} & \mathbf{0} \\
 \mathbf{0}  & \boldsymbol{\rho}^{-1} \\
 \end{pmatrix} 
 %%%%%%%
  \widetilde{\mathbf{H}}_y
 \widetilde{\mathbf{D}}_y
 \begin{bmatrix} 
\widetilde{\boldsymbol{p}}^{1} \\
%%%%%%%
\widetilde{\boldsymbol{p}}^{2}
 \end{bmatrix} 
%%%%%%%
%%%%%%%
\\
&+ 
\frac{\left(sS_z\right)^*}{S_z} 
%%%%%%%
\widetilde{\mathbf{D}}_z^{\dagger}  
\begin{pmatrix} 
\boldsymbol{\rho}^{-1} & \mathbf{0} \\
 \mathbf{0}  & \boldsymbol{\rho}^{-1} \\
 \end{pmatrix} 
 %%%%%%%
  \widetilde{\mathbf{H}}_z
 \widetilde{\mathbf{D}}_z
 \begin{bmatrix} 
\widetilde{\boldsymbol{p}}^{1} \\
%%%%%%%
\widetilde{\boldsymbol{p}}^{2}
 \end{bmatrix} 
 +
 \frac{1}{2S_xZ}
%%%%%
\begin{pmatrix} 
\mathbf{H}_y\mathbf{H}_z & \mathbf{0} \\
 \mathbf{0}  & \mathbf{H}_y\mathbf{H}_z\\
 \end{pmatrix} 
 %%%%%%%
 \begin{pmatrix} 
\mathbf{B}_x{(1,1)}  & -\mathbf{B}_x{(1,-1)}  \\
-\mathbf{B}_x^T{(1,-1)}   & \mathbf{B}_x{(-1,-1)}\\
 \end{pmatrix} 
 %%%%%%%
  %%%%%%%
 %%%%%%%
 \begin{bmatrix} 
\widetilde{\boldsymbol{p}}^{1} \\
%%%%%%%
\widetilde{\boldsymbol{p}}^{2}
 \end{bmatrix} 
%%%%%%%
  \\
%%%%%%%
&
\frac{1+r_x}{2ZS_x}
%%%%%
\begin{pmatrix} 
\mathbf{H}_y\mathbf{H}_z & \mathbf{0} \\
 \mathbf{0}  & \mathbf{H}_y\mathbf{H}_z\\
 \end{pmatrix}  
 %%%%%%%
 \begin{pmatrix} 
\mathbf{B}_z{(-1,-1)}  & \mathbf{0} \\
 \mathbf{0}  & \mathbf{B}_z{(1,1)}\\
 \end{pmatrix} 
 %%%%%%%
  %%%%%%%
 %%%%%%%
 \begin{bmatrix} 
\widetilde{\boldsymbol{p}}^{1} \\
%%%%%%%
\widetilde{\boldsymbol{p}}^{2}
 \end{bmatrix} 
%%%%%%%
+
\frac{1+r_y}{2ZS_y}
%%%%%
\begin{pmatrix} 
\mathbf{H}_x\mathbf{H}_z & \mathbf{0} \\
 \mathbf{0}  & \mathbf{H}_x\mathbf{H}_z\\
 \end{pmatrix} 
 %%%%%%%
 \begin{pmatrix} 
\mathbf{B}_y{(-1,-1)} & \mathbf{0} \\
 \mathbf{0}  & \mathbf{B}_y{(1,1)}\\
 \end{pmatrix} 
 %%%%%%%
  %%%%%%%
 %%%%%%%
 \begin{bmatrix} 
\widetilde{\boldsymbol{p}}^{1} \\
%%%%%%%
\widetilde{\boldsymbol{p}}^{2}
 \end{bmatrix} 
%%%%%%%
\\
&
+
\frac{1+r_z}{2ZS_z}
%%%%%
\begin{pmatrix} 
\mathbf{H}_x\mathbf{H}_y & \mathbf{0} \\
 \mathbf{0}  & \mathbf{H}_x \mathbf{H}_y\\
 \end{pmatrix} 
 %%%%%%%
 \begin{pmatrix} 
\mathbf{B}_z{(-1,-1)}  & \mathbf{0} \\
 \mathbf{0}  & \mathbf{B}_z{(1,1)}\\
 \end{pmatrix} 
 %%%%%%%
  %%%%%%%
 %%%%%%%
 \begin{bmatrix} 
\widetilde{\boldsymbol{p}}^{1} \\
%%%%%%%
\widetilde{\boldsymbol{p}}^{2}
 \end{bmatrix} 
%%%%%%%
\\
&
=   
%%%%%%%
 \begin{pmatrix} 
  %%%%%%%
\mathbf{H} {\boldsymbol{\kappa}}^{-1}& \mathbf{0} \\
 \mathbf{0}  & \mathbf{H} {\boldsymbol{\kappa}}^{-1}\\
 %%%%%%%
 \end{pmatrix} 
 %%%%%%%
   \begin{bmatrix}
 \left( \frac{1}{S_x}{\boldsymbol{f}}_p- \frac{{\boldsymbol{\kappa}}}{sS_yS_x}{\boldsymbol{f}}_{\sigma} -  \frac{{\boldsymbol{\kappa}}}{sS_zS_x}{\boldsymbol{f}}_{\psi}  \right) \\
 %%%%%%
  \left( \frac{1}{S_x}{\boldsymbol{f}}_p- \frac{{\boldsymbol{\kappa}}}{sS_yS_x}{\boldsymbol{f}}_{\sigma} -  \frac{{\boldsymbol{\kappa}}}{sS_zS_x}{\boldsymbol{f}}_{\psi}  \right)
  \end{bmatrix} 
  +
\frac{(s^*S_x^*)}{S_x} 
%%%%%%%
\widetilde{\mathbf{D}}_x^{\dagger}  
\begin{pmatrix} 
\boldsymbol{\rho}^{-1} & \mathbf{0} \\
 \mathbf{0}  & \boldsymbol{\rho}^{-1} \\
 \end{pmatrix} 
 %%%%%%%
  \widetilde{\mathbf{H}}_x
 \begin{bmatrix} 
\frac{\boldsymbol{\rho}}{S_x} \boldsymbol{f}_u \\
%%%%%%%
\frac{\boldsymbol{\rho}}{S_x} \boldsymbol{f}_u
 \end{bmatrix} 
%%%%%%%
\\
&
+ 
\frac{(s^*S_y^*)}{S_y} 
%%%%%%%
\widetilde{\mathbf{D}}_y^{\dagger}  
\begin{pmatrix} 
\boldsymbol{\rho}^{-1} & \mathbf{0} \\
 \mathbf{0}  & \boldsymbol{\rho}^{-1} \\
 \end{pmatrix} 
 %%%%%%%
  \widetilde{\mathbf{H}}_y
 \begin{bmatrix} 
\frac{\boldsymbol{\rho}}{S_y} \boldsymbol{f}_v \\
%%%%%%%
\frac{\boldsymbol{\rho}}{S_y} \boldsymbol{f}_v
 \end{bmatrix} 
%%%%%%%
+
\frac{(s^*S_z^*)}{S_z} 
%%%%%%%
\widetilde{\mathbf{D}}_z^{\dagger} 
\begin{pmatrix} 
\boldsymbol{\rho}^{-1} & \mathbf{0} \\
 \mathbf{0}  & \boldsymbol{\rho}^{-1} \\
 \end{pmatrix} 
 %%%%%%%
  \widetilde{\mathbf{H}}_x
 \begin{bmatrix} 
\frac{\boldsymbol{\rho}}{S_z} \boldsymbol{f}_w\\
%%%%%%%
\frac{\boldsymbol{\rho}}{S_z} \boldsymbol{f}_w
 \end{bmatrix} 
%%%%%%%
\end{split}
\end{equation}
}
%%%%%%%%% %%%%%%%%%
%%%%%%%%% %%%%%%%%%
Note that it can be showed that $\Re{\left(\widetilde{\mathbf{H}}_j\left(s, d_j\right) \right)}  = \Re{\left(\widetilde{\mathbf{H}}_j\left(s, d_j\right) \right)}^T > 0$,  for all $j = x, y, z$.
%%%%%%%%% %%%%%%%%%

%%%%%%%%% %%%%%%%%%
%%%%%%%%% %%%%%%%%%
Let $ \widetilde{\mathbf{U}}\left(s\right)  = [\widetilde{\mathbf{U}}^1\left(s\right), \widetilde{\mathbf{U}}^2\left(s\right) ]^T$, $ \widetilde{\mathbf{V}}\left(s\right)  = [\widetilde{\mathbf{V}}^1\left(s\right), \widetilde{\mathbf{V}}^2\left(s\right) ]^T$, with 
$
\widetilde{\mathbf{U}}^1\left(s\right)  = \left(\widetilde{\mathbf{u}}_1^1\left(s\right), \widetilde{\mathbf{u}}_2^1\left(s\right), \widetilde{\mathbf{u}}^1_3\left(s\right), \widetilde{\mathbf{u}}_4^1\left(s\right)\right), \quad 
\widetilde{\mathbf{U}}^2\left(s\right)  = \left(\widetilde{\mathbf{u}}_1^2\left(s\right), \widetilde{\mathbf{u}}_2^2\left(s\right), \widetilde{\mathbf{u}}^2_3\left(s\right), \widetilde{\mathbf{u}}_4^2\left(s\right)\right),
$
%%%%%%%%% %%%%%%%%%
$
\widetilde{\mathbf{V}}^1\left(s\right)  = \left(\widetilde{\mathbf{v}}_1^1\left(s\right), \widetilde{\mathbf{v}}_2^1\left(s\right), \widetilde{\mathbf{v}}^1_3\left(s\right), \widetilde{\mathbf{v}}_4^1\left(s\right)\right), \quad 
\widetilde{\mathbf{V}}^2\left(s\right)  = \left(\widetilde{\mathbf{v}}_1^2\left(s\right), \widetilde{\mathbf{v}}_2^2\left(s\right), \widetilde{\mathbf{v}}^2_3\left(s\right), \widetilde{\mathbf{v}}_4^2\left(s\right)\right).
$
Introduce the discrete scalar product
{ 
 \small
\begin{equation}\label{eq:scalar_product_discrete_laplace}
\begin{split}
\Big\langle \widetilde{\mathbf{U}}\left(s\right) ,   \widetilde{\mathbf{V}}\left(s\right)\Big\rangle_h &:= a
 \begin{bmatrix} 
\widetilde{\boldsymbol{u}}^{1}_1 \\
%%%%%%%
\widetilde{\boldsymbol{u}}^{2}_1
 \end{bmatrix}^{\boldsymbol{\dagger}}
 %%%%%%%
 \begin{pmatrix} 
\mathbf{H} {\boldsymbol{\kappa}}^{-1}& \mathbf{0} \\
 \mathbf{0}  & \mathbf{H} {\boldsymbol{\kappa}}^{-1}\\
 \end{pmatrix}  
 %%%%%%%
 \begin{bmatrix} 
\widetilde{\boldsymbol{v}}^{1}_1 \\
%%%%%%%
\widetilde{\boldsymbol{v}}^{2}_1
 \end{bmatrix} 
%%%%%%%
+ 
\gamma_x
%%%%%%%
 \begin{bmatrix} 
\widetilde{\boldsymbol{u}}^{1}_2 \\
%%%%%%%
\widetilde{\boldsymbol{u}}^{2}_2
 \end{bmatrix}^{\boldsymbol{\dagger}}
\begin{pmatrix} 
\boldsymbol{\rho}^{-1} & \mathbf{0} \\
 \mathbf{0}  & \boldsymbol{\rho}^{-1} \\
 \end{pmatrix} 
 %%%%%%%
 \Re\left( \widetilde{\mathbf{H}}_x\right)
 \begin{bmatrix} 
\widetilde{\boldsymbol{v}}^{1}_2 \\
%%%%%%%
\widetilde{\boldsymbol{v}}^{2}_2
 \end{bmatrix}
%%%%%%%
\\
&
+ 
\gamma_y
%%%%%%%
 \begin{bmatrix} 
\widetilde{\boldsymbol{u}}^{1}_3 \\
%%%%%%%
\widetilde{\boldsymbol{u}}^{2}_3
 \end{bmatrix}^{\boldsymbol{\dagger}}
\begin{pmatrix} 
\boldsymbol{\rho}^{-1} & \mathbf{0} \\
 \mathbf{0}  & \boldsymbol{\rho}^{-1} \\
 \end{pmatrix} 
 %%%%%%%
 \Re\left( \widetilde{\mathbf{H}}_y\right)
 \begin{bmatrix} 
\widetilde{\boldsymbol{v}}^{1}_3 \\
%%%%%%%
\widetilde{\boldsymbol{v}}^{2}_3
 \end{bmatrix}
%%%%%%%
+ 
\gamma_z
%%%%%%%
 \begin{bmatrix} 
\widetilde{\boldsymbol{u}}^{1}_4 \\
%%%%%%%
\widetilde{\boldsymbol{u}}^{2}_4
 \end{bmatrix}^{\boldsymbol{\dagger}}
\begin{pmatrix} 
\boldsymbol{\rho}^{-1} & \mathbf{0} \\
 \mathbf{0}  & \boldsymbol{\rho}^{-1} \\
 \end{pmatrix} 
 %%%%%%%
 \Re\left( \widetilde{\mathbf{H}}_z\right)
 \begin{bmatrix} 
\widetilde{\boldsymbol{v}}^{1}_4 \\
%%%%%%%
\widetilde{\boldsymbol{v}}^{2}_4
 \end{bmatrix},
%%%%%%%
 \end{split}
\end{equation}
}
and the corresponding norm
{ 
 \small
\begin{equation}\label{eq:norm_discrete_laplace}
\begin{split}
\widetilde{\mathcal{E}}\left(\widetilde{\mathbf{U}}\left(s\right) \right): &= \sqrt{\Big\langle \widetilde{\mathbf{U}}\left(s\right) ,   \widetilde{\mathbf{U}}\left(s\right) \Big\rangle_h},
 \end{split}
\end{equation}
}
where $\Re{s} = a \ge 0$, $0\le \gamma_j \le |s|$ are defined in \eqref{eq:scaled_variables_1} and \eqref{eq:scaled_variables_2}.

\begin{theorem}\label{Theo:Numerical_Stability_PML_Laplace}
%%%%%%%%% %%%%%%%%%
Consider the  semi-discrete DGSEM approximation \eqref{eq:disc_elemental_pde1_pml_laplace}--\eqref{eq:disc_elemental_pde4_pml_laplace} of  the PML equation \eqref{eq:acoustic_pml_3D_Laplace1},  in the Laplace space with $\Re{s}  = a \ge 0$,   subject to boundary conditions \eqref{eq:boundary_condition_acoustic_laplace}. Let 
%%%%%%%%% %%%%%%%%%
$ \widetilde{\mathbf{U}}\left(s\right)  = [\widetilde{\mathbf{U}}^1\left(s\right), \widetilde{\mathbf{U}}^2\left(s\right) ]^T$, $ \widetilde{\mathbf{F}}\left(s\right)  = [\widetilde{\mathbf{F}}^1\left(s\right), \widetilde{\mathbf{F}}^2\left(s\right) ]^T$, where $ \widetilde{\mathbf{F}}^l\left(s\right),  = [\widetilde{\boldsymbol{f}}_1^{l}, \widetilde{\boldsymbol{f}}_2^{l}, \widetilde{\boldsymbol{f}}_3^{l}, \widetilde{\boldsymbol{f}}_4^{l} ]^T$, $l = 1, 2$, with 
%%%%%%%%% %%%%%%%%%
\[
 \begin{bmatrix} 
\widetilde{\boldsymbol{u}}_1^{1} \\
%%%%%%%
\widetilde{\boldsymbol{u}}_1^{2}
 \end{bmatrix}
 =
 \begin{bmatrix} 
s\widetilde{\boldsymbol{p}}^{1} \\
%%%%%%%
s\widetilde{\boldsymbol{p}}^{2}
 \end{bmatrix},
 %%%%%%
 \quad
 %%%%%
 \begin{bmatrix} 
\widetilde{\boldsymbol{u}}_2^{1} \\
%%%%%%%
\widetilde{\boldsymbol{u}}_2^{2}
 \end{bmatrix}
 =
\widetilde{\mathbf{D}}_x
 \begin{bmatrix} 
\widetilde{\boldsymbol{p}}^{1} \\
%%%%%%%
\widetilde{\boldsymbol{p}}^{2}
 \end{bmatrix},
 %%%%
 \quad
 %%%%
  \begin{bmatrix} 
\widetilde{\boldsymbol{u}}_3^{1} \\
%%%%%%%
\widetilde{\boldsymbol{u}}_3^{2}
 \end{bmatrix}
 =
\widetilde{\mathbf{D}}_y
 \begin{bmatrix} 
\widetilde{\boldsymbol{p}}^{1} \\
%%%%%%%
\widetilde{\boldsymbol{p}}^{2}
 \end{bmatrix},
 %%%%%
 \quad
 %%%%%
  \begin{bmatrix} 
\widetilde{\boldsymbol{u}}_4^{1} \\
%%%%%%%
\widetilde{\boldsymbol{u}}_4^{2}
 \end{bmatrix}
 =
\widetilde{\mathbf{D}}_z
 \begin{bmatrix} 
\widetilde{\boldsymbol{p}}^{1} \\
%%%%%%%
\widetilde{\boldsymbol{p}}^{2}
 \end{bmatrix}
\]
%%%%%%%%
%%%%%%%%
\[
 \begin{bmatrix} 
\widetilde{\boldsymbol{f}}_1^{1} \\
%%%%%%%
\widetilde{\boldsymbol{f}}_1^{2}
 \end{bmatrix}
 =
 \begin{bmatrix} 
\widetilde{\boldsymbol{f}}_p^{1} \\
%%%%%%%
\widetilde{\boldsymbol{f}}_p^{2}
 \end{bmatrix},
 %%%%%%
 \quad
 %%%%%
 \begin{bmatrix} 
\widetilde{\boldsymbol{f}}_2^{1} \\
%%%%%%%
\widetilde{\boldsymbol{f}}_2^{2}
 \end{bmatrix}
 =
 \begin{bmatrix} 
\frac{\boldsymbol{\rho}_1}{S_x}\widetilde{\boldsymbol{f}}_u^{1} \\
%%%%%%%
\frac{\boldsymbol{\rho}_2}{S_x}\widetilde{\boldsymbol{f}}_u^{2}
 \end{bmatrix},
 %%%%
 \quad
 %%%%
  \begin{bmatrix} 
\widetilde{\boldsymbol{f}}_3^{1} \\
%%%%%%%
\widetilde{\boldsymbol{f}}_3^{2}
 \end{bmatrix}
 =
   \begin{bmatrix}
 \frac{\boldsymbol{\rho}_1}{S_y}\widetilde{\boldsymbol{f}}_v^{1} \\
%%%%%%%
\frac{\boldsymbol{\rho}_2}{S_y}\widetilde{\boldsymbol{f}}_v^{2}
 \end{bmatrix}
 %%%%%
 \quad
 %%%%%
  \begin{bmatrix} 
\widetilde{\boldsymbol{f}}_4^{1} \\
%%%%%%%
\widetilde{\boldsymbol{f}}_4^{2}
 \end{bmatrix}
 =
   \begin{bmatrix}
\frac{\boldsymbol{\rho}_1}{S_z}\widetilde{\boldsymbol{f}}_w^{1} \\
%%%%%%%
\frac{\boldsymbol{\rho}_2}{S_z}\widetilde{\boldsymbol{f}}_w^{2}
 \end{bmatrix},
\]

%\[ 
%\widetilde{\mathbf{U}}^k\left(s\right)  = \left(s\widetilde{p}, \frac{1}{S_x} \frac{\partial \widetilde{p} }{\partial x}, \frac{1}{S_y} \frac{\partial \widetilde{p} }{\partial y}, \frac{1}{S_z} \frac{\partial \widetilde{p} }{\partial z}\right)^T, \quad 
%%\]
%%%%%%%%%% %%%%%%%%%
%%\[ 
%\widetilde{\mathbf{F}}\left(s\right)  = \left(\widetilde{F}_p(x,y,z), \frac{\rho}{S_x}   f_u(x,y,z),  \frac{\rho}{S_y}   f_v(x,y,z), \frac{\rho}{S_z}   f_w(x,y,z)\right)^T,
%\]
and
%%%%%
\[
\widetilde{\mathbf{f}}_p = \frac{1}{S_x}{\mathbf{f}}_p - \frac{\kappa}{sS_yS_x}{\mathbf{f}}_{\sigma}  -  \frac{\kappa}{sS_zS_x}{\mathbf{f}}_{\psi}.
\]
%%%%%%%%% %%%%%%%%%
If $\widetilde{\mathcal{E}}\left(\widetilde{\mathbf{U}}\left(s\right) \right) > 0$, and $\omega_y = 1$, $\omega_z = 1$, with element-wise constant damping  $d_x  \ge 0, d_y \ge 0,  d_z \ge 0$,  then we have
%%%%%%%%% %%%%%%%%%
{\small
\begin{equation}\label{eq:energy_estimate_pml_laplace_corner}
 \widetilde{\mathcal{E}}^2\left(\widetilde{\mathbf{U}}\left(s\right) \right) +\widetilde{\mathrm{BT}}(s) \le   \widetilde{\mathcal{E}}\left(\widetilde{\mathbf{U}}\left(s\right) \right)\widetilde{\mathcal{E}}\left(\widetilde{\mathbf{F}}\left(s\right) \right), \quad \widetilde{\mathrm{BT}}(s)   = {\Re\left(\frac{1}{S_x}\right)\widetilde{\mathbf{BT}}^{(x)} + \Re\left(\frac{1}{S_y}\right) \widetilde{\mathbf{BT}}^{(y)}  + \Re\left(\frac{1}{S_z}\right)\widetilde{\mathbf{BT}}^{(z)}  + \Re\left(\frac{1}{S_x}\right) \widetilde{\mathbf{IT}}^{(x)}} \ge 0,
\end{equation}
}
%%%%%%%%% %%%%%%%%%
where 
{
\small
\[
\widetilde{\mathbf{IT}}^{(x)}  = |s|\frac{1}{2Z}\sum_{j=1}^{N+1}\sum_{k=1}^{N+1}\left(\left|\widetilde{\boldsymbol{p}}^{1}\left(1, \eta_j, \theta_k, s\right) -\widetilde{\boldsymbol{p}}^{2}\left(-1, \eta_j, \theta_k, s\right)\right|^2\right)h_jh_k, 
\]
%%%%%%%%% %%%%%%%%%
%%%%%%%%% %%%%%%%%%
\begin{align*}
\widetilde{\mathbf{BT}}^{(x)}  = |s|\frac{1+r_x}{2Z}\sum_{j=1}^{N+1}\sum_{k=1}^{N+1}\left(\left|\widetilde{\boldsymbol{p}}^{1}\left(-1, \eta_j, \theta_k, s\right)|^2 + |\widetilde{\boldsymbol{p}}^{2}\left(1, \eta_j, \theta_k, s\right)\right|^2\right)h_jh_k,  
\end{align*}
%%%%%%
\[ 
\widetilde{\mathbf{BT}}^{(y)}  = |s|\frac{1+r_y}{2Z}\sum_{i=1}^{N+1}\sum_{k=1}^{N+1}\left(\left|\widetilde{\boldsymbol{p}}^{1}\left(\xi_i, -1, \theta_k, s\right)\right|^2 + \left|\widetilde{\boldsymbol{p}}^{1}\left(\xi_i, 1, \theta_k, s\right)\right|^2 + \left|\widetilde{\boldsymbol{p}}^{2}\left(\xi_i, -1, \theta_k, s\right)\right|^2 + \left|\widetilde{\boldsymbol{p}}^{2}\left(\xi_i, 1, \theta_k, s\right)\right|^2\right)h_ih_k,
\]
\[ 
\widetilde{\mathbf{BT}}^{(z)}  = |s|\frac{1+r_z}{2Z}\sum_{i=1}^{N+1}\sum_{j=1}^{N+1}\left(\left|\widetilde{\boldsymbol{p}}^{1}\left(\xi_i, \eta_j, -1, s\right)\right|^2 + \left|\widetilde{\boldsymbol{p}}^{1}\left(\xi_i, \eta_j, 1, s\right)\right|^2 + \left|\widetilde{\boldsymbol{p}}^{2}\left(\xi_i, \eta_j, 1, s\right)\right|^2 + \left|\widetilde{\boldsymbol{p}}^{2}\left(\xi_i, \eta_j, 1, s\right)\right|^2\right)h_ih_j.
\]
}
There are no nontrivial solutions with $\Re{s} = a > 0$.
\end{theorem}
%%%%%%%%% %%%%%%%%%
The  proof of the Theorem \ref{Theo:Numerical_Stability_PML_Laplace} can be found in \ref{sec:proof_disc_pml}.

%%%%%%%%% %%%%%%%%%
By Theorem \ref{Theo:Numerical_Stability_PML_Laplace},  with $\omega_y = 1, \omega_z = 1$, the semi-discrete approximation is asymptotically stable for all  element-wise constant PML damping  $d_x, d_y, d_z \ge 0$.
When $\omega_y \ne1, \omega_z \ne 1$, we could not find a discrete energy estimate, and there is no guarantee that the numerical approximation is stable. In the next section, we will perform numerical experiments to verify accuracy and the  stability analysis. With the standard choice $\omega_y = 0, \omega_z = 0$, the numerical experiments also demonstrate that the numerical approximation for the PML is unstable. 

%%%%%%%%%%%%%%%%%%%
%%%%%%%%%%%%%%%%%%%
%%%%%%%%%%%%%%%%%%%

  \section{Numerical experiments}
%%%%%%%%%%%%%
%%%%%%%%%%%%%
In this section, we present numerical experiments.  The experiments performed are aimed at quantifying numerical errors introduced by discretizing the PML as well as verifying the stability analysis of the last section, and demonstrating the power of the PML stabilizing flux fluctuations. 
%%%%%%%%%%%%%
We use  constant acoustic wave speed $c = 1.484$ km/s   and constant medium density $\rho = 1$ g/cm$^3$. Tensor product bases of Lagrange polynomials, of degrees $P = \{ 2, 4, 6, 8\}$,  are used with GLL,  GL, and GLR quadrature nodes, separately. Time-integration is performed using the high order  ADER  scheme \cite{DumbserPeshkovRomenski}  of the same order of accuracy with the spatial discretization. We will consider first a 2D problem and perform detailed numerical experiments, including showing convergence  for h- and p-refinement.   We will  proceed later to  3D numerical simulations of acoustic waves in an unbounded domain.
%Thus, for polynomial approximations of degree $P$, we will expect optimal asymptotic convergence rate of $P+1$. 
%%%%%%%%%%%%%

 We consider a 2D problem with $\partial /\partial z = 0$, $d_z(z) = 0$, the velocity component $w$ and the auxiliary variable $\psi$ drop out. 
%%%%%%%%%%%%%
%%%%%%%%%%%%%
%%%%%%%%%%%%%
%%%%%%%%%%%%%
 We will first consider the  vertical PML  strip problem (with $d_x(x) \ge 0$, $d_y(y) = 0$), that is a PML in the $x$-direction truncating the left and right boundaries, and proceed later to simulate a whole space problem surrounded by the PML. 
 %%%%%%%%%%%%%
%%%%%%%%%%%%%
 The later situation involves both the vertical (with $d_x(x) \ge 0$, $d_y(y) = 0$) and horizontal (with $d_x(x) = 0$, $d_y(y) \ge 0$) PML layers, and  PML corners  (with $d_x(x) \ge 0$, $d_y(y) \ge 0$) where both layers are simultaneously active.  In all experimental setups, we will demonstrate the  power of the PML stabilizing parameter  $\omega_y:= \omega = 1$.
%%%%%%%%%%%%%
%%%%%%%%%%%%%
%%%%%%%%%%%%%
\subsection{The vertical strip PML problem}
%%%%%%%%%%%%%
%%%%%%%%%%%%%
To begin with, consider a 2D $100$ \text{km} $\times$  $50$ \text{km} rectangular domain, with $\left(x, y\right)$ = $[-50$ \text{km}, $50$ \text{km}$]\times[0$ \text{km}, $50$ \text{km}$]$. 
In the $x$-direction we  introduce two additional layers, each of width $\delta = 10$ km, having $50\le |x| \le 60$ in which the PML equations are solved.  At all boundaries of the domain we set the absorbing boundary conditions, with  ${r}_x = {r}_y = 0$ in \eqref{eq:boundary_condition_acoustic}.  The absorbing boundary boundary conditions are used so that we can make comparisons with results in the literature \cite{KDuru2016}.  However, it is possible to use any well-posed boundary condition with $ |{r}_x| \le 1$,   $ |{r}_y| \le 1$.   
%%%%%%%%%%%%%%%
%%%%%%%%%%%%%%%

%%%%%%%%%%%%%%%
%%%%%%%%%%%%%%%
We set the initial condition
%%%%%%%%%%%%%
\begin{align}
p(x, y, t=0) = e^{-\log\left({2}\right)\frac{x^2 + (y-25)^2}{9}},
\end{align}
%%%%%%%%%%%%%
for the pressure field, and zero initial condition for the velocity fields and the auxiliary variable.
%%%%%%%%%%%%%
The damping profile is a cubic monomial
%%%%%%%%%%%%%
\begin{equation}\label{eq:damping_func}
\begin{split}
&d_x\left(x\right) = \left \{
\begin{array}{rl}
0 \quad {}  \quad {}& \text{if} \quad |x| \le 50 \quad \text{km},\\
d_0\Big(\frac{|x|-50}{\delta}\Big)^3  & \text{if}  \quad |x| \ge 50 \quad \text{km},
\end{array} \right.
\end{split}
\end{equation}
%%%%%%%%%%%%%
%%%%%%%%%%%%%
where $d_0\ge 0$ is the damping strength.
%%%%%%%%%%%%%
%%%%%%%%%%%%%

%%%%%%%%%%%%%%%%%%%%
%%%%%%%%%%%%%%%%%%%%
\subsubsection{Numerical stability}
\begin{figure}[h!]
\begin{subfigure}
    \centering
        %%%%%%%
%\stackunder[5pt]{\includegraphics[width=0.49\textwidth]{wave_field_pml_t100s.eps}}{$\theta_x = 0$.}%
\stackunder[5pt]{\includegraphics[width=0.325\linewidth]{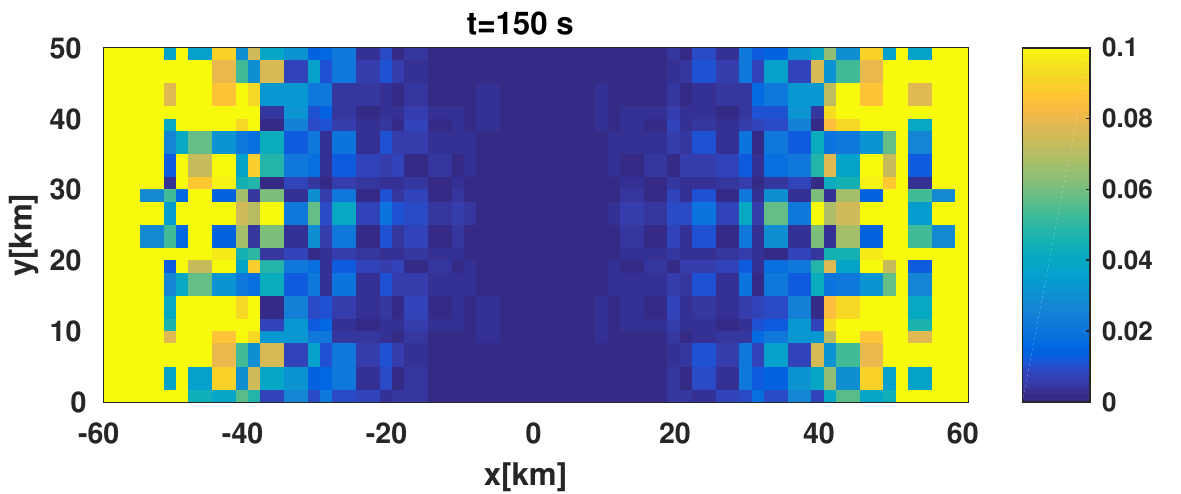}}{GLL.}%
\hspace{0.0cm}%
\stackunder[5pt]{\includegraphics[width=0.325\linewidth]{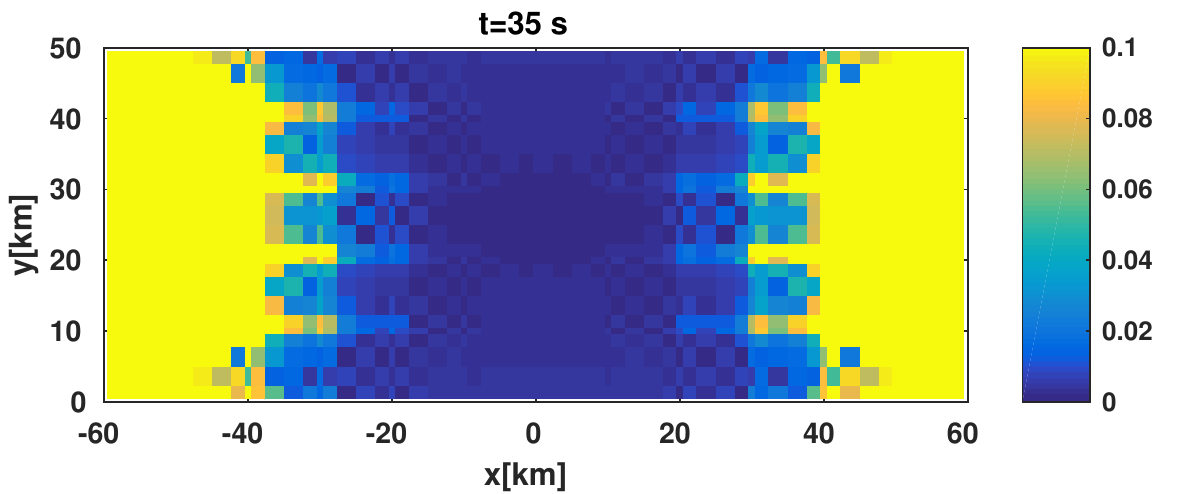}}{GL.}%
\hspace{0.0cm}%
\stackunder[5pt]{\includegraphics[width=0.325\linewidth]{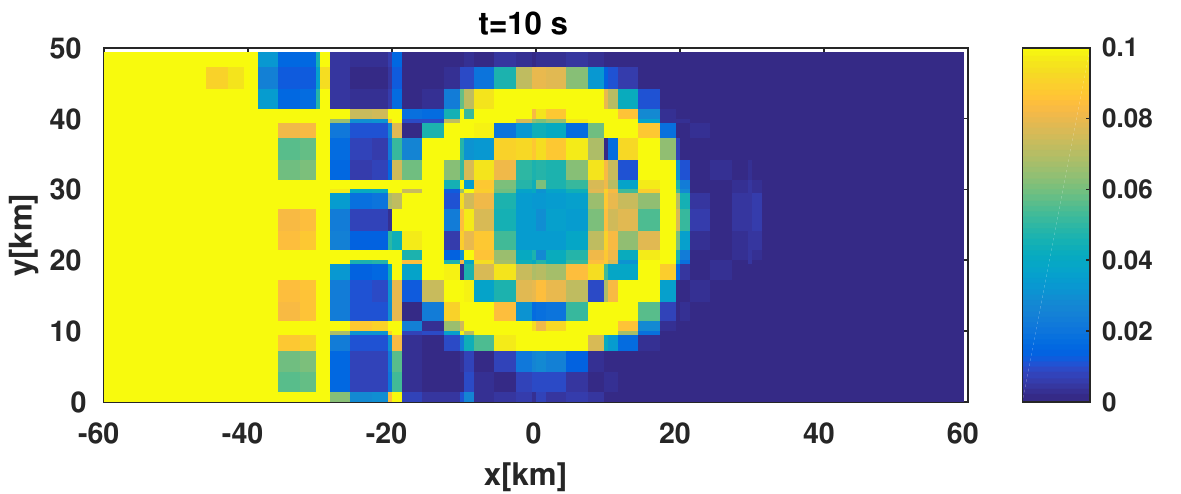}}{GLR.}%
%wave_field_pml_theta1_t100s.eps
%\stackunder[5pt]{\includegraphics[width=0.49\textwidth]{wave_field_t100s.eps}}{$\theta_x = 1$.}%
     \end{subfigure}
    \caption{Snapshots of the absolute pressure $|p|$ at $t= 150, 35$ and $10$ s, without PML flux fluctuation stabilization, $\omega_y = 0$. Numerical instabilities are generated by the PML boundaries.}
    \label{fig:Omega_0_strip}
\end{figure}

\begin{figure}[h!]
\begin{subfigure}
    \centering
        %%%%%%%
\stackunder[5pt]{\includegraphics[width=0.325\linewidth]{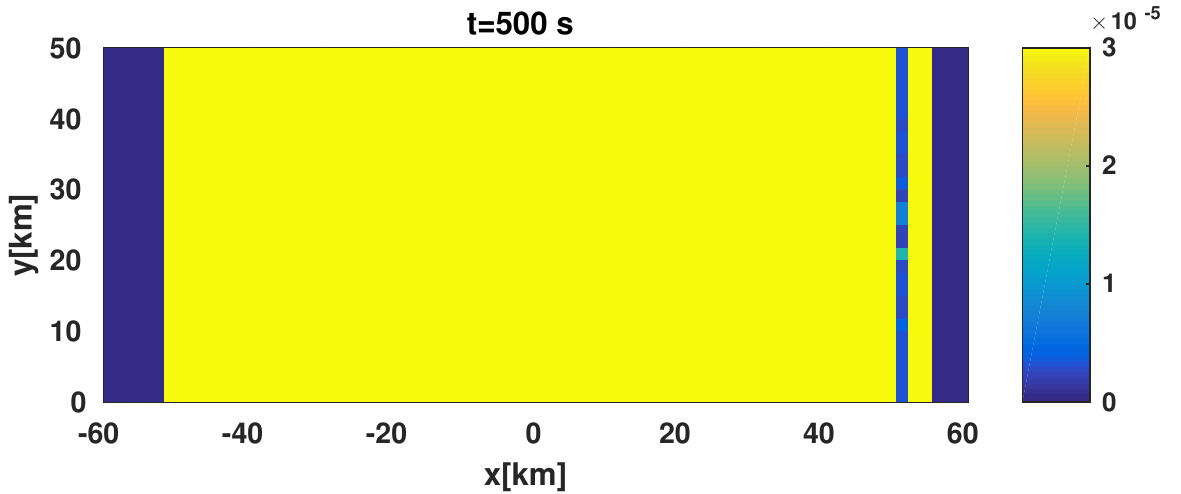}}{GLL.}%
\hspace{0.0cm}%
\stackunder[5pt]{\includegraphics[width=0.325\linewidth]{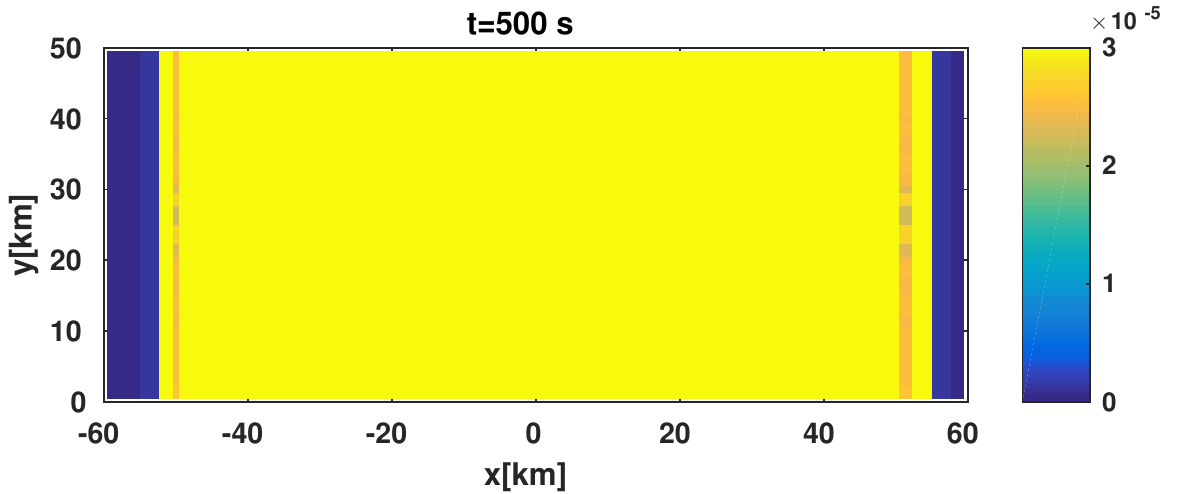}}{GL.}%
\hspace{0.0cm}%
\stackunder[5pt]{\includegraphics[width=0.325\linewidth]{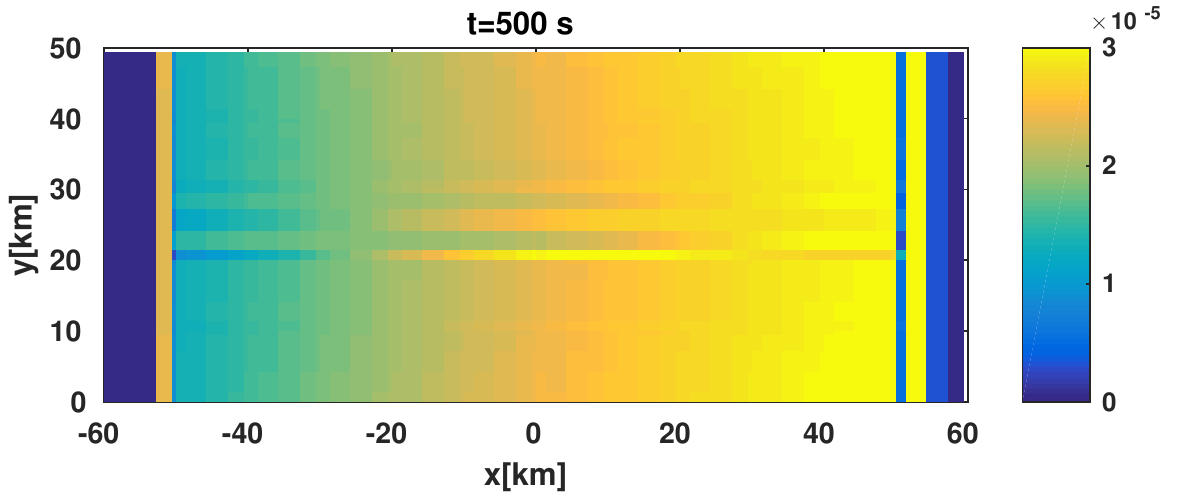}}{GLR.}%
%%%%
     \end{subfigure}
    \caption{Snapshots  of the absolute pressure $|p|$ at the final time $t= 500$ s, with PML flux fluctuation stabilization, $\omega_y = 1$. There are no instabilities. Solutions are stable after a very long time.}
    \label{fig:Omega_1_strip}
\end{figure}

\begin{figure}[h!]
\begin{subfigure}
    \centering
        %%%%%%%
\stackunder[5pt]{\includegraphics[width=0.33\linewidth]{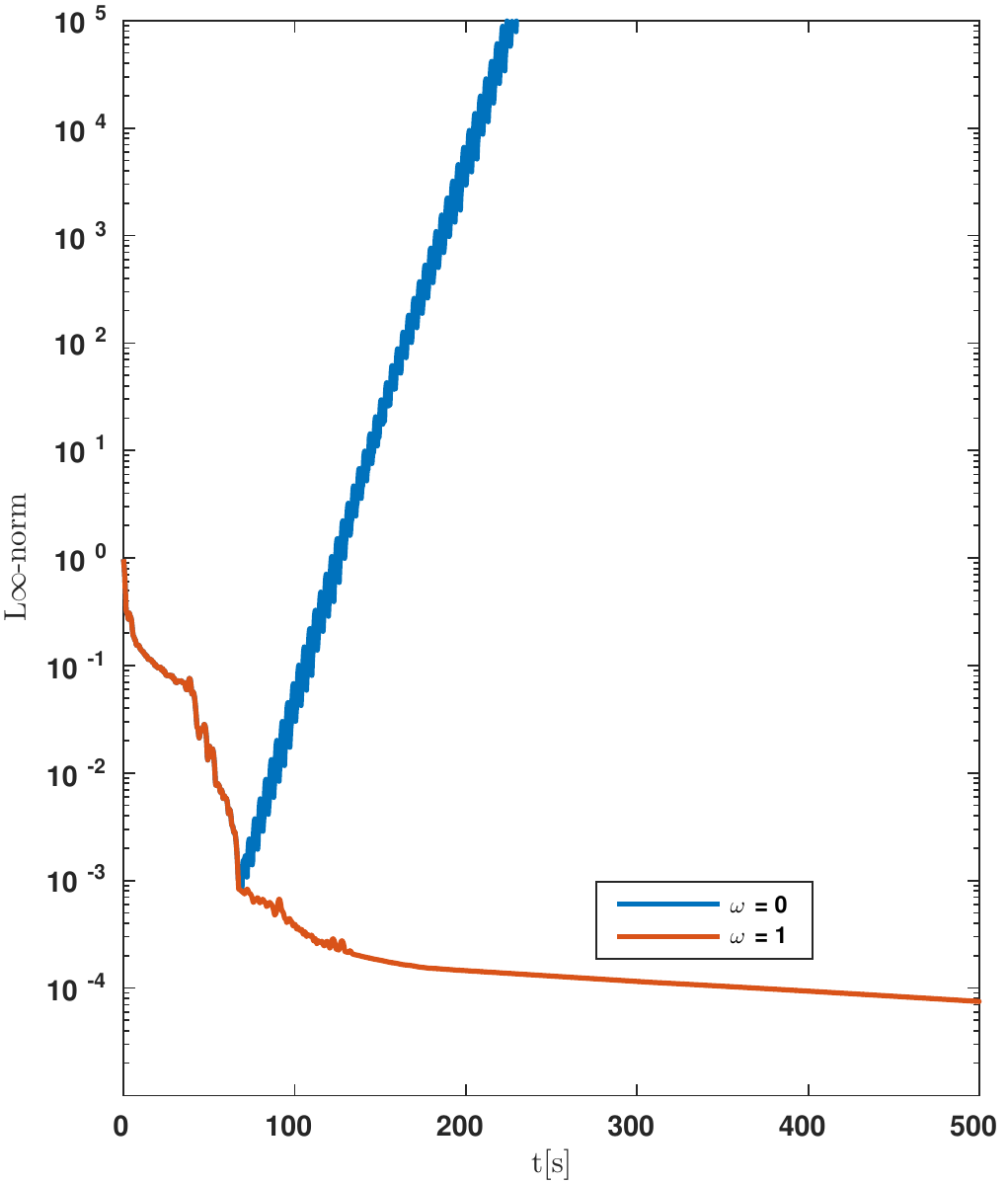}}{GLL.}%
\hspace{0.0cm}%
\stackunder[5pt]{\includegraphics[width=0.33\linewidth]{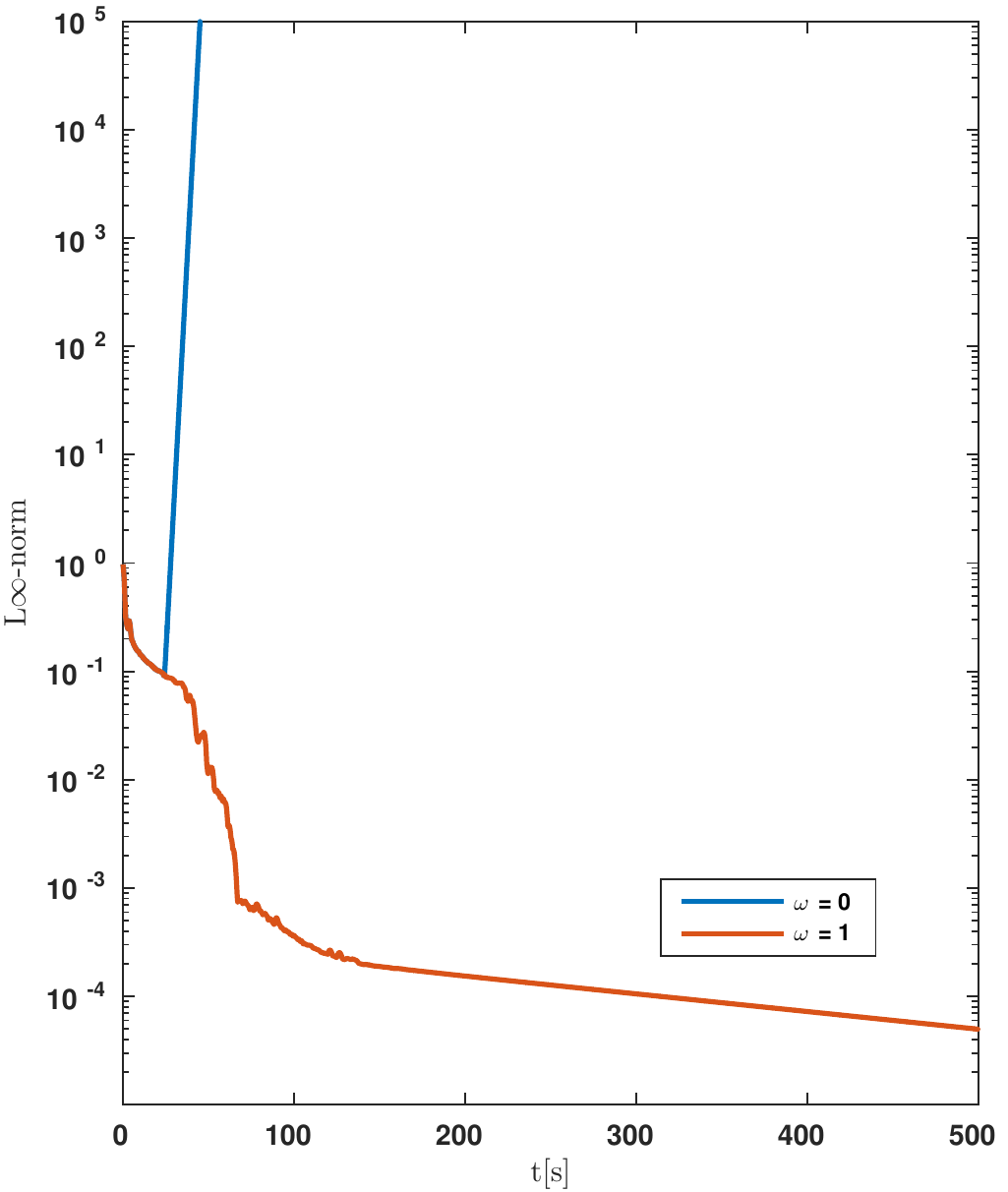}}{GL.}%
\hspace{0.0cm}%
\stackunder[5pt]{\includegraphics[width=0.33\linewidth]{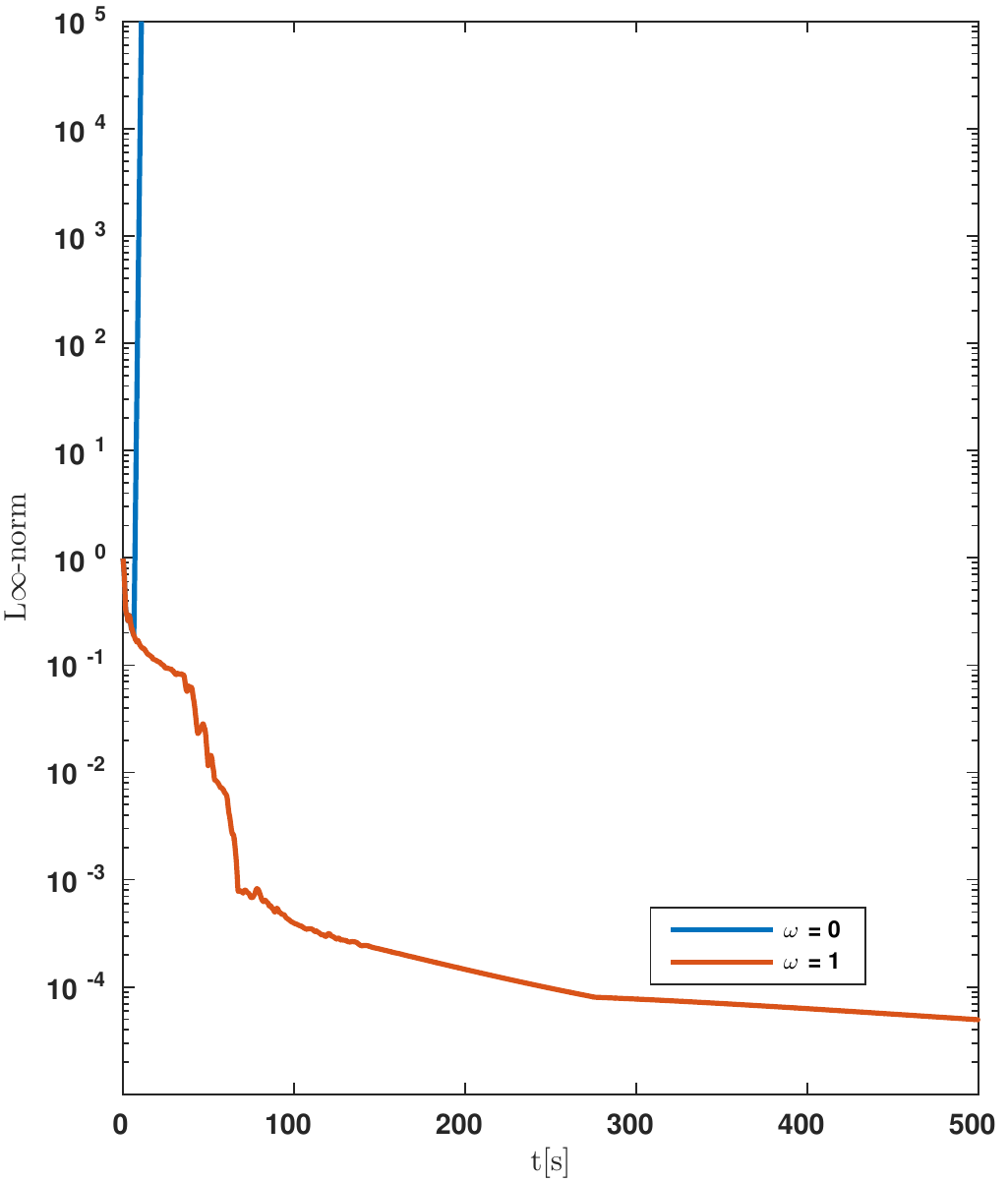}}{GLR.}%
 %%%%%%%
     \end{subfigure}
    \caption{The time history of the L$_{\infty}$-norm of the pressure field for $\omega_y = 0, 1$.}
    \label{fig:Time_series_error_strip}
\end{figure}

%%%%%%%%%%%%%%%%%%%%
%%%%%%%%%%%%%%%%%%%%
Here, we investigate numerical stability. Set the damping strength $d_0 = 8$. We discretize the domain with a uniform element size $\Delta{x} = \Delta{y} = 10$ km, spanning the PML with only one DGSEM element, and approximate the solution by a polynomial of degree $P = 4$. We use the time step
\[
\Delta{t} = \frac{CFL}{\left(2P+1\right) c} \min{\left(\Delta{x}, \Delta{y}\right)},
\]
with the $CFL = 0.35$ number and the acoustic wave speed $c = 1.484$ km/s. The final time is $t = 500 $ s.  The snapshots of the absolute pressure fields are plotted in Figure \ref{fig:Omega_0_strip} and Figure  \ref{fig:Omega_1_strip} and the time history  L$_{\infty}$-norm of the pressure field plotted in Figure  \ref{fig:Time_series_error_strip}. First, we compute the solutions by directly appending the PML terms, without the PML stabilizing term, that is $\omega_y:= \omega = 0$, in the auxiliary differential equation. From Figure \ref{fig:Omega_0_strip},  the solution in the PML explodes after some time steps. See also Figure \ref{fig:Time_series_error_strip}.  The initiation time of the explosive numerical mode depends on the quadrature rule used and the mesh resolution. Note that the GLR quadrature nodes are asymmetric. It is not surprising that for the GLR rule, numerical instability is also asymmetric, appearing at the left boundary first before spreading to the entire computational domain.  On a finer mesh the growth persists, but it starts a much later time.
%%%%%%%%%%%%%%%%%%%%
%%%%%%%%%%%%%%%%%%%%

%%%%%%%%%%%%%%%%%%%%
%%%%%%%%%%%%%%%%%%%%
When the PML stabilizing term is present, $\omega = 1$, the solution is stable, after a very long time, $t = 500$ s. This is clearly demonstrated by the small amplitude $\sim 10^{-5}$ of the pressure,  in Figure \ref{fig:Omega_1_strip} and Figure \ref{fig:Time_series_error_strip}, at the final time $t = 500$ s.
%%%%%%%%%%%%%%%%%%%%
%%%%%%%%%%%%%%%%%%%%

\subsubsection{Accuracy and convergence}
Next we will verify the accuracy of the DGSEM approximation of the PML. From the last experiments we know that the PML stabilization parameter $\omega = 1$ is critical for numerical stability. Therefore, we set $\omega = 1$ to ensure stability.
 We will use 
%%%%%%%%%%%%%
%%%%%%%%%%%%%
\begin{equation}\label{eq:damping_coef}
d_0 = \frac{4c}{2\delta}\ln{\frac{1}{\mathrm{tol}}},
\end{equation}
%%%%%%%%%%%%%
%%%%%%%%%%%%%
where $c = 1.484$ km/s is the acoustic wave speed, $\delta = 10$ km is the width of the PML, and  $\mathrm{tol}$ is the magnitude of the relative PML error \cite{KDuru2016}. We choose 
%%%%%%%%%%%%%
%%%%%%%%%%%%%
\begin{equation}\label{eq:tol_coef}
\mathrm{tol} = C_0\left[\frac{1}{\delta}\frac{\Delta{x}}{P+1}\right]^{P+1},
\end{equation}
%%%%%%%%%%%%%
%%%%%%%%%%%%%
where $C_0>0$ is an empirically determined constant.  We will use $C_0 = 10$ through the experiments. There are three parameters that control the PML error, the element size $\Delta{x}$, the width of the PML $\delta$, and the degree of the DGSEM polynomial approximation $P$. If we fix the PML width $\delta =  \delta_0 > 0$ and the mesh size $\Delta{x} = \Delta{x}_0 > 0$, and increase the polynomial degree $P \to \infty$, then the PML error, $\mathrm{tol}$,  will converge to zero exponentially. This is called $p$-convergence. On the other hand, we have $h$-convergence, if we fix the PML width $\delta = \delta_0 > 0$,   the polynomial degree $P = P_0 \ge 0$,  and decrease mesh size $\Delta{x} \to 0$, then the PML error will converge to zero at the rate $\mathrm{tol}\sim O(\Delta{x}^{P+1})$.  Finally,  if we fix the mesh size  $\Delta{x} = \Delta{x}_0 > 0$,  the polynomial degree $P = P_0 \ge 0$,  and increase the  PML width $\delta \to \infty$, then the PML error will  also converge to zero at the rate $\mathrm{tol}\sim O((1/\delta)^{P+1})$.
%%%%%%%%%%%%%

%%%%%%%%%%%%%
In the coming experiments, we will fix the PML width $\delta =  10 $ km and verify both the $p$-convergence and $h$-convergence. We will begin with $p$-convergence. To do this we fix the mesh-size $\Delta{x} = 5$ km and vary the polynomial degree $P = 2, 4, 6, 8$. In order to evaluate errors we compute a reference solution in large domain so that the reflection from the boundaries do not reach the interior at the final time. We chose the final time $t = 60$ s, so that the waves reach the PML boundaries,  and the PML errors propagate back into the computational domain. By comparing the reference solution to the PML solution in the interior, $|x| \le 50 $ km, in the L$_{\infty}$-norm, we obtain an accurate measure of the total PML error.
%%%%%%%%%%%%%

\begin{figure}[h!]
\begin{subfigure}
    \centering
        %%%%%%%
%\stackunder[5pt]{\includegraphics[width=0.49\textwidth]{wave_field_pml_t100s.eps}}{$\theta_x = 0$.}%
\stackunder[5pt]{\includegraphics[width=0.33\linewidth]{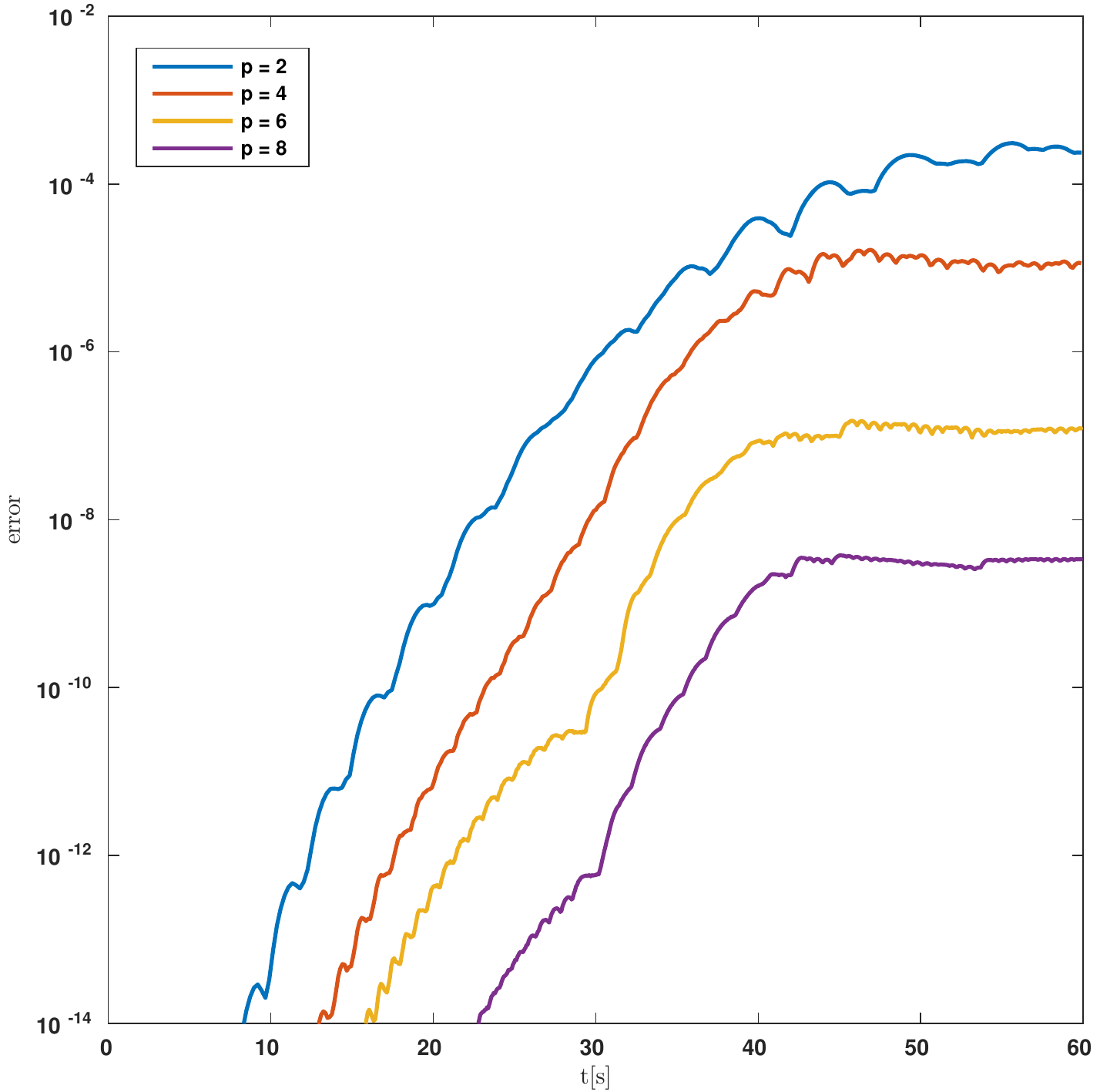}}{GLL.}%
\hspace{0.0cm}%
\stackunder[5pt]{\includegraphics[width=0.33\linewidth]{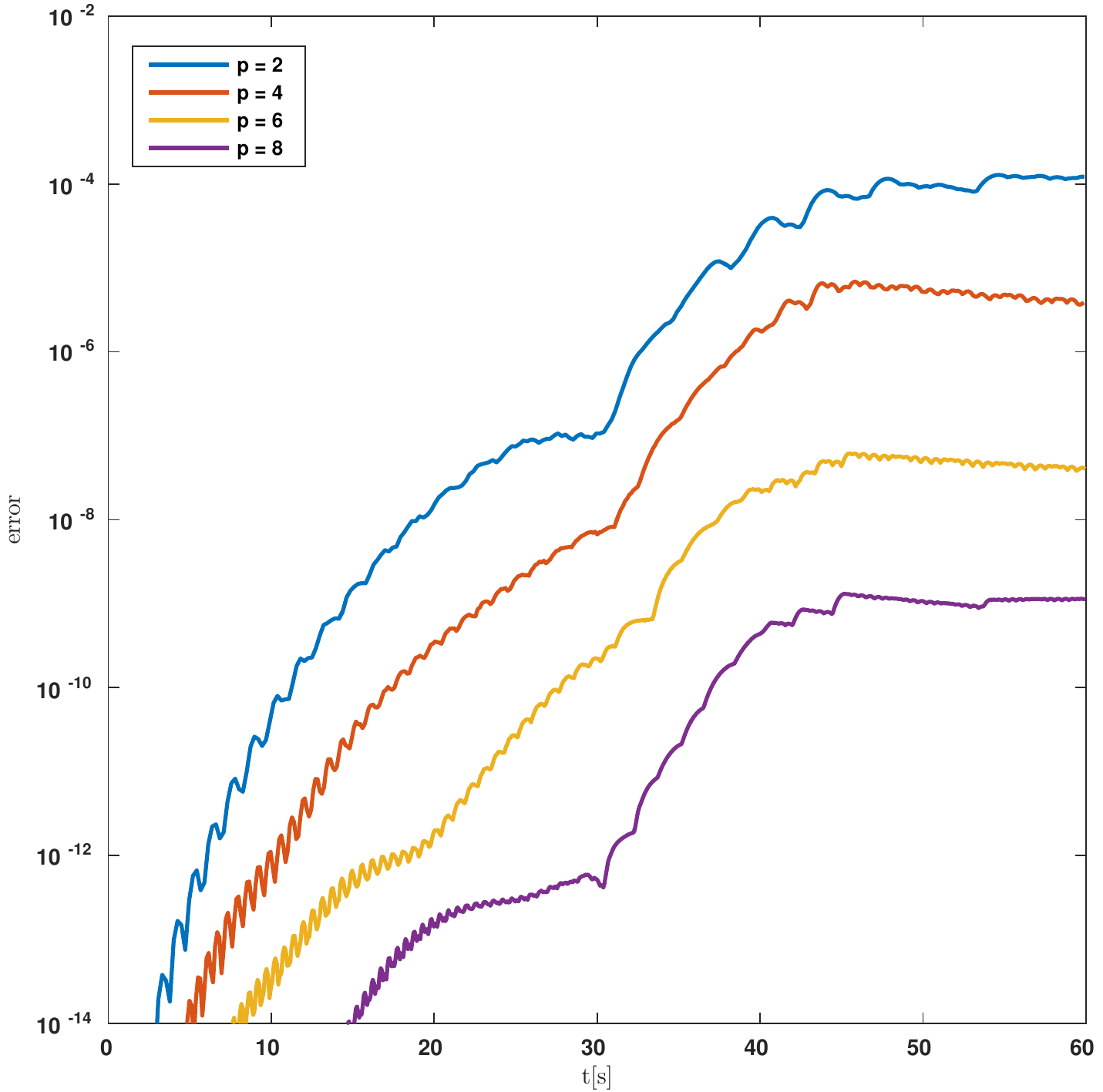}}{GL.}%
\hspace{0.0cm}%
\stackunder[5pt]{\includegraphics[width=0.33\linewidth]{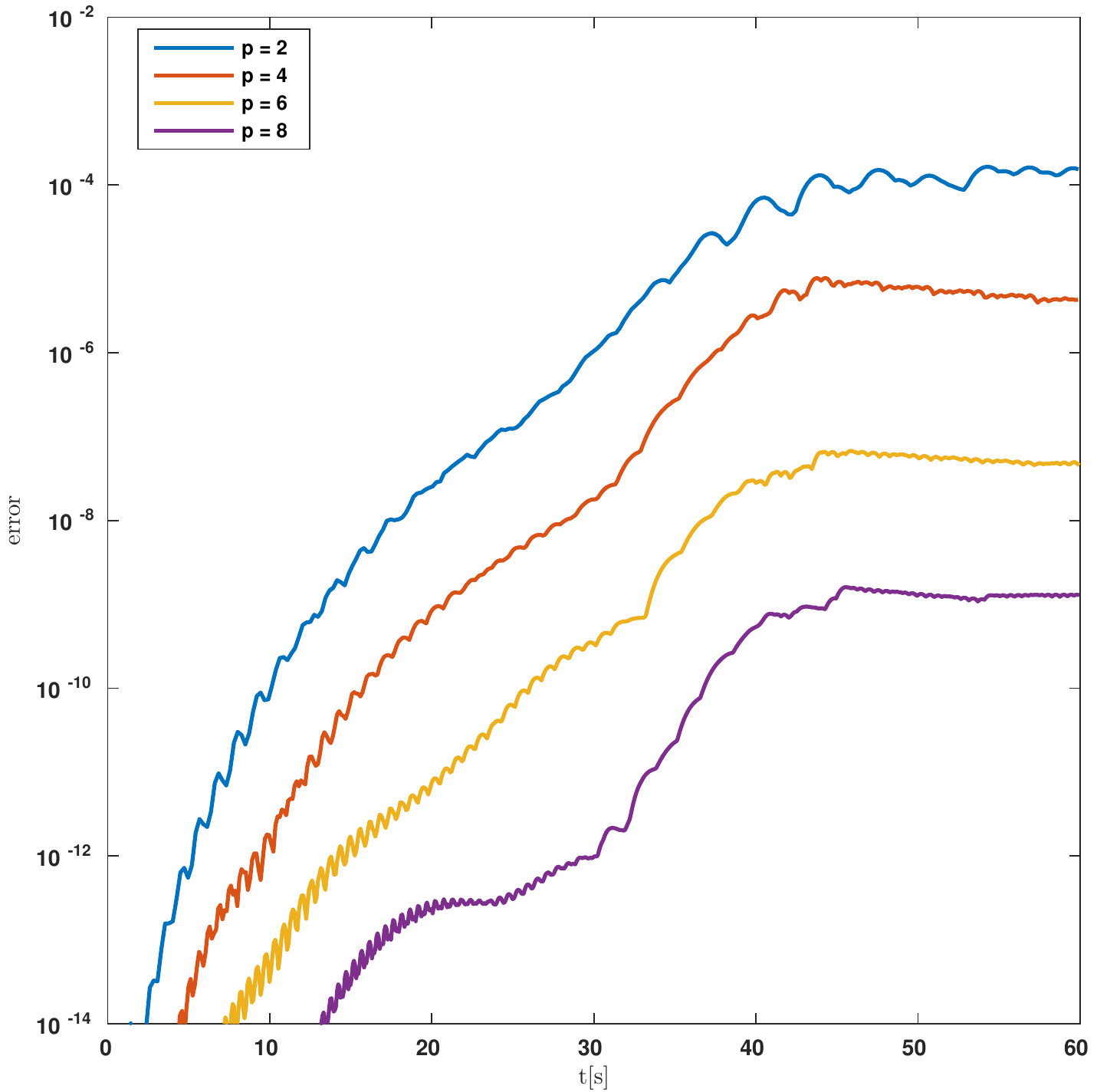}}{GLR.}%
%wave_field_pml_theta1_t100s.eps
%\stackunder[5pt]{\includegraphics[width=0.49\textwidth]{wave_field_t100s.eps}}{$\theta_x = 1$.}%
     \end{subfigure}
    \caption{ Time history of  PML errors using element size $\Delta{x} = 5$ and  polynomial  degrees $P = \{ 2, 4, 6, 8\}$.}
    \label{fig:Time_series_error_strip_p_convergence}
\end{figure}

\begin{figure} [htb!]
 \centering
%%%%%%%%%%%%%%%%%%
{\includegraphics[width=0.5\linewidth]{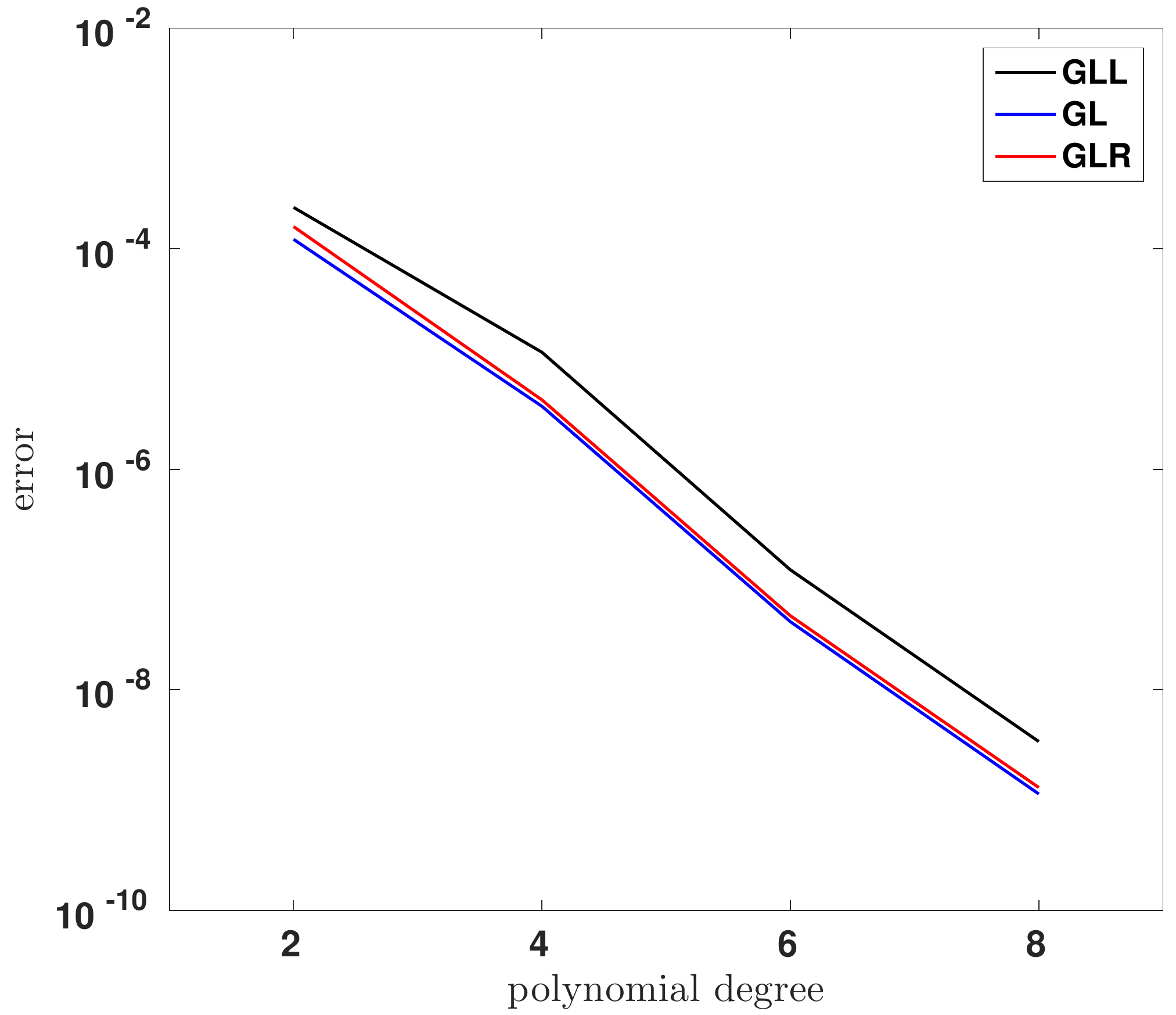}\label{fig:EnergyVsTime_hf}}
%{\includegraphics[width=0.495\linewidth]{pml_h_convergence.eps}\label{fig:EnergyVsTime_hf}}
%%%%%%%%%%%%%%%%%%
%%%
%%%%%%%%%%%%%%%%%%%
%%%%%%%%%%%%%%%%%%%
 \caption{\textit{Spectral  accurate convergence rates of PML errors  for p-refinement}}
   \label{fig:Final_error_strip_p}
\end{figure}

Time histories of the PML errors are shown in Figure \ref{fig:Time_series_error_strip_p_convergence} for  the polynomial  degrees $P = \{ 2, 4, 6, 8\}$ and quadrature rules GLL, GL, and GLR. The error at the final time $t = 60$ s are plotted against the polynomial degree $P$, in Figure \ref{fig:Final_error_strip_p}. Note that the PML error converges spectrally to zero.
%%%%%%%%%%%%%
%%%%%%%%%
%%%%%%%%%
%%%%%%%%%

\begin{figure} [htb!]
 \centering
%%%%%%%%%%%%%%%%%%
%{\includegraphics[width=0.45\linewidth]{pml_spectral_convergence.eps}\label{fig:EnergyVsTime_hf}}
{\includegraphics[width=0.5\linewidth]{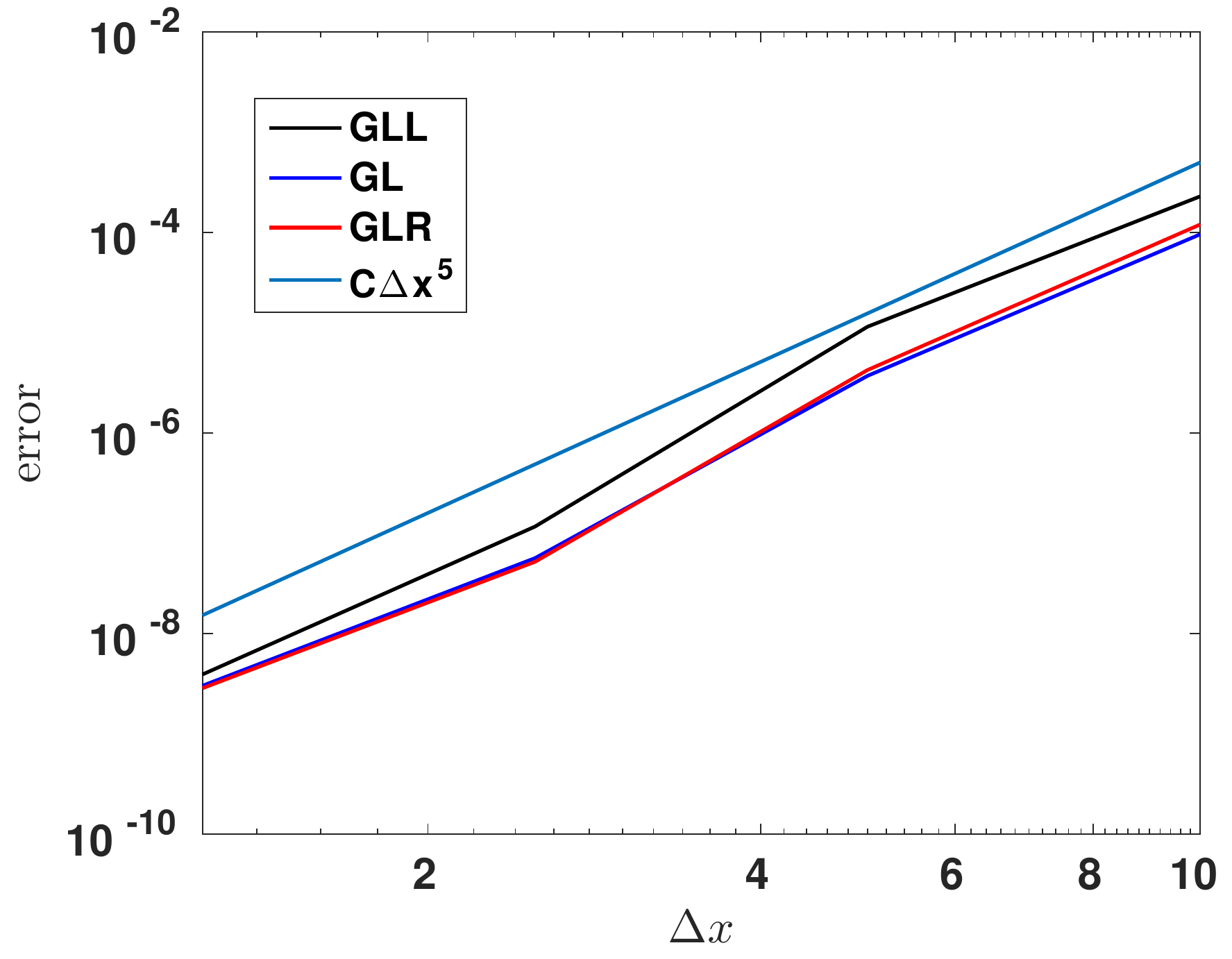}\label{fig:EnergyVsTime_hf}}
%%%%%%%%%%%%%%%%%%
%%%%%%%%%%%%%%%%%%
%%%%%%%%%%%%%%%%%%
 \caption{\textit{High order accurate convergence rates of PML errors  for h-refinement.}}
   \label{fig:Final_error_strip_h}
\end{figure}

%%%%%%%%%
Next we consider $h$-convergence. We now fix the polynomial degree $P = 4$ and decrease the mesh size, $\Delta{x} = \{ 10, 5, 2.5, 1.25\}$, and again use the quadrature rules, GLL, GL, and GLR. The PML error is expected to converge to zero at the rate $\mathrm{tol}\sim O(\Delta{x}^{5})$.
Time histories of the error are plotted in Figure \ref{fig:Time_series_error_strip_p_convergence}. The error at the final time $t = 60$ s are plotted against the mesh size $\Delta{x}$, in Figure \ref{fig:Final_error_strip_h}. Note that the error is parallel to the theoretical convergence  rate $\mathrm{tol}\sim O(\Delta{x}^{5})$. Therefore, the PML error converges zero  optimally,  at the rate $O(\Delta{x}^{5})$.
%%%%%%%%%
%%%%%%%%%
%\begin{figure} [htb!]
% \centering
%%%%%%%%%%%%%%%%%%%
\begin{figure}[h!]
\begin{subfigure}
    \centering
        %%%%%%%
%\stackunder[5pt]{\includegraphics[width=0.49\textwidth]{wave_field_pml_t100s.eps}}{$\theta_x = 0$.}%
\stackunder[5pt]{\includegraphics[width=0.33\linewidth]{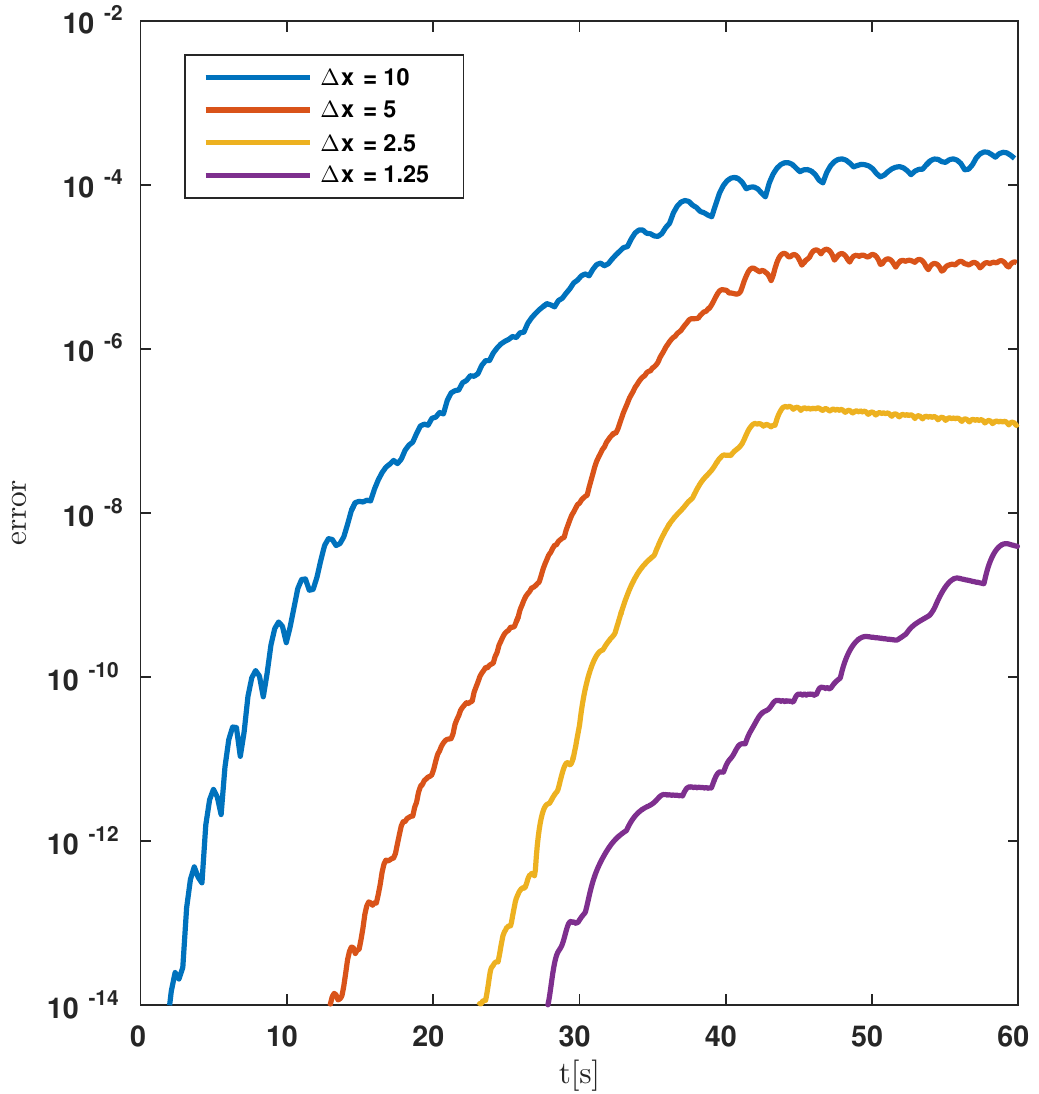}}{GLL.}%
\hspace{0.0cm}%
\stackunder[5pt]{\includegraphics[width=0.33\linewidth]{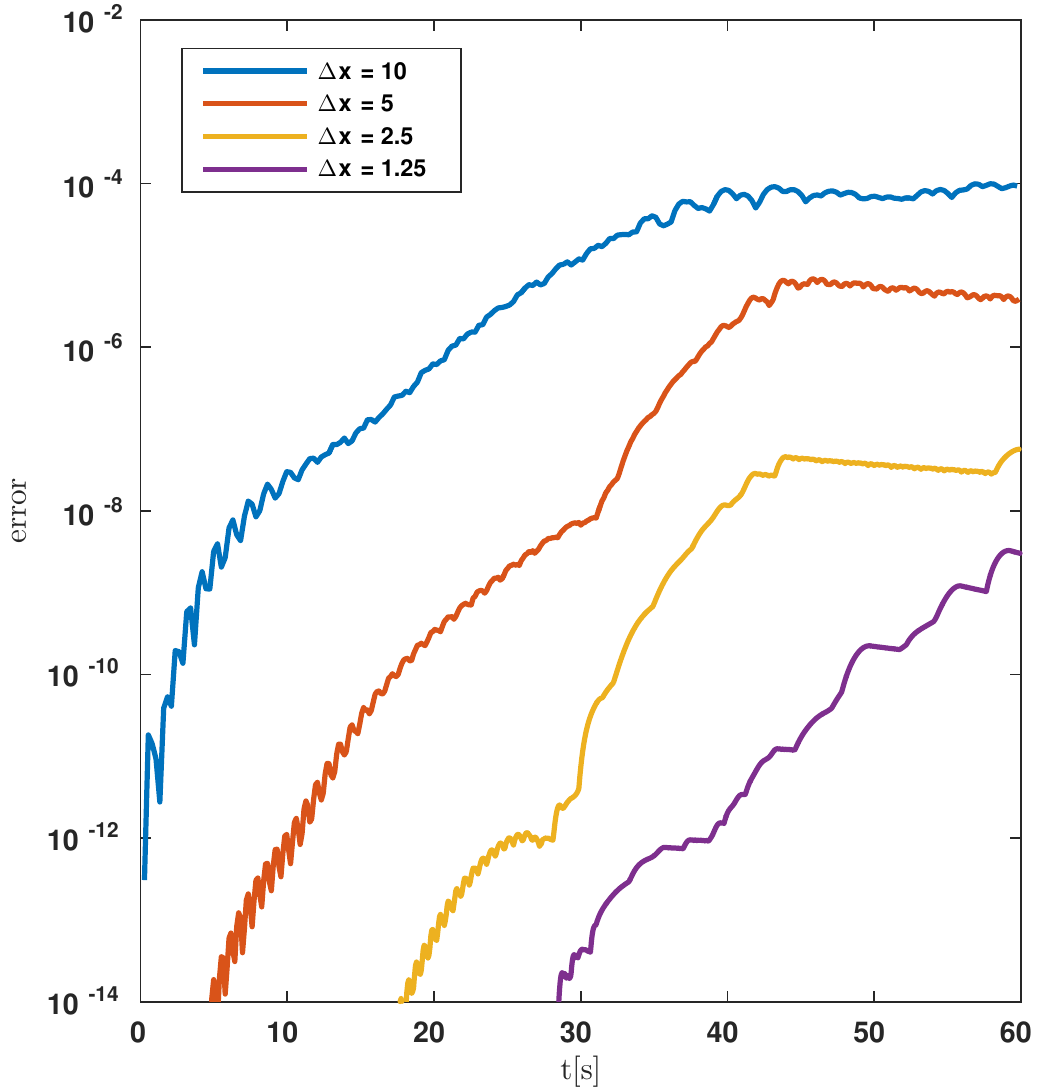}}{GL.}%
\hspace{0.0cm}%
\stackunder[5pt]{\includegraphics[width=0.33\linewidth]{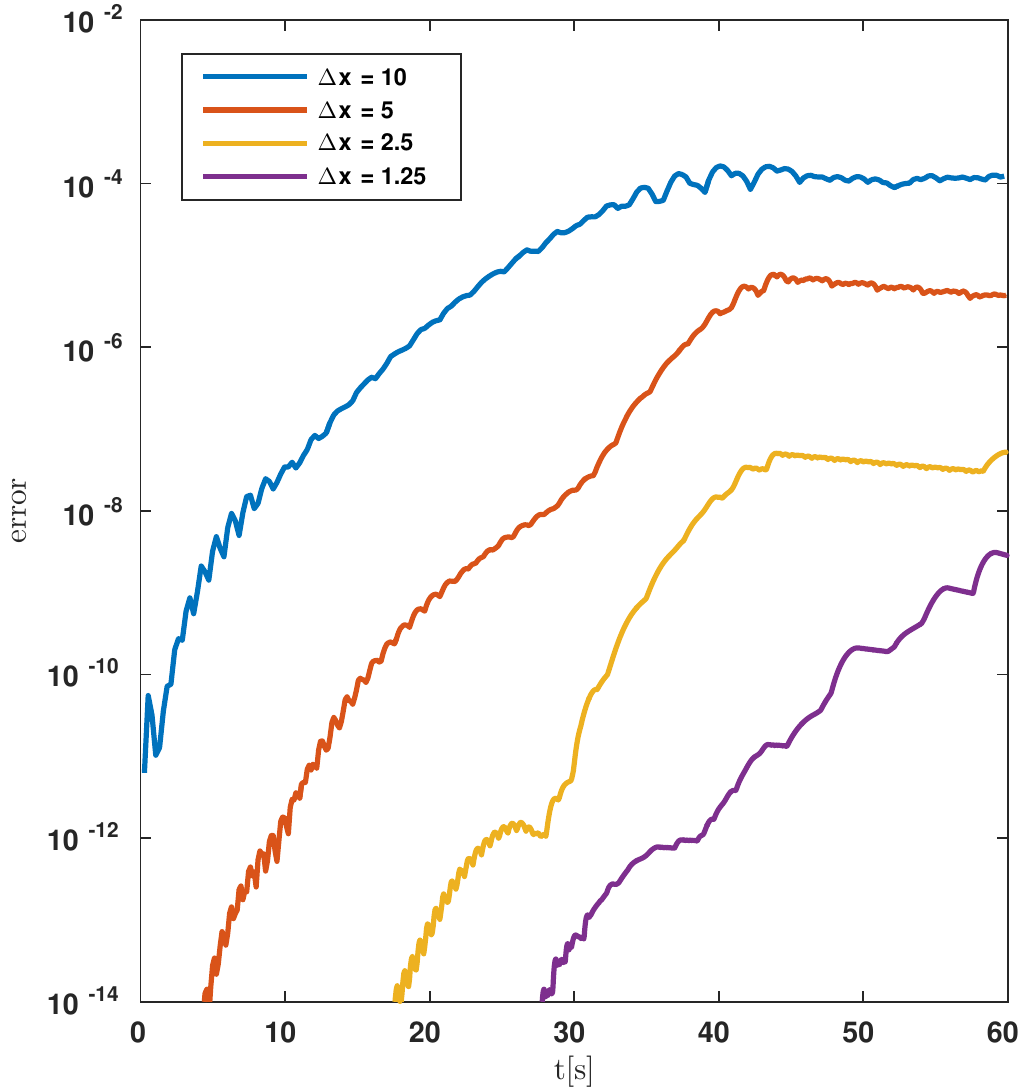}}{GLR.}%
%wave_field_pml_theta1_t100s.eps
%\stackunder[5pt]{\includegraphics[width=0.49\textwidth]{wave_field_t100s.eps}}{$\theta_x = 1$.}%
     \end{subfigure}
    \caption{Time history of PML errors using $P=4$ degree polynomial and element sizes $\Delta{x} = \{ 10, 5, 2.5, 1.25\}$. }
    \label{fig:Time_series_error_strip_h_convergence}
\end{figure}
%%%%%%%%%%%%%%%%%%%
%%%%%%%%%%%%%%%%%%%

\subsection{Whole space problem}
\begin{figure}[h!]
\begin{subfigure}
    \centering
        %%%%%%%
\stackunder[5pt]{\includegraphics[width=0.325\linewidth]{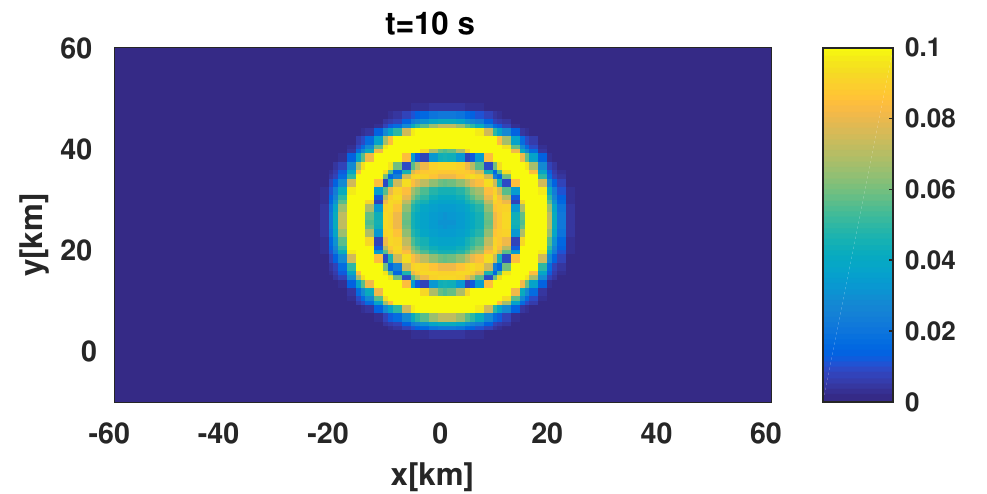}}{}%
\hspace{0.0cm}%
\stackunder[5pt]{\includegraphics[width=0.325\linewidth]{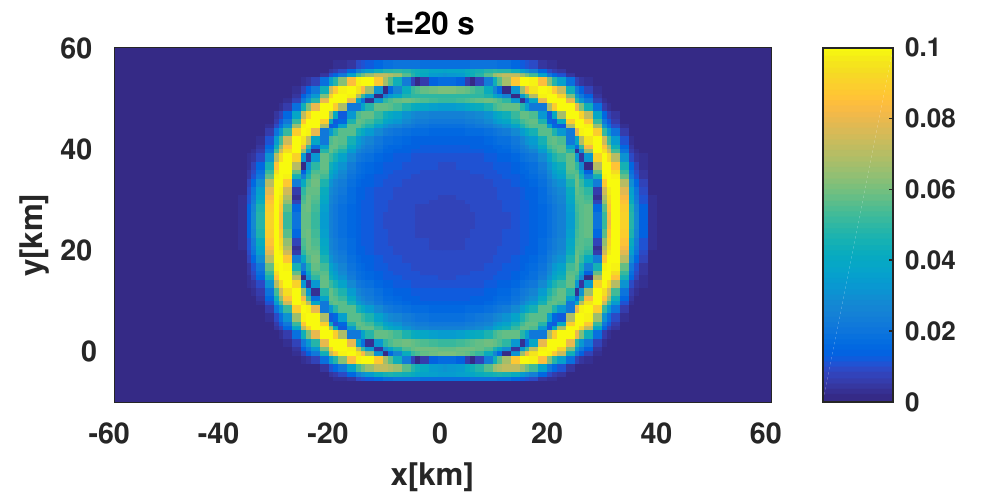}}{}%
\hspace{0.0cm}%
\stackunder[5pt]{\includegraphics[width=0.325\linewidth]{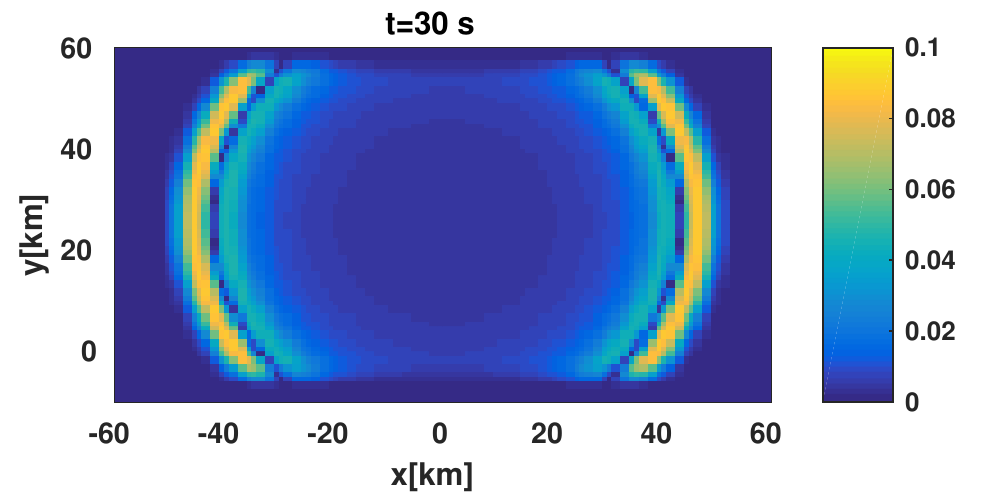}}{}%
\\
%%%%
\stackunder[5pt]{\includegraphics[width=0.325\linewidth]{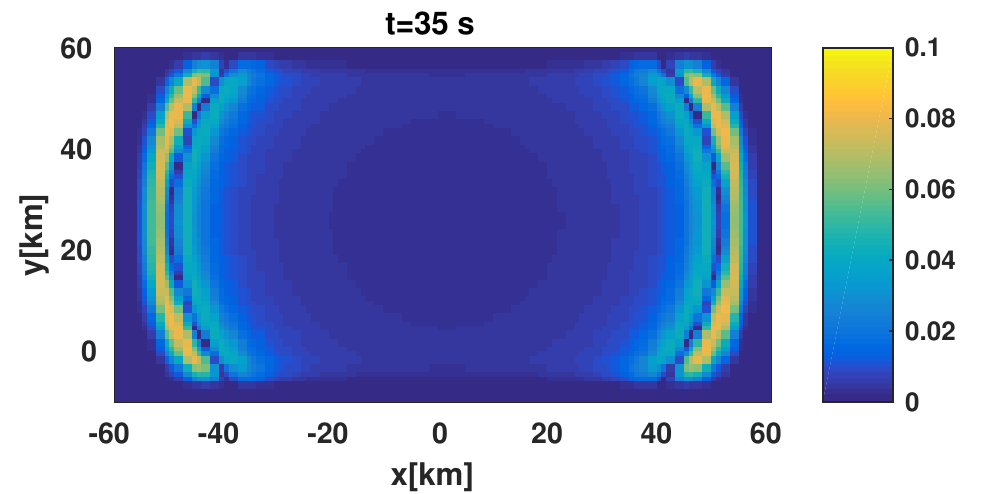}}{}%
\hspace{0.0cm}%
\stackunder[5pt]{\includegraphics[width=0.325\linewidth]{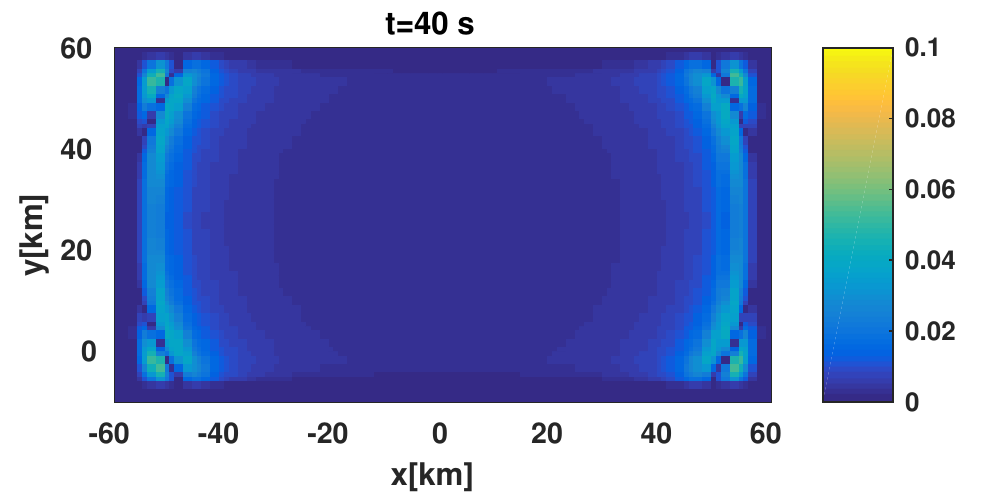}}{}%
\hspace{0.0cm}%
\stackunder[5pt]{\includegraphics[width=0.325\linewidth]{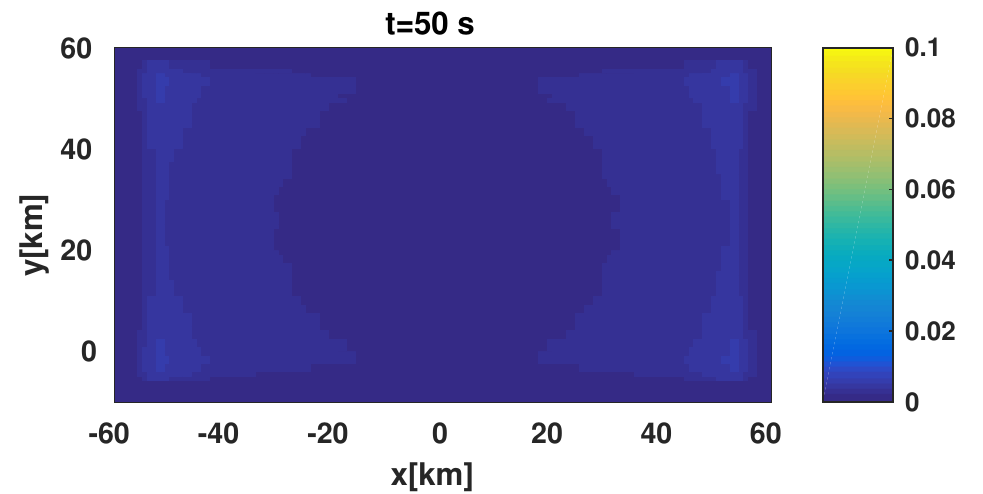}}{}%
     \end{subfigure}
    \caption{Snapshots of the absolute pressure $|p|$ in a 2D rectangular domain completely surrounded by the  PML, with PML flux fluctuation stabilization, $\omega_y = 1$. Note that   there are  corner regions where both PML damping functions,  $d_x\left(x\right)>0$, $d_y\left(y\right)>0$,  are simultaneously active.}
    \label{fig:Omega_1_whole_space}
\end{figure}

\begin{figure} [htb!]
 \centering
%%%%%%%%%%%%%%%%%%
{\includegraphics[width=0.5\linewidth]{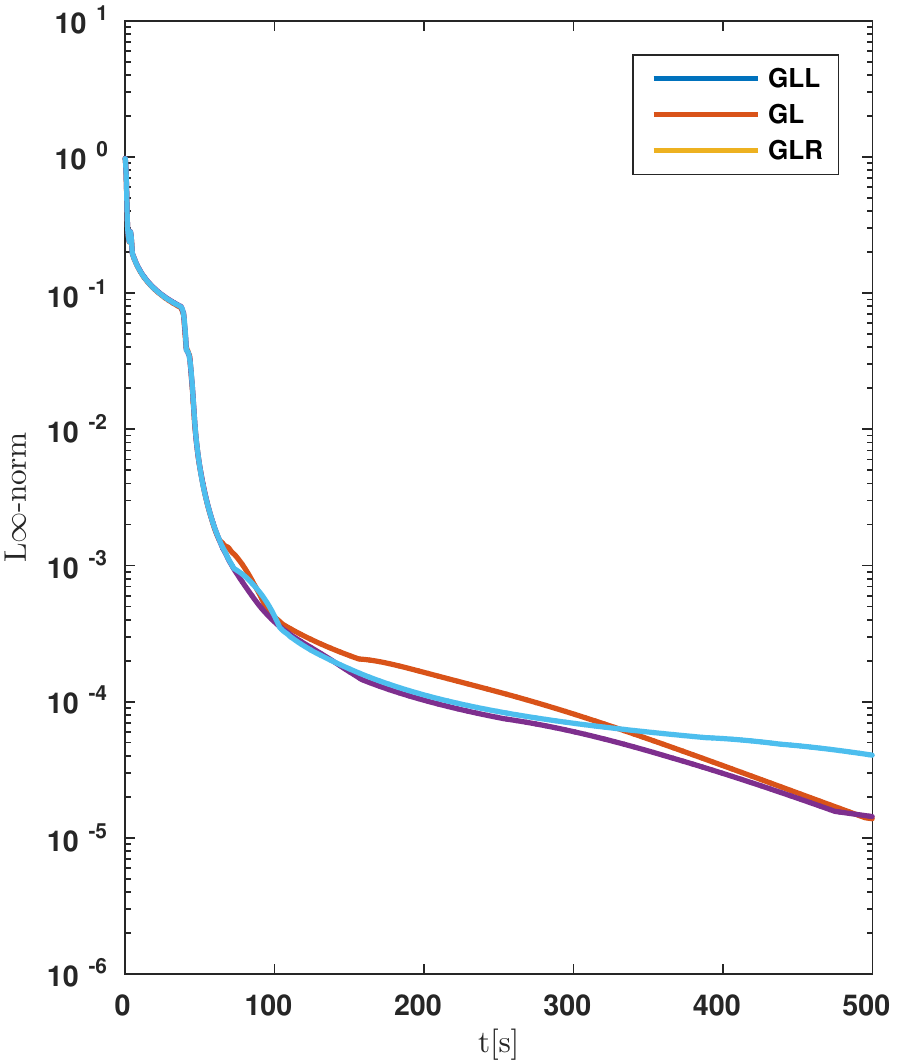}\label{fig:EnergyVsTime_wholespace}}
%%%%%%%%%%%%%%%%%%
%%%%%%%%%%%%%%%%%%%
%%%%%%%%%%%%%%%%%%%
 \caption{\textit{Time histories of the L$_{\infty}$-norm of the pressure field in a 2D rectangular domain completely surrounded by the  PML, with $\omega_y = 1$.}}
 \label{fig:time_history_whole_space}
\end{figure}

%%%%%%%%%%%%%%%%%%%
We will simulate a whole space problem, that is a computational domain surrounded completely by the PML. We consider the rectangular domain above and include the PML of width $\delta = 10$ km in both $x$-direction and $y$-direction. Thus there is a vertical layer defined by $d_x(x) > 0, d_y(y) = 0$ and a horizontal layer defined by $d_x(x) = 0, d_y(y) > 0$. Note that there are corner regions where both the vertical layer and the horizontal layer are simultaneously active, that is $ d_x(x) > 0, d_y(y) > 0$. Here we will investigate again numerical stability. From the analysis of previous sections and the numerical experiments of the last subsection we know that the PML stabilizing parameter $\omega = 1$ is critical for numerical stability. Therefore, we set $\omega = 1$ and run the simulation for a long time, $t = 500$ s. We have used  the mesh size $\Delta{x} = 5$ km, giving two DGSEM elements for the PML. The three  quadrature rules, GLL, GL and GLR, are used separately, and they give equivalent results. Snapshots of the absolute pressure field $|p|$ are displayed in Figure \ref{fig:Omega_1_whole_space}, showing how the initial Gaussian pulse spreads, and its being absorbed by  the PML.  At $t = 50$ s, the waves have been completely absorbed by the PML. The time history of the L$_{\infty}$-norm of the pressure field is plotted in  Figure \ref{fig:time_history_whole_space}. The solution decays through out the simulation. This again verifies the stability analysis of the previous sections. Note that  our PML and the numerical approximations do not suffer numerical stiffness of corners, as demonstrated in \cite{BecachePrieto2012} for some PML models.
%%%%%%%%%%%%%%%%%%%
%%%%%%%%%%%%%%%%%%%
%%%%%%%%%%%%%%%%%%%

%%%%%%%%%%%%%%%%%%%%
%%%%%%%%%%%%%%%%%%%%
\subsection{3D acoustics wave with the PML}
%%%%%%%%%%%%%%%%%%%
Here, we present 3D numerical experiments. Consider the 3D acoustic wave  equation and simulate an unbounded domain with  the PML.
%%%%%%%%%%%%%%%%%%%
We generate waves by  adding  the point source
%%%%%%%%%%%%%%%%%%%
\[
f(x,y,z,t) =   \delta_x(x-x_0)\delta_y(y-y_0)\delta_z(z-z_0)g(t), \quad g(t) = \frac{1}{\sigma_0\sqrt{2\pi}}e^{-\frac{(t-t_0)^2}{2\sigma_0^2}}, \quad \sigma_0 = 0.1149, \quad t_0 = 0.7,
\]
%%%%%%%%%%%%%%%%%%%
to the pressure equation, where $ \delta_{\xi}(\xi)$ is the one dimensional Dirac delta function. Note that $\int_{-\infty}^{\infty} g(t) dt = 1$.
%%%%%%%%%%%%%%%%%%%
In an unbounded domain, with constant  density $\rho$ and  wave speed $c$, the pressure field has the exact solution
%%%%%%%%%%%%%%%%%%%
%%%%%%%%%%%%%%%%%%%
%%%%%%%%%%%%%%%%%%%
\begin{align}\label{eq:analytica_sol}
p(x, y, z, t) = \frac{-1}{4 \pi r c^2} g^{\prime}(t-r/c), \quad r = \sqrt{(x-x_0)^2 + (y-y_0)^2 + (z-z_0)^2}.
\end{align}
%%%%%%%%%%%%%%%%%%%
%%%%%%%%%%%%%%%%%%%

\begin{figure}[h!]
\begin{subfigure}
    \centering
        %%%%%%%
%\stackunder[5pt]{\includegraphics[width=0.49\textwidth]{wave_field_pml_t100s.eps}}{$\theta_x = 0$.}%
\stackunder[5pt]{\includegraphics[width=0.35\linewidth]{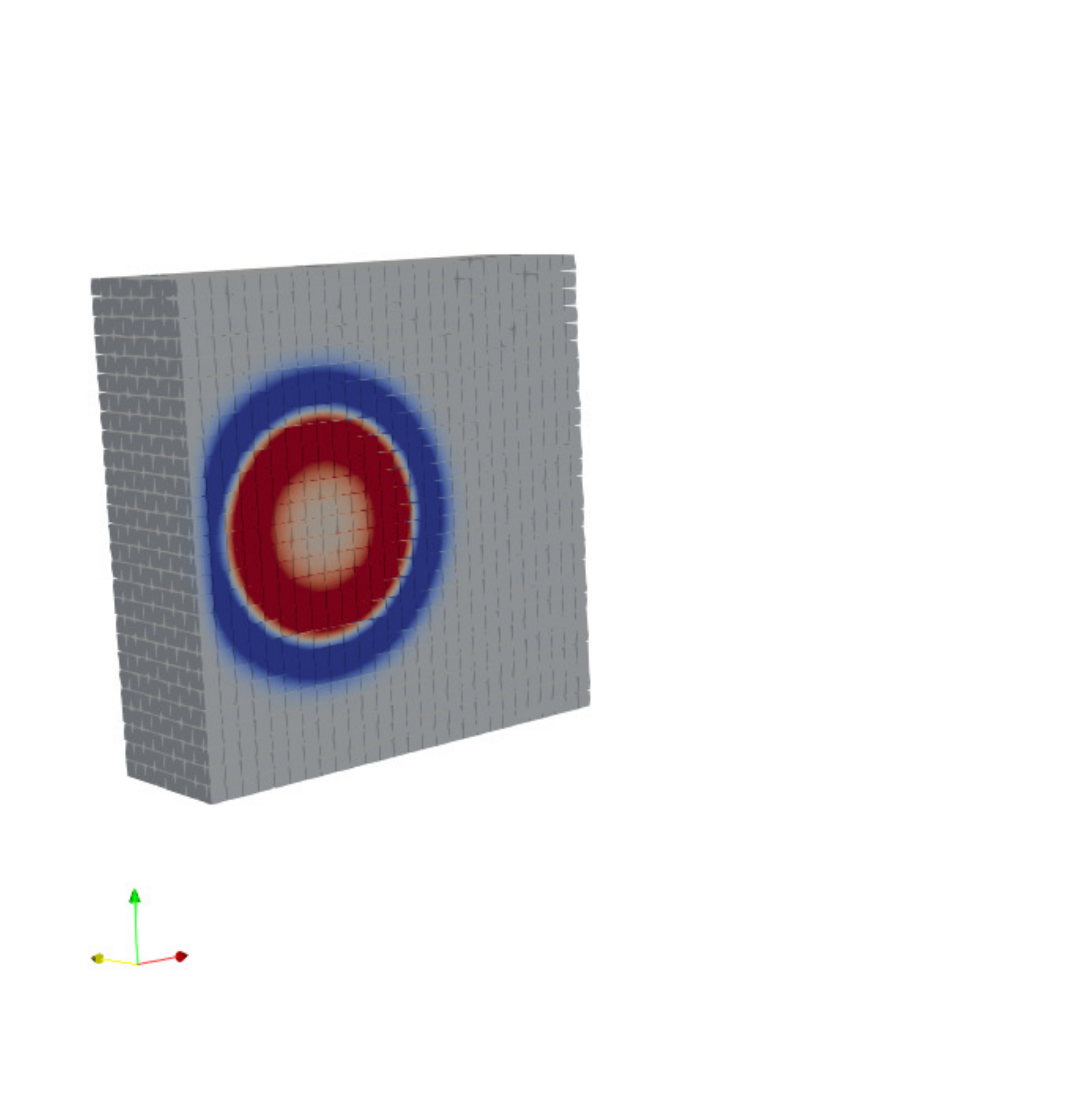}}{$t = 1.7$ s}%
\hspace{0.0cm}%
\stackunder[5pt]{\includegraphics[width=0.35\linewidth]{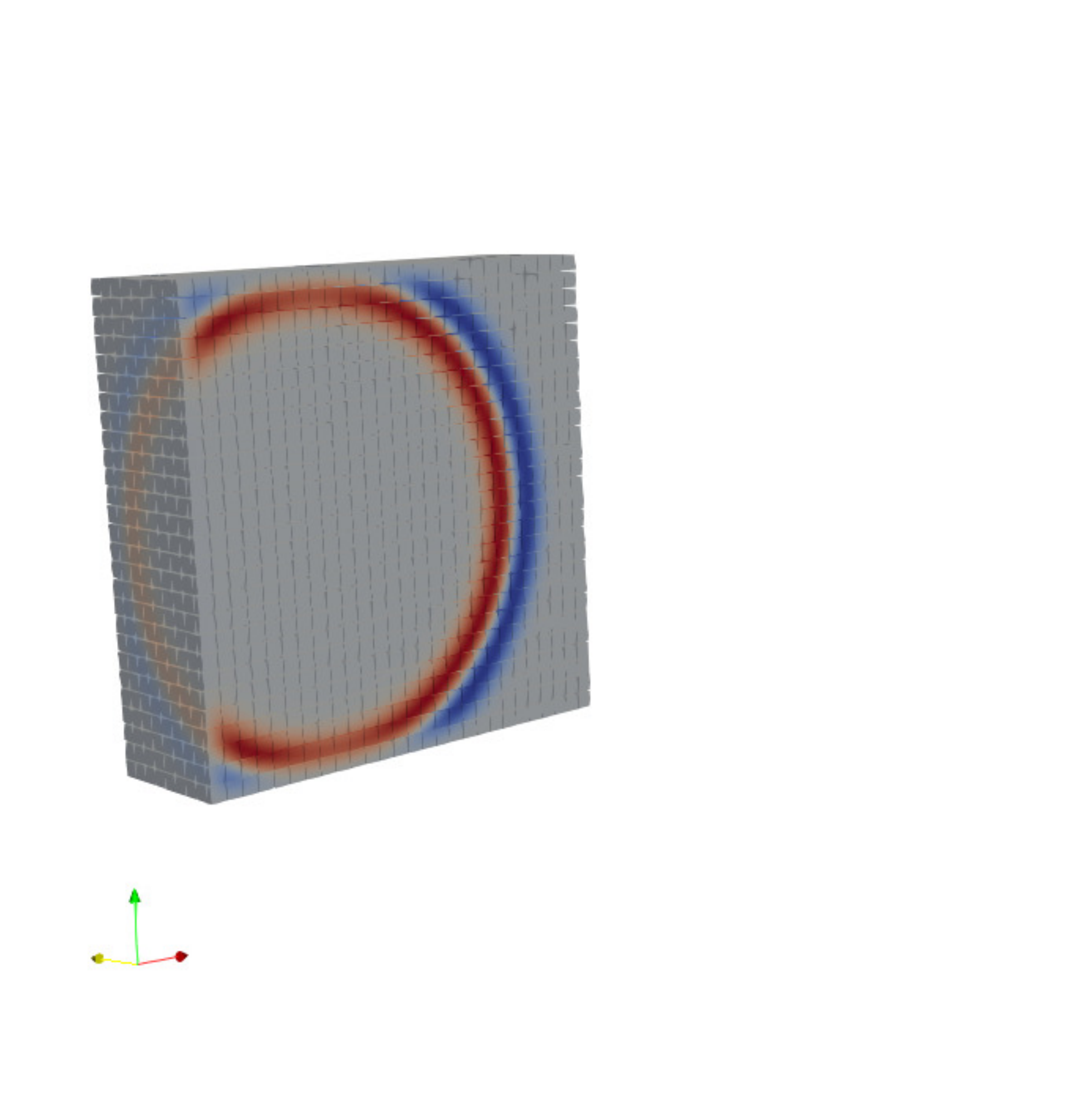}}{$t = 2.5$ s}%
\hspace{0.0cm}%
\stackunder[5pt]{\includegraphics[width=0.35\linewidth]{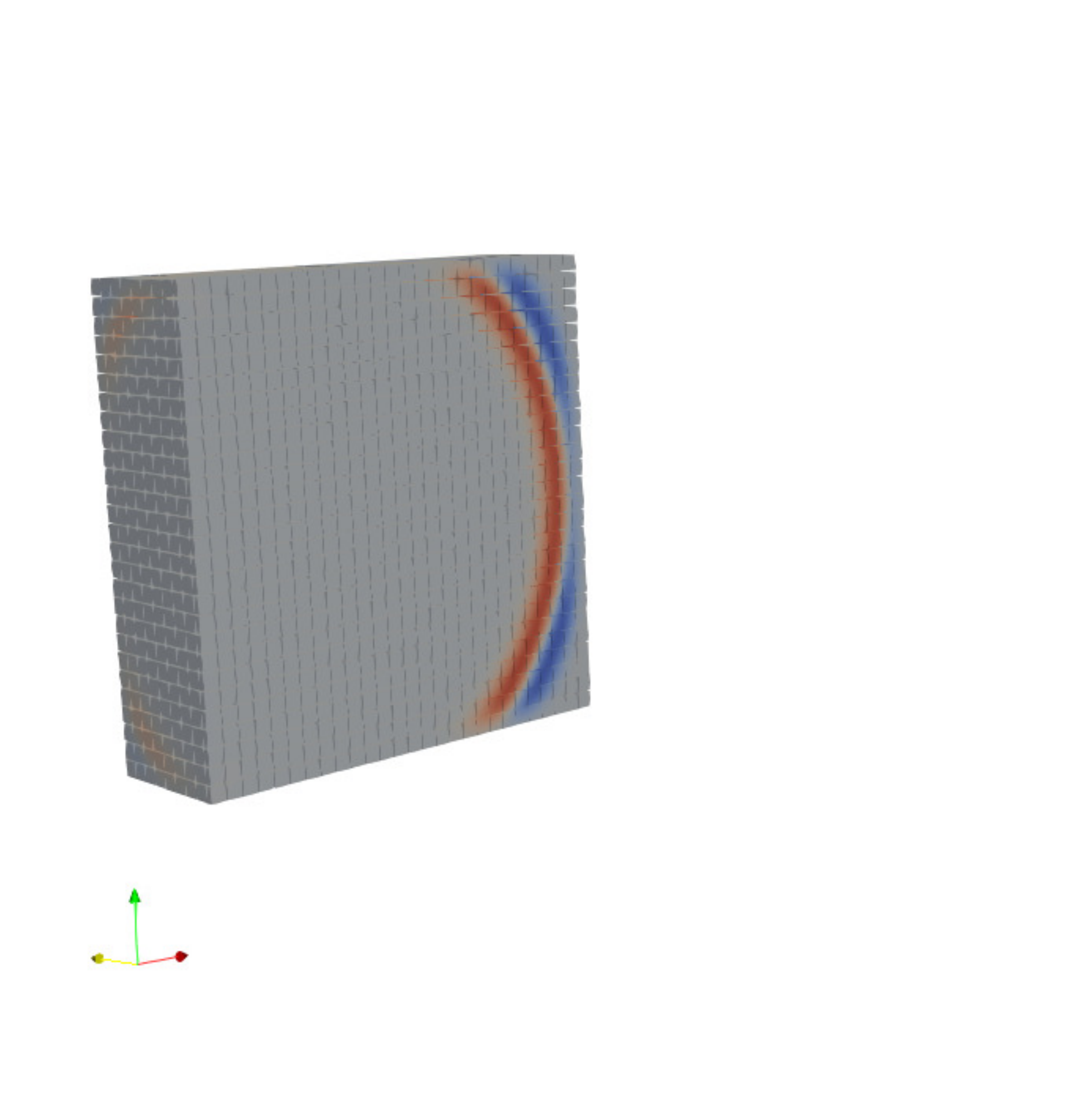}}{$t = 3.0$ s}%
%wave_field_pml_theta1_t100s.eps
%\stackunder[5pt]{\includegraphics[width=0.49\textwidth]{wave_field_t100s.eps}}{$\theta_x = 1$.}%
     \end{subfigure}
    \caption{Snapshots of the pressure $p$ in a 3D cubic domain completely surrounded by the  PML, with PML flux fluctuation stabilization, $\omega_y = 1$, $\omega_z = 1$.}
    \label{fig:Omega_1_whole_space_3D}
\end{figure}

\begin{figure}[h!]
\begin{subfigure}
    \centering
        %%%%%%%
%\stackunder[5pt]{\includegraphics[width=0.49\textwidth]{wave_field_pml_t100s.eps}}{$\theta_x = 0$.}%
\stackunder[5pt]{\includegraphics[width=0.5\linewidth]{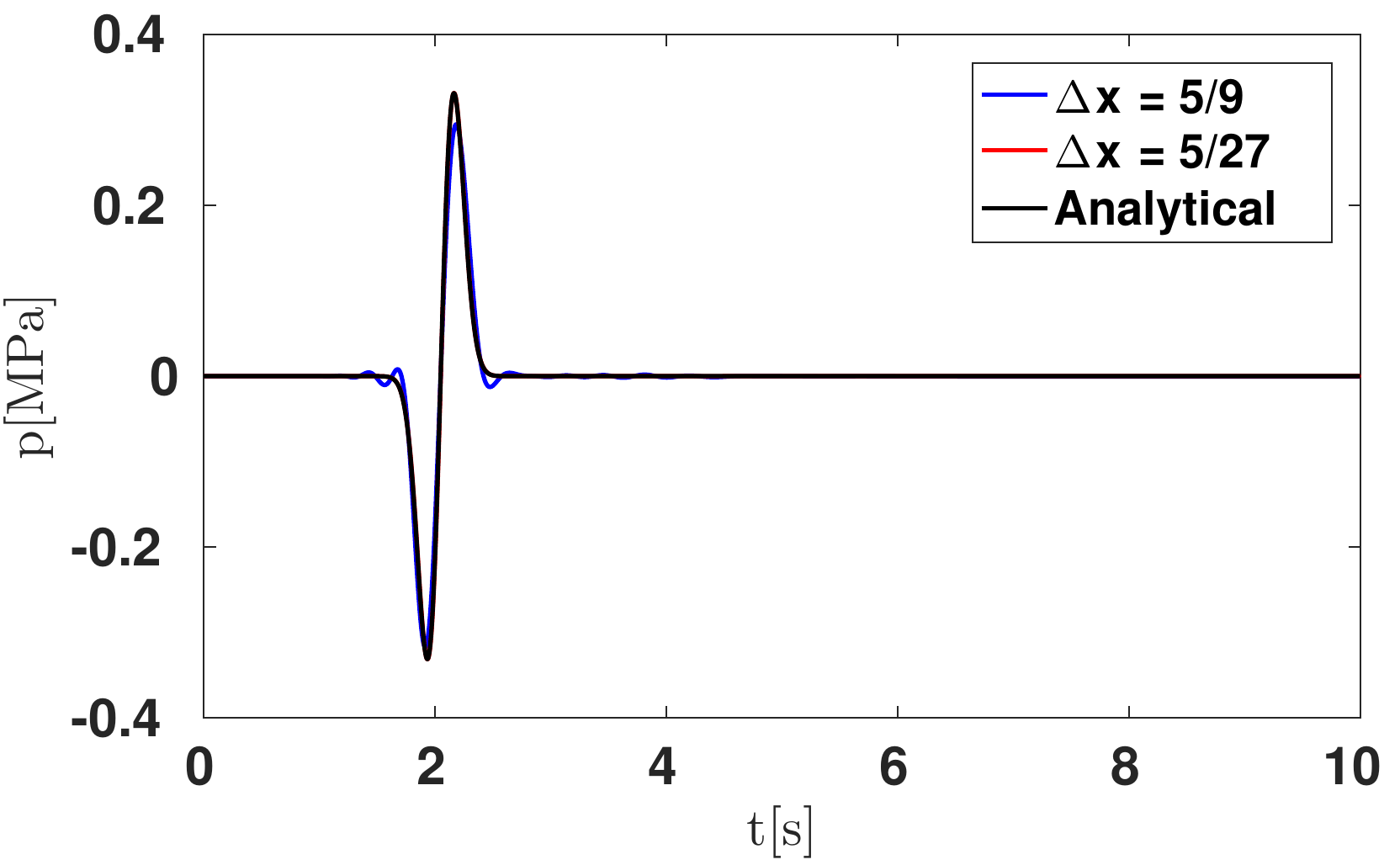}}{PML}%
\hspace{0.0cm}%
\stackunder[5pt]{\includegraphics[width=0.5\linewidth]{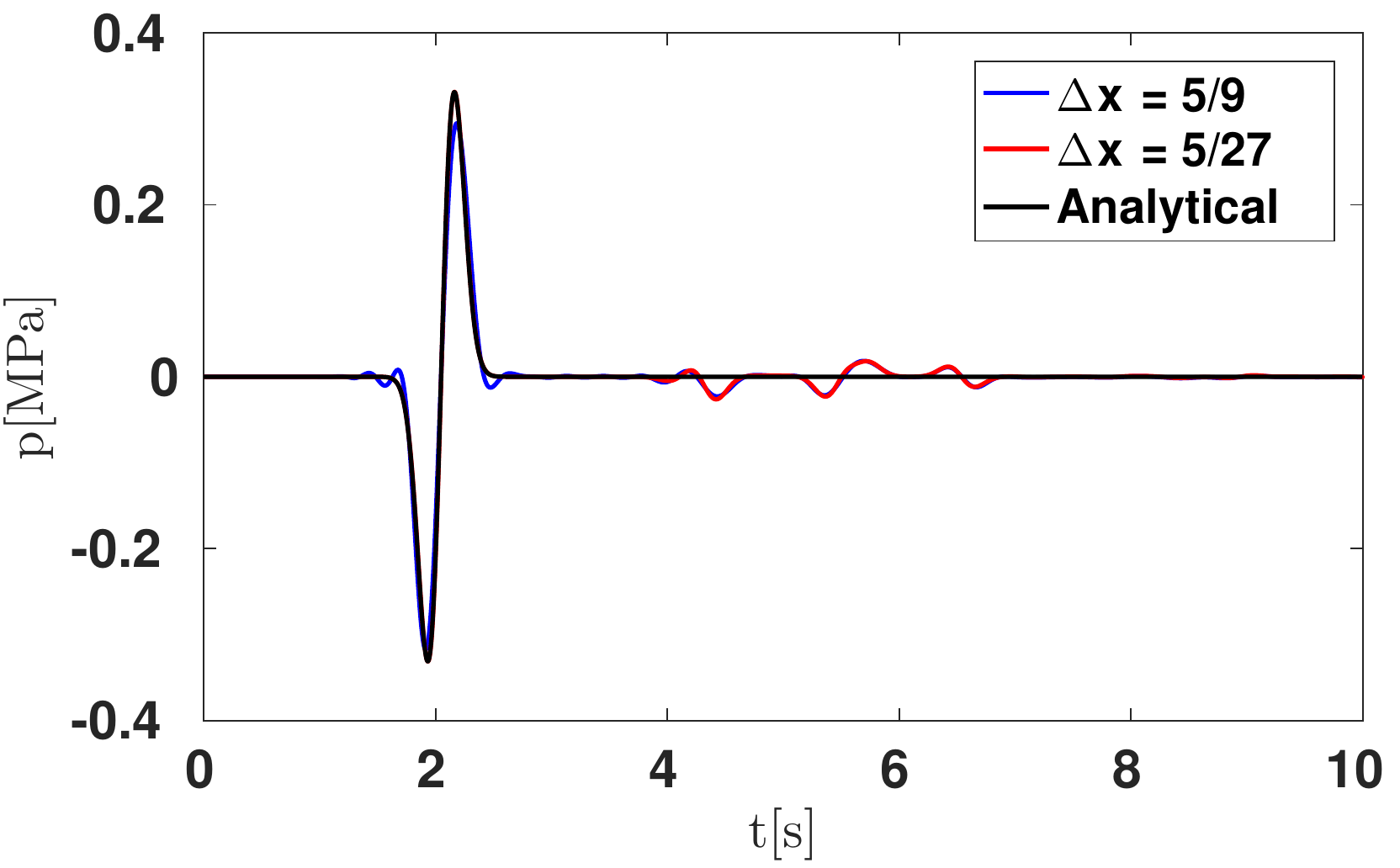}}{Absorbing boundary condition}%
%\hspace{0.0cm}%
%\stackunder[5pt]{\includegraphics[width=0.35\linewidth]{snap_shot_t30.eps}}{$t = 3.0$ s}%
%wave_field_pml_theta1_t100s.eps
%\stackunder[5pt]{\includegraphics[width=0.49\textwidth]{wave_field_t100s.eps}}{$\theta_x = 1$.}%
     \end{subfigure}
   \caption{A comparison of the  pressure field against the analytical solution  at  $(x_r, y_r, z_r) = (3.5, 2.5, 2.5)$. Note without the PML artificial reflections become prominent after $t \ge 3.7$ s. }
    \label{fig:Seismogram_Omega_1_whole_space_3D}
\end{figure}

We consider the computational cube, $(x, y, z) \in [0$ km, $5$ km$]$ $\times$ $[0$ km, $5$ km$]$ $\times$ $[0$ km, $5$ km$]$, and place the source at $(x_0, y_0, z_0) = (1.5, 2.5, 2.5)$ and a receiver, $2$ km away  from the source, at  $(x_r, y_r, z_r) = (3.5, 2.5, 2.5)$. The boundaries are closed with the absorbing boundary condition. That is $r_x = r_y = r_z = 0$ in \eqref{eq:boundary_condition_acoustic}. We discretize the domain with two levels of mesh resolution, $9\times9\times9$  and $27\times27\times27$ DGSEM elements, separately,  with $P =4$ polynomial approximation on GL nodes. The PML width is $\delta = 0.550 $ km.  At the initial mesh resolution $9\times9\times9$ the PML is contained in a single element, at the higher mesh resolution $27\times27\times27$ the PML  spans three DGSEM elements at each domain boundary.  From the analysis in the previous sections, and 2D numerical experiments presented in the last subsection, we know that the PML flux fluctuations stabilization parameters $\omega_y = 1$, $\omega_z = 1$ are critical for numerical stability of the PML. Therefore, we set $\omega_y = 1$, $\omega_z = 1$ and use the relative PML error, $tol = 0.1\%$.  We run the simulation for $t = 10$ s. Snapshots of the pressure field are displayed in Figure \ref{fig:Omega_1_whole_space_3D}, showing how the pressure spreads and the absorption of waves  by the PML. We have also performed numerical simulations with only the absorbing boundary condition, that is with zero PML parameters, $d_x=d_y=d_z =0$. In Figure \ref{fig:Seismogram_Omega_1_whole_space_3D}, we compare the seismograms at  $(x_r, y_r, z_r) = (3.5, 2.5, 2.5)$ with the analytical solution \eqref{eq:analytica_sol}. With the absorbing boundary condition only, the numerical solution matches the analytical solution very well before reflections from domain boundaries arrive.  At about $t \ge 3.7$ s, the numerical reflections have corrupted the solution everywhere. The solutions can never converge by p- or h-refinement, see Figure \ref{fig:Seismogram_Omega_1_whole_space_3D}. With PML the numerical solution matches the analytical solution excellently for all times. We can make very accurate numerical simulation at any future time.  As demonstrated in the last subsection, see also Figure \ref{fig:Seismogram_Omega_1_whole_space_3D}, by the appropriate choice of PML parameters we can make all errors converge to zero by the convergence rate of the DGSEM approximation.

%%%%%%%%%%%%%%%%%%%

%\begin{figure}[h!]
%\begin{subfigure}
%    \centering
%        %%%%%%%
%         {\includegraphics[width=0.5\linewidth]{PML_Seismogram_3D.eps}\label{fig:EnergyVsTime_wholespace}}
%         {\includegraphics[width=0.5\linewidth]{Seismogram_3D.eps}\label{fig:EnergyVsTime_wholespace}}
%        %%%%%%%
%     \end{subfigure}
%    \caption{A comparison of the  pressure field against the analytical solution  at  $(x_r, y_r, z_r) = (3.5, 2.5, 2.5)$. Note without the PML artificial reflections become prominent after $t \ge 3.7$ s }
%    \label{fig:Seismogram_Omega_1_whole_space_3D}
%\end{figure}
%%%%%%%%%%%%%%%%%%%
%%%%%%%%%%%%%%%%%%%%

%\begin{figure}[h!]
%\begin{subfigure}
%    \centering
%        %%%%%%%
%         {\includegraphics[width=0.35\linewidth]{snap_shot_t17.eps}\label{fig:EnergyVsTime_wholespace}}
%          {\includegraphics[width=0.35\linewidth]{snap_shot_t25.eps}\label{fig:EnergyVsTime_wholespace}}
%           {\includegraphics[width=0.35\linewidth]{snap_shot_t30.eps}\label{fig:EnergyVsTime_wholespace}}       
%     \end{subfigure}
%    \caption{Time history of PML errors using $P=4$ degree polynomial and element sizes $\Delta{x} = \{ 10, 5, 2.5, 1.25\}$. }
%    \label{fig:Time_series_error_strip_h_convergence}
%\end{figure}
%
%%%%%%%%%%%%%%%%%%%%
%
%%%%%%%%%%%%%%%%%%%%
%%%%%%%%%%%%%%%%%%%%
\section{Summary and outlook}
%%%%%%%%%%%%%%%%%%%%
The well-posedness and stability of  the PML for the acoustic wave equation have been established in previous works, see  \cite{DuKrSIAM, HalpernPetit-BergezRauch2011}. However, previous attempts to effectively include the PML in many modern numerical methods such the DGSEM proved to be a nightmare for practitioners. Exponential and/or linear growth is often seen in numerical simulations. In this paper,  we have developed a provably  DGSEM for  the acoustic wave equation truncated by the PML.  Our approach is rooted in a rigorous mathematical analysis, beginning from the continuous model down to the discrete problem. We begin by deriving  continuous energy estimates for the 3D PML in the Laplace space.  By mimicking the energy estimate in the discrete setting we construct asymptotically stable DGSEM approximation of the PML for the acoustic wave equation.  The analysis here is performed for the 3D PML problem. However, the results are also valid for the corresponding 2D PML model. We performed numerical experiments, in 2D and 3D, verifying stability and accuracy of the PML.

%%%%%%%%%%%%%%%%%%%%
%%%%%%%%%%%%%%%%%%%%
%%%%%%%%%%%%%%%%%%%%
The 2D and 3D production code, for acoustic waves simulation, have been  implemented on ExaHyPE (www.exahype.eu), a simulation engine for hyperbolic PDEs, on adaptive Cartesian meshes, for exa-scale supercomputers. This software, ExaHyPE,  is open source and publicly available.
%%%%%%%%%%%%%%%%%%%%
%%%%%%%%%%%%%%%%%%%%

%
%
The theory and techniques developed  in this paper can be extended to other problems such as the Maxwell's equations and linear elastic wave equation. However,  for certain problems such as elastic wave propagation, where more than one wave type and wave speed are simultaneously present, the current  theory may not be extended in straightforward manner. Therefore, we will require further assumptions and simplifications for such problems. In a forthcoming paper \cite{DuruGabrielKreiss2017}  we will extend the theory developed in this paper to more complex problems such as linear elastic wave.
%%%%%%%%%%%%%%%%%%%%%
%
\section*{Acknowledegments}
The project has received funding from the European Union's Horizon 2020 research and innovation program under grant agreement No 671698.
%\newline
%{}\hfill{\includegraphics[angle=90, width=0.15\textwidth]{EU_logo.pdf}}

%%%%%%%%%%%%%%%%%
%%%%%%%%%%%%%%%%%
%%%%%%%%%%%%%%%%%%%%
%
\appendix
%%%%%%%%%%%%%%%%%%%%%
\section{Proof of Theorem \ref{Theo:Stability_PML_Laplace}}\label{sec:appendix1}
%%%%%%%%%%%
Eliminate $\sigma$, $\psi$ in \eqref{eq:acoustic_pml_3D_Laplace1} by using ${\phi}_{\sigma} = {\phi}_p/(s^*S_y^*)$, ${\phi}_{\psi} = {\phi}_p/(s^*S_z^*)$ and
%%%%%%%%%%%
%%%%%%%%%%% 
\[
 \int_{\Omega}{{\phi}_p^*}\widetilde{\sigma}dxdydz = -\int_{\Omega}\frac{{\phi}_p^*}{sS_y}\left(\left(d_y - d_x\right)\frac{\partial \widetilde{v}}{\partial y} - f_{\sigma}(x,y,z)\right)dxdydz,
\]
%%%%%%%%%%%
%\quad 
%%%%%%%%%%%
\[
 \int_{\Omega}{{\phi}_p^*}\widetilde{\psi}dxdydz = -\int_{\Omega}\frac{{\phi}_p^*}{sS_y}\left(\left(d_z - d_x\right)\frac{\partial \widetilde{w}}{\partial z} - f_{\psi}(x,y,z) \right)dxdydz.
\]
%%%%%%%%%%%
%%%%%%%%%%%
Introducing ${\phi}_p = \widetilde{\phi}_p/S_x^*$, ${\phi}_u = \widetilde{\phi}_u/(\rho S_x^*)$, ${\phi}_v = \widetilde{\phi}_v/(\rho S_y^*)$, ${\phi}_w = \widetilde{\phi}_w/(\rho S_z^*)$ we have
%%%%%%%%%%%
%%%%%%%%%%%
\begin{equation}\label{eq:elastic_pml_2D_Laplace1_appendix}
\begin{split}
&\int_{\Omega}\widetilde{\phi}_p^*\frac{s}{\kappa} \widetilde{p} dxdydz = - \int_{\Omega} \widetilde{\phi}_p^*\left(\frac{1}{S_x}\frac{\partial \widetilde{u}}{\partial x} +  \frac{1}{S_y}\frac{\partial \widetilde{v}}{\partial y} + \frac{1}{S_z}\frac{\partial \widetilde{w}}{\partial z} \right) dxdydz\\
& + \int_{\Omega} \widetilde{\phi}_p^*\left( \frac{1}{S_x}f_{p}(x,y,z) - \frac{1}{sS_yS_x}f_{\sigma}(x,y,z)  -  \frac{1}{sS_zS_x}f_{\psi}(x,y,z)  \right) dxdydz, \\
&\int_{\Omega}\widetilde{\phi}_u^*s\widetilde{u} dxdydz = -\int_{\Omega}\frac{\widetilde{\phi}_u^*}{\rho S_x}\left(\frac{\partial \widetilde{p}}{\partial x} - \rho f_u(x,y,z)\right) dxdydz, \\
&\int_{\Omega}\widetilde{\phi}_v^*s\widetilde{v} dxdydz = -\int_{\Omega}\frac{\widetilde{\phi}_v^*}{\rho S_y}\left(\frac{\partial \widetilde{p}}{\partial y} -  \rho f_v(x,y,z)\right) dxdydz,\\
&\int_{\Omega}\widetilde{\phi}_w^*s\widetilde{w} dxdydz = -\int_{\Omega}\frac{\widetilde{\phi}_w^*}{\rho S_z}\left(\frac{\partial \widetilde{p}}{\partial z} -\rho  f_w(x,y,z)\right)dxdydz.
 %%
%\frac{1}{\kappa} s \widetilde{p}  &= -\frac{1}{S_x}\frac{\partial \widetilde{u}}{\partial x} -  \frac{1}{S_y}\frac{\partial \widetilde{v} }{\partial y} + \frac{1}{\kappa}f_p(x,y), \\
%\rho{ s\widetilde{u}} &= -\frac{1}{S_x}\frac{\partial \widetilde{p}}{\partial x} + \rho f_u(x,y), \\
%\rho s{\widetilde{v}} &= -\frac{1}{S_y}\frac{\partial \widetilde{p}}{\partial y} + \rho f_v(x,y),
%\mathbf{P} s\widetilde{\mathbf{u}} = \frac{1}{S_x} \mathbf{A}\frac{\partial{\widetilde{\mathbf{u}} }}{\partial x} + \frac{1}{S_y} \mathbf{B}\frac{\partial{\widetilde{\mathbf{u}} }}{\partial y} + \mathbf{P}{\mathbf{u}}^0, \quad \Re{s} > 0,
  \end{split}
  \end{equation}
  Integrating the first equation in \eqref{eq:elastic_pml_2D_Laplace1_appendix} by parts, using 
  \[
   \left(\widetilde{\phi}_p, \widetilde{\phi}_u, \widetilde{\phi}_v, \widetilde{\phi}_w\right) = \left(s^*s\widetilde{p}, \frac{s}{S_x^*}\frac{\partial \widetilde{p}}{\partial x}, \frac{s}{S_y^*}\frac{\partial \widetilde{p}}{\partial y}, \frac{s}{S_z^*}\frac{\partial \widetilde{p}}{\partial z}\right),
  \]
   and eliminating the velocity equations,
   we have

    \begin{equation}\label{eq:acoustic_pml_2D_Laplace1}
\begin{split}
  &s\int_{\Omega} (s\widetilde{p})^*\frac{1}{\kappa} (s\widetilde{p}) dxdydz +  \frac{s^* S_x^*}{S_x}\int_{\Omega}\frac{1}{S_x^*} \frac{\partial \widetilde{p}^*}{\partial x}  \frac{1}{\rho}  \frac{1} {S_x}\frac{\partial \widetilde{p}}{\partial x} dxdydz +  \frac{s^* S_y^*}{S_y}\int_{\Omega} \frac{1}{S_y^*} \frac{\partial \widetilde{p}^*}{\partial y}   \frac{1}{\rho} \frac{1}{S_y} \frac{\partial \widetilde{p}}{\partial y} dxdydz \\
  &+  \frac{s^* S_z^*}{S_z}\int_{\Omega} \frac{1}{S_z^*} \frac{\partial \widetilde{p}^*}{\partial z}  \frac{1}{\rho}   \frac{1}{S_z} \frac{\partial \widetilde{p}}{\partial z} dxdydz 
    + \int_{-y_b}^{y_b} \int_{-z_b}^{z_b}\frac{\widetilde{u}}{S_x}{  s^*s \widetilde{p}^*}\Big|_{-x_b}^{x_b} dydz + \int_{-x_b}^{x_b} \int_{-z_b}^{z_b}\frac{\widetilde{v}}{S_y}{  s^* s \widetilde{p}^*}\Big|_{-y_b}^{y_b} dxdz \\
    &+ \int_{-x_b}^{x_b} \int_{-y_b}^{y_b}\frac{\widetilde{w}}{S_z}{  s^* s \widetilde{p}^*}\Big|_{-z_b}^{z_b} dxdy =  \\
    & + s\int_{\Omega} \left(s\widetilde{p}\right)^*\frac{1}{\kappa}\left( \frac{1}{S_x}f_{p}(x,y,z) - \frac{\kappa}{sS_yS_x}f_{\sigma}(x,y,z)  -  \frac{\kappa}{sS_zS_x}f_{\psi}(x,y,z)  \right) dxdydz \\
    &+
   \frac{s^* S_x^*}{S_x} \int_{\Omega}\frac{1}{\rho} \frac{1}{S_x^*} \frac{\partial \widetilde{p}^*}{\partial x}\frac{1}{S_z}\left( \rho  f_u(x,y,z)\right)  dxdydz 
    + \frac{s^* S_y^*}{S_y}\int_{\Omega}\frac{1}{\rho} \frac{1}{S_y^*} \frac{\partial \widetilde{p}^*}{\partial y}\frac{1}{S_y}\left( \rho  f_v(x,y,z)\right)  dxdydz \\
    &+ \frac{s^* S_z^*}{S_z}\int_{\Omega}\frac{1}{\rho} \frac{1}{S_z^*} \frac{\partial \widetilde{p}^*}{\partial z}\frac{1}{S_z}\left( \rho  f_w(x,y,z)\right)  dxdydz.
  \end{split}
  \end{equation}
  %%%%%%%%%%%%%%%%%%%
  %%%%%%%%%%%%%%%%%%%
  Introducing the inner product \eqref{eq:inner_product_complex} yields
  %%%%%%%%%%%%%%%%%%%
  %%%%%%%%%%%%%%%%%%%
  \begin{equation}\label{eq:acoustic_pml_2D_2}
\begin{split}
&{s}\left({s}  \widetilde{p}, {s}  \widetilde{p}\right)_{1/\kappa}+\left(\frac{s^*S_x^*}{S_x}\right)\left(\frac{1}{S_x} \frac{\partial \widetilde{p} }{\partial x}, \frac{1}{S_x}\frac{\partial \widetilde{p}}{\partial x}\right)_{1/\rho} + \left(\frac{s^*S_y^*}{S_y}\right)\left(\frac{1}{S_y} \frac{\partial \widetilde{p} }{\partial y}, \frac{1}{S_y}\frac{\partial \widetilde{p}}{\partial y}\right)_{1/\rho}   + \Re\left(\frac{s^*S_z^*}{S_z}\right)\left(\frac{1}{S_z} \frac{\partial \widetilde{p} }{\partial z}, \frac{1}{S_z}\frac{\partial \widetilde{p}}{\partial z}\right)_{1/\rho}  \\
& + \left(\frac{1}{S_x}\right)\widetilde{BT}^{(x)} + \left(\frac{1}{S_y}\right)\widetilde{BT}^{(y)}  + \left(\frac{1}{S_z}\right)\widetilde{BT}^{(z)}  =  \\
&{s}\left({s}  \widetilde{p}, \widetilde{F}_p \right)_{1/\kappa}+\left(\frac{s^*S_x^*}{S_x}\right)\left(\frac{1}{S_x} \frac{\partial \widetilde{p} }{\partial x}, \frac{1}{S_z}\left( \rho  f_u(x,y,z)\right)\right)_{1/\rho} + \left(\frac{s^*S_y^*}{S_y}\right)\left(\frac{1}{S_y} \frac{\partial \widetilde{p} }{\partial y}, \frac{1}{S_y}\left( \rho  f_v(x,y,z)\right)\right)_{1/\rho}   \\
&+ \left(\frac{s^*S_z^*}{S_z}\right)\left(\frac{1}{S_z} \frac{\partial \widetilde{p} }{\partial z}, \frac{1}{S_z}\left( \rho  f_w(x,y,z)\right)\right)_{1/\rho} ,
%\mathbf{P} s\widetilde{\mathbf{u}} = \frac{1}{S_x} \mathbf{A}\frac{\partial{\widetilde{\mathbf{u}} }}{\partial x} + \frac{1}{S_y} \mathbf{B}\frac{\partial{\widetilde{\mathbf{u}} }}{\partial y} + \mathbf{P}{\mathbf{u}}^0, \quad \Re{s} > 0,
  \end{split}
  \end{equation}
%%%%%%%%%%%%%%%%%%%
where the forcing is
\[
\widetilde{F}_p(x,y,z) = \frac{1}{S_x}f_{p}(x,y,z) - \frac{\kappa}{sS_yS_x}f_{\sigma}(x,y,z)  -  \frac{\kappa}{sS_zS_x}f_{\psi}(x,y,z),
\]
with  the boundary terms
%%%%%%%%%%%%%%%%%%%
{\small
\[
\widetilde{BT}^{(x)}  = |s| \int_{-y_b}^{y_b} \int_{-z_b}^{z_b}{\widetilde{u}}{\widetilde{p}^*}\Big|_{-x_b}^{x_b} dydz \ge 0, \quad
%\]
%%%%%%%%%%%%%%%%%%%
%%%%%%%%%%%%%%%%%%%
%\[
\widetilde{BT}^{(y)}  =|s| \int_{-x_b}^{x_b} \int_{-z_b}^{z_b}{\widetilde{v}}{\widetilde{p}^*}\Big|_{-y_b}^{y_b} dxdz \ge 0, \quad
%\]
%%%%%%%%%%%%%%%%%%%
%%%%%%%%%%%%%%%%%%%
%\[
\widetilde{BT}^{(z)}  =|s| \int_{-x_b}^{x_b} \int_{-y_b}^{y_b}{\widetilde{w}}{\widetilde{p}^*}\Big|_{-z_b}^{z_b} dxdy \ge 0.
\]
}
%%%%%%%%%%%%%%%%%%%
%%%%%%%%%%%%%%%%%%%
 Adding the complex of conjugate of \eqref{eq:acoustic_pml_2D_2}  gives
 %%%%%%%%%%%%%%%%%%%
%%%%%%%%%%%%%%%%%%%
\begin{equation}\label{eq:elastic_pml_2D_2}
\begin{split}
&\Re{(s)}\left({s}  \widetilde{p}, {s}  \widetilde{p}\right)_{1/\kappa}+\Re\left(\frac{s^*S_x^*}{S_x}\right)\left(\frac{1}{S_x} \frac{\partial \widetilde{p} }{\partial x}, \frac{1}{S_x}\frac{\partial \widetilde{p}}{\partial x}\right)_{1/\rho} + \Re\left(\frac{s^*S_y^*}{S_y}\right)\left(\frac{1}{S_y} \frac{\partial \widetilde{p} }{\partial y}, \frac{1}{S_y}\frac{\partial \widetilde{p}}{\partial y}\right)_{1/\rho}   \\
&+ \Re\left(\frac{s^*S_z^*}{S_z}\right)\left(\frac{1}{S_z} \frac{\partial \widetilde{p} }{\partial z}, \frac{1}{S_z}\frac{\partial \widetilde{p}}{\partial z}\right)_{1/\rho}  
 + |s|\left(\Re\left(\frac{1}{S_x}\right)\widetilde{BT}^{(x)} + \Re\left(\frac{1}{S_y}\right)\widetilde{BT}^{(y)}  + \Re\left(\frac{1}{S_z}\right)\widetilde{BT}^{(z)} \right) =   \\
&\Re{(s)}\left(\frac{1}{2}  \left(\widetilde{p} + \widetilde{p}^*\right), \frac{1}{2}  \left(\widetilde{F}_p +\widetilde{F}_p^*\right) \right)_{1/\kappa}+ \Re\left(\frac{s^*S_x^*}{S_x}\right)\left(\frac{1}{2S_x} \frac{\partial \left(\widetilde{p} + \widetilde{p}^*\right)}{\partial x}, \frac{1}{S_x}\left( \rho  f_u(x,y,z)\right)\right)_{1/\rho} \\
&+ \Re\left(\frac{s^*S_y^*}{S_y}\right)\left(\frac{1}{2S_y} \frac{\partial \left(\widetilde{p} + \widetilde{p}^*\right)}{\partial y}, \frac{1}{S_y}\left( \rho  f_v(x,y,z)\right)\right)_{1/\rho}   + \Re\left(\frac{s^*S_z^*}{S_z}\right)\left(\frac{1}{2S_z} \frac{\partial  \left(\widetilde{p} + \widetilde{p}^*\right) }{\partial z}, \frac{1}{S_z}\left( \rho  f_w(x,y,z)\right)\right)_{1/\rho} .
%\mathbf{P} s\widetilde{\mathbf{u}} = \frac{1}{S_x} \mathbf{A}\frac{\partial{\widetilde{\mathbf{u}} }}{\partial x} + \frac{1}{S_y} \mathbf{B}\frac{\partial{\widetilde{\mathbf{u}} }}{\partial y} + \mathbf{P}{\mathbf{u}}^0, \quad \Re{s} > 0,
  \end{split}
  \end{equation}
  
  Introducing $\Re{s} = a$, $\Re\left(\frac{\left(sS_j\right)^*}{S_j}\right) =  \gamma_j $ and using Cauchy-Schwarz inequality yields
  %%%%%%%%%%%%%%%%%%%
%%%%%%%%%%%%%%%%%%%
  \[
      \widetilde{E}^2\left( \widetilde{\mathbf{U}}\left(s\right)  \right) + \mathrm{BT}(s)    \le    \widetilde{E}\left( \widetilde{\mathbf{U}}\left(s\right)  \right)\widetilde{E}\left(\widetilde{\mathbf{F}}\left(s\right) \right), \quad \mathrm{BT}(s)   = {\Re\left(\frac{1}{S_x}\right)\widetilde{BT}^{(x)} + \Re\left(\frac{1}{S_y}\right)\widetilde{BT}^{(y)}  + \Re\left(\frac{1}{S_z}\right)\widetilde{BT}^{(z)} } \ge 0.
  \]
   $\hfill\square$
\section{Proof of Theorem \ref{Theo:Numerical_Stability_PML_Laplace}}\label{sec:proof_disc_pml}
%%%%%%%%%%%%%%%%%%%%%
%%%%%%%%%%%%%%%%%%%%%
%%%%%%%%%%%%%%%%%%%%%
Consider the DGSEM approximation \eqref{eq:disc_elemental_pde1_pml_laplace}--\eqref{eq:disc_elemental_pde4_pml_laplace} in two elements, separated at $x = x_2$. Note that $L = 2$, $M = 1$, $N = 1$, so there is only one internal element boundary. Let us denote the solutions in the first element by $\widetilde{p}^{1}(\xi,\eta,\theta,s)$, $\widetilde{u}^{1}(\xi,\eta,\theta, s)$, $\widetilde{v}^{1}(\xi,\eta,\theta, s)$, $\widetilde{w}^{1}(\xi,\eta,\theta, s)$ and the solution the second element denoted by $\widetilde{p}^{2}(\xi,\eta,\theta, s)$, $\widetilde{u}^{2}(\xi,\eta,\theta, s)$, $\widetilde{v}^{2}(\xi,\eta,\theta, s)$, $\widetilde{w}^{2}(\xi,\eta,\theta, s)$. With $\omega_y = 1, \omega_z = 1$, the last two terms in \eqref{eq:disc_elemental_pde1_pml_laplace} vanish, we have
 %%%%%%%
 %%%%%%%
{ 
 \small
\begin{equation}\label{eq:disc_elemental_pde1_pml}
\begin{split}
 &
  s
  %%%%%%%
  \begin{pmatrix} 
  %%%%%%%
 {\boldsymbol{\kappa}}^{-1} & \mathbf{0}\\
 \mathbf{0}  & {\boldsymbol{\kappa}}^{-1}\\
 %%%%%%%
 \end{pmatrix} 
 %%%%%%%
 %%%%%%%
 \begin{bmatrix} 
\widetilde{\boldsymbol{p}}^{1} \\
%%%%%%%
\widetilde{\boldsymbol{p}}^{2}
 \end{bmatrix} 
%%%%%%%
+ \frac{1}{S_x}
%%%%%%%
\begin{pmatrix} 
\mathbf{D}_x & \mathbf{0} \\
 \mathbf{0}  & \mathbf{D}_x\\
 \end{pmatrix} 
 %%%%%%%
 \begin{bmatrix} 
\widetilde{\boldsymbol{u}}^{1} \\
%%%%%%%
\widetilde{\boldsymbol{u}}^{2}
 \end{bmatrix} 
%%%%%%%
%%%%%%%
+ \frac{1}{S_y}
%%%%%%%
\begin{pmatrix} 
\mathbf{D}_y & \mathbf{0} \\
 \mathbf{0}  & \mathbf{D}_y\\
 \end{pmatrix} 
 %%%%%%%
 \begin{bmatrix} 
\widetilde{\boldsymbol{v}}^{1} \\
%%%%%%%
\widetilde{\boldsymbol{v}}^{2}
 \end{bmatrix} 
%%%%%%%
%%%%%%%
+ \frac{1}{S_z}
%%%%%%%
\begin{pmatrix} 
\mathbf{D}_z & \mathbf{0} \\
 \mathbf{0}  & \mathbf{D}_z\\
 \end{pmatrix} 
 %%%%%%%
 \begin{bmatrix} 
\widetilde{\boldsymbol{w}}^{1} \\
%%%%%%%
\widetilde{\boldsymbol{w}}^{2}
 \end{bmatrix} \\
%%%%%%%
&=   
%%%%%%%
 \begin{pmatrix} 
  %%%%%%%
{\boldsymbol{\kappa}}^{-1} & \mathbf{0}\\
\mathbf{0}  & {\boldsymbol{\kappa}}^{-1}\\
 %%%%%%%
 \end{pmatrix} 
 %%%%%%%
   \begin{bmatrix}
 \left( \frac{1}{S_x}{\boldsymbol{f}}_p- \frac{{\boldsymbol{\kappa}}}{sS_yS_x}{\boldsymbol{f}}_{\sigma} -  \frac{{\boldsymbol{\kappa}}}{sS_zS_x}{\boldsymbol{f}}_{\psi}  \right) \\
 %%%%%%
  \left( \frac{1}{S_x}{\boldsymbol{f}}_p- \frac{{\boldsymbol{\kappa}}}{sS_yS_x}{\boldsymbol{f}}_{\sigma} -  \frac{{\boldsymbol{\kappa}}}{sS_zS_x}{\boldsymbol{f}}_{\psi}  \right)
  \end{bmatrix} \\
- &
%%%%%
\frac{1+r_x}{2ZS_x}
%%%%%
\begin{pmatrix} 
\mathbf{H}_x^{-1} & \mathbf{0} \\
 \mathbf{0}  & \mathbf{H}_x^{-1}\\
 \end{pmatrix} 
 %%%%%%%
 \begin{pmatrix} 
\mathbf{B}_x{(-1,-1)}  & \mathbf{0} \\
 \mathbf{0}  & \mathbf{B}_x{(1,1)}\\
 \end{pmatrix} 
 %%%%%%%
  %%%%%%%
 %%%%%%%
 \begin{bmatrix} 
\widetilde{\boldsymbol{p}}^{1} \\
%%%%%%%
\widetilde{\boldsymbol{p}}^{2}
 \end{bmatrix} 
%%%%%%%
-
%%%%%
\frac{1-r_x}{2S_x}
%%%%%
\begin{pmatrix} 
\mathbf{H}_x^{-1} & \mathbf{0} \\
 \mathbf{0}  & \mathbf{H}_x^{-1}\\
 \end{pmatrix} 
 %%%%%%%
 \begin{pmatrix} 
-\mathbf{B}_x{(-1,-1)}  & \mathbf{0} \\
 \mathbf{0}  & \mathbf{B}_x{(1,1)}\\
 \end{pmatrix} 
 %%%%%%%
  %%%%%%%
 %%%%%%%
 \begin{bmatrix} 
\widetilde{\boldsymbol{u}}^{1} \\
%%%%%%%
\widetilde{\boldsymbol{u}}^{2}
 \end{bmatrix} 
%%%%%%%
 %%%%%%%
\\
 %%%%%
 - &
%%%%%
\frac{1+r_y}{2ZS_y}
%%%%%
\begin{pmatrix} 
\mathbf{H}_y^{-1} & \mathbf{0} \\
 \mathbf{0}  & \mathbf{H}_y^{-1}\\
 \end{pmatrix} 
 %%%%%%%
 \begin{pmatrix} 
\mathbf{B}_y{(-1,-1)}  & \mathbf{0} \\
 \mathbf{0}  & \mathbf{B}_y{(1,1)}\\
 \end{pmatrix} 
 %%%%%%%
  %%%%%%%
 %%%%%%%
 \begin{bmatrix} 
\widetilde{\boldsymbol{p}}^{1} \\
%%%%%%%
\widetilde{\boldsymbol{p}}^{2}
 \end{bmatrix} 
%%%%%%%
-
%%%%%
\frac{1-r_y}{2S_y}
%%%%%
\begin{pmatrix} 
\mathbf{H}_y^{-1} & \mathbf{0} \\
 \mathbf{0}  & \mathbf{H}_y^{-1}\\
 \end{pmatrix} 
 %%%%%%%
 \begin{pmatrix} 
-\mathbf{B}_y{(-1,-1)}  & \mathbf{0} \\
 \mathbf{0}  & \mathbf{B}_y{(1,1)}\\
 \end{pmatrix} 
 %%%%%%%
  %%%%%%%
 %%%%%%%
 \begin{bmatrix} 
\widetilde{\boldsymbol{u}}^{1} \\
%%%%%%%
\widetilde{\boldsymbol{u}}^{2}
 \end{bmatrix} 
%%%%%%%
 %%%%%%%
\\
 %%%%%
- &
%%%%%
\frac{1+r_z}{2ZS_z}
%%%%%
\begin{pmatrix} 
\mathbf{H}_z^{-1} & \mathbf{0} \\
 \mathbf{0}  & \mathbf{H}_z^{-1}\\
 \end{pmatrix} 
 %%%%%%%
 \begin{pmatrix} 
\mathbf{B}_z{(-1,-1)}  & \mathbf{0} \\
 \mathbf{0}  & \mathbf{B}_z{(1,1)}\\
 \end{pmatrix} 
 %%%%%%%
  %%%%%%%
 %%%%%%%
 \begin{bmatrix} 
\widetilde{\boldsymbol{p}}^{1} \\
%%%%%%%
\widetilde{\boldsymbol{p}}^{2}
 \end{bmatrix} 
%%%%%%%
-
%%%%%
\frac{1-r_z}{2S_z}
%%%%%
\begin{pmatrix} 
\mathbf{H}_z^{-1} & \mathbf{0} \\
 \mathbf{0}  & \mathbf{H}_z^{-1}\\
 \end{pmatrix} 
 %%%%%%%
 \begin{pmatrix} 
-\mathbf{B}_z{(-1,-1)}  & \mathbf{0} \\
 \mathbf{0}  & \mathbf{B}_z{(1,1)}\\
 \end{pmatrix} 
 %%%%%%%
  %%%%%%%
 %%%%%%%
 \begin{bmatrix} 
\widetilde{\boldsymbol{u}}^{1} \\
%%%%%%%
\widetilde{\boldsymbol{u}}^{2}
 \end{bmatrix} 
%%%%%%%
 %%%%%%%
\\
 %%%%%
 - &
%%%%%
\frac{1}{2ZS_x}
%%%%%
\begin{pmatrix} 
\mathbf{H}_x^{-1} & \mathbf{0} \\
 \mathbf{0}  & \mathbf{H}_x^{-1}\\
 \end{pmatrix} 
 %%%%%%%
 \begin{pmatrix} 
\mathbf{B}_x{(1,1)}  & -\mathbf{B}_x{(1,-1)} \\
-\mathbf{B}_x^T{(1,-1)}   & \mathbf{B}_x{(-1,-1)}\\
 \end{pmatrix} 
 %%%%%%%
  %%%%%%%
 %%%%%%%
 \begin{bmatrix} 
\widetilde{\boldsymbol{p}}^{1} \\
%%%%%%%
\widetilde{\boldsymbol{p}}^{2}
 \end{bmatrix} 
%%%%%%%
-
%%%%%
\frac{1}{2S_x}
%%%%%
\begin{pmatrix} 
\mathbf{H}_x^{-1} & \mathbf{0} \\
 \mathbf{0}  & \mathbf{H}_x^{-1}\\
 \end{pmatrix} 
 %%%%%%%
 \begin{pmatrix} 
-\mathbf{B}_x{(1,1)}  & \mathbf{B}_x{(1,-1)} \\
-\mathbf{B}_x^T{(1,-1)}  & \mathbf{B}_x{(-1,-1)}\\
 \end{pmatrix} 
 %%%%%%%
  %%%%%%%
 %%%%%%%
 \begin{bmatrix} 
\widetilde{\boldsymbol{u}}^{1} \\
%%%%%%%
\widetilde{\boldsymbol{u}}^{2}
 \end{bmatrix} 
%%%%%%%
 %%%%%%%
%\\
% %%%%%
% -&
% \frac{1}{S_y}\mathbf{H}_y^{-1}\left(\mathbf{B}_y(-1) \boldsymbol{Z}^{-1} \left(\frac{1-r_x}{2}\boldsymbol{Z} \widetilde{\boldsymbol{v}}( s) + \frac{1+r_x}{2}\widetilde{\boldsymbol{p}}( s)  \right) - \mathbf{B}_y(1) \boldsymbol{Z}^{-1}\left( \frac{1-r_x}{2}\boldsymbol{Z} \widetilde{\boldsymbol{v}}( s) - \frac{1+r_x}{2}\widetilde{\boldsymbol{p}}( s)  \right) \right)\\
%%%%%%
% -&\frac{1}{S_z}\mathbf{H}_z^{-1}\left(\mathbf{B}_z(-1) \boldsymbol{Z}^{-1} \left(\frac{1-r_z}{2}\boldsymbol{Z} \widetilde{\boldsymbol{w}}( s) + \frac{1+r_z}{2}\widetilde{\boldsymbol{p}}( s)  \right) - \mathbf{B}_z(1) \boldsymbol{Z}^{-1}\left( \frac{1-r_z}{2}\boldsymbol{Z} \widetilde{\boldsymbol{w}}( s) - \frac{1+r_z}{2}\widetilde{\boldsymbol{p}}( s)  \right) \right)\\
\end{split}
\end{equation}
%%%%%%%
%%%%%%%
 \begin{equation}\label{eq:disc_elemental_pde2_pml}
\begin{split}
&
s
  %%%%%%%
  \begin{pmatrix} 
  %%%%%%%
 {\boldsymbol{\rho}} & \mathbf{0}\\
 \mathbf{0}  & {\boldsymbol{\rho}}\\
 %%%%%%%
 \end{pmatrix} 
 %%%%%%%
 %%%%%%%
 \begin{bmatrix} 
\widetilde{\boldsymbol{u}}^{1} \\
%%%%%%%
\widetilde{\boldsymbol{u}}^{2}
 \end{bmatrix} 
%%%%%%%
+ \frac{1}{S_x}
%%%%%%%
\begin{pmatrix} 
\mathbf{D}_x & \mathbf{0} \\
 \mathbf{0}  & \mathbf{D}_x\\
 \end{pmatrix} 
 %%%%%%%
 \begin{bmatrix} 
\widetilde{\boldsymbol{p}}^{1} \\
%%%%%%%
\widetilde{\boldsymbol{p}}^{2}
 \end{bmatrix}  
%%%%%%%
=
%%%%%
\begin{bmatrix} 
\frac{\boldsymbol{\rho}}{S_x} \boldsymbol{f}_u \\
%%%%%%%
\frac{\boldsymbol{\rho}}{S_x} \boldsymbol{f}_u
 \end{bmatrix} 
 %%%%%
 %%%%%
 \\
 %%%%%
 %%%%%
& 
-\frac{1+r_x}{2S_x}
%%%%%
\begin{pmatrix} 
\mathbf{H}_x^{-1} & \mathbf{0} \\
 \mathbf{0}  & \mathbf{H}_x^{-1}\\
 \end{pmatrix} 
 %%%%%%%
 \begin{pmatrix} 
\mathbf{B}_x{(-1,-1)}  & \mathbf{0} \\
 \mathbf{0}  & -\mathbf{B}_x{(1,1)}\\
 \end{pmatrix} 
 %%%%%%%
  %%%%%%%
 %%%%%%%
 \begin{bmatrix} 
\widetilde{\boldsymbol{p}}^{1} \\
%%%%%%%
\widetilde{\boldsymbol{p}}^{2}
 \end{bmatrix} 
%%%%%%%
-
%%%%%
Z\frac{1-r_x}{2S_x}
%%%%%
\begin{pmatrix} 
\mathbf{H}_x^{-1} & \mathbf{0} \\
 \mathbf{0}  & \mathbf{H}_x^{-1}\\
 \end{pmatrix} 
 %%%%%%%
 \begin{pmatrix} 
\mathbf{B}_x{(-1,-1)}  & \mathbf{0} \\
 \mathbf{0}  & \mathbf{B}_x{(1,1)}\\
 \end{pmatrix} 
 %%%%%%%
  %%%%%%%
 %%%%%%%
 \begin{bmatrix} 
\widetilde{\boldsymbol{u}}^{1} \\
%%%%%%%
\widetilde{\boldsymbol{u}}^{2}
 \end{bmatrix} 
%%%%%%%
 %%%%%%%
\\
 &- 
%%%%%
\frac{1}{2S_x}
%%%%%
\begin{pmatrix} 
\mathbf{H}_x^{-1} & \mathbf{0} \\
 \mathbf{0}  & \mathbf{H}_x^{-1}\\
 \end{pmatrix} 
 %%%%%%%
 \begin{pmatrix} 
-\mathbf{B}_x{(1,1)} & \mathbf{B}_x{(1,-1)}  \\
-\mathbf{B}_x^T{(1,-1)}  & \mathbf{B}_x{(-1,-1)}\\
 \end{pmatrix} 
 %%%%%%%
  %%%%%%%
 %%%%%%%
 \begin{bmatrix} 
\widetilde{\boldsymbol{p}}^{1} \\
%%%%%%%
\widetilde{\boldsymbol{p}}^{2}
 \end{bmatrix} 
%%%%%%%
-
%%%%%
\frac{Z}{2S_x}
%%%%%
\begin{pmatrix} 
\mathbf{H}_x^{-1} & \mathbf{0} \\
 \mathbf{0}  & \mathbf{H}_x^{-1}\\
 \end{pmatrix} 
 %%%%%%%
 \begin{pmatrix} 
\mathbf{B}_x{(1,1)} & -\mathbf{B}_x{(1,-1)} \\
-\mathbf{B}_x^T{(1,-1)}  & \mathbf{B}_x{(-1,-1)} \\
 \end{pmatrix} 
 %%%%%%%
  %%%%%%%
 %%%%%%%
 \begin{bmatrix} 
\widetilde{\boldsymbol{u}}^{1} \\
%%%%%%%
\widetilde{\boldsymbol{u}}^{2}
 \end{bmatrix} 
%%%%%%%
 %%%%%%%
\end{split}
\end{equation}
%%%%%%
\begin{equation}\label{eq:disc_elemental_pde3_pml}
\begin{split}
&
s
  %%%%%%%
  \begin{pmatrix} 
  %%%%%%%
 {\boldsymbol{\rho}} & \mathbf{0}\\
 \mathbf{0}  & {\boldsymbol{\rho}}\\
 %%%%%%%
 \end{pmatrix} 
 %%%%%%%
 %%%%%%%
 \begin{bmatrix} 
\widetilde{\boldsymbol{v}}^{1} \\
%%%%%%%
\widetilde{\boldsymbol{v}}^{2}
 \end{bmatrix} 
%%%%%%%
+ \frac{1}{S_y}
%%%%%%%
\begin{pmatrix} 
\mathbf{D}_y & \mathbf{0} \\
 \mathbf{0}  & \mathbf{D}_y\\
 \end{pmatrix} 
 %%%%%%%
 \begin{bmatrix} 
\widetilde{\boldsymbol{p}}^{1} \\
%%%%%%%
\widetilde{\boldsymbol{p}}^{2}
 \end{bmatrix}  
%%%%%%%
=
%%%%%
\begin{bmatrix} 
\frac{\boldsymbol{\rho}}{S_y} \boldsymbol{f}_v \\
%%%%%%%
\frac{\boldsymbol{\rho}}{S_y} \boldsymbol{f}_v
 \end{bmatrix} 
 %%%%%
 %%%%%
 \\
 %%%%%
 %%%%%
& 
-\frac{1+r_y}{2S_y}
%%%%%
\begin{pmatrix} 
\mathbf{H}_x^{-1} & \mathbf{0} \\
 \mathbf{0}  & \mathbf{H}_x^{-1}\\
 \end{pmatrix} 
 %%%%%%%
 \begin{pmatrix} 
\mathbf{B}_y{(-1,-1)} & \mathbf{0} \\
 \mathbf{0}  & -\mathbf{B}_y{(1,1)}\\
 \end{pmatrix} 
 %%%%%%%
  %%%%%%%
 %%%%%%%
 \begin{bmatrix} 
\widetilde{\boldsymbol{p}}^{1} \\
%%%%%%%
\widetilde{\boldsymbol{p}}^{2}
 \end{bmatrix} 
%%%%%%%
-
%%%%%
Z\frac{1-r_y}{2S_y}
%%%%%
\begin{pmatrix} 
\mathbf{H}_y^{-1} & \mathbf{0} \\
 \mathbf{0}  & \mathbf{H}_y^{-1}\\
 \end{pmatrix} 
 %%%%%%%
 \begin{pmatrix} 
\mathbf{B}_y{(-1,-1)}  & \mathbf{0} \\
 \mathbf{0}  & \mathbf{B}_y{(1,1)}\\
 \end{pmatrix} 
 %%%%%%%
  %%%%%%%
 %%%%%%%
 \begin{bmatrix} 
\widetilde{\boldsymbol{v}}^{1} \\
%%%%%%%
\widetilde{\boldsymbol{v}}^{2}
 \end{bmatrix} 
%%%%%%%
 %%%%%%%
\end{split}
\end{equation}
%%%%%%
\begin{equation}\label{eq:disc_elemental_pde4_pml}
\begin{split}
&
s
  %%%%%%%
  \begin{pmatrix} 
  %%%%%%%
 {\boldsymbol{\rho}} & \mathbf{0}\\
 \mathbf{0}  & {\boldsymbol{\rho}}\\
 %%%%%%%
 \end{pmatrix} 
 %%%%%%%
 %%%%%%%
 \begin{bmatrix} 
\widetilde{\boldsymbol{w}}^{1} \\
%%%%%%%
\widetilde{\boldsymbol{w}}^{2}
 \end{bmatrix} 
%%%%%%%
+ \frac{1}{S_z}
%%%%%%%
\begin{pmatrix} 
\mathbf{D}_z & \mathbf{0} \\
 \mathbf{0}  & \mathbf{D}_z\\
 \end{pmatrix} 
 %%%%%%%
 \begin{bmatrix} 
\widetilde{\boldsymbol{p}}^{1} \\
%%%%%%%
\widetilde{\boldsymbol{p}}^{2}
 \end{bmatrix}  
%%%%%%%
=
%%%%%
\begin{bmatrix} 
\frac{\boldsymbol{\rho}}{S_z} \boldsymbol{f}_w \\
%%%%%%%
\frac{\boldsymbol{\rho}}{S_z} \boldsymbol{f}_w
 \end{bmatrix} 
 %%%%%
 %%%%%
 \\
 %%%%%
 %%%%%
& 
-\frac{1+r_z}{2S_z}
%%%%%
\begin{pmatrix} 
\mathbf{H}_z^{-1} & \mathbf{0} \\
 \mathbf{0}  & \mathbf{H}_z^{-1}\\
 \end{pmatrix} 
 %%%%%%%
 \begin{pmatrix} 
\mathbf{B}_z{(-1,-1)}  & \mathbf{0} \\
 \mathbf{0}  & -\mathbf{B}_z{(1,1)}\\
 \end{pmatrix} 
 %%%%%%%
  %%%%%%%
 %%%%%%%
 \begin{bmatrix} 
\widetilde{\boldsymbol{p}}^{1} \\
%%%%%%%
\widetilde{\boldsymbol{p}}^{2}
 \end{bmatrix} 
%%%%%%%
-
%%%%%
Z\frac{1-r_z}{2S_z}
%%%%%
\begin{pmatrix} 
\mathbf{H}_z^{-1} & \mathbf{0} \\
 \mathbf{0}  & \mathbf{H}_z^{-1}\\
 \end{pmatrix} 
 %%%%%%%
 \begin{pmatrix} 
\mathbf{B}_z{(-1,-1)}& \mathbf{0} \\
 \mathbf{0}  & \mathbf{B}_z{(1,1)}\\
 \end{pmatrix} 
 %%%%%%%
  %%%%%%%
 %%%%%%%
 \begin{bmatrix} 
\widetilde{\boldsymbol{w}}^{1} \\
%%%%%%%
\widetilde{\boldsymbol{w}}^{2}
 \end{bmatrix} 
%%%%%%%
 %%%%%%%
 %%%%%%%
\end{split}
\end{equation}
%%%%%%%
% \begin{equation}\label{eq:disc_elemental_pde3_pml}
%\begin{split}
%& {\boldsymbol{\rho}} s \widetilde{\boldsymbol{v}}(s)   + \frac{1}{S_y}\mathbf{D}_y \widetilde{\boldsymbol{p}}(s) =  \frac{\boldsymbol{\rho}}{S_y} \boldsymbol{f}_v  \\
% -&\frac{1}{S_y}\mathbf{H}_y^{-1}\left(\mathbf{B}_y(-1)  \left(\frac{1-r_x}{2}\boldsymbol{Z} \widetilde{\boldsymbol{v}}( s) + \frac{1+r_x}{2}\widetilde{\boldsymbol{p}}( s)  \right) 
% + 
% \mathbf{B}_y(1) \left( \frac{1-r_x}{2}\boldsymbol{Z} \widetilde{\boldsymbol{v}}( s) - \frac{1+r_x}{2}\widetilde{\boldsymbol{p}}( s)  \right) \right)\\
%\end{split}
%\end{equation}
%%%%%%%
%%%%%%%
%\begin{equation}\label{eq:disc_elemental_pde4_pml}
%\begin{split}
%& {\boldsymbol{\rho}} s \widetilde{\boldsymbol{w}}(s)   + \frac{1}{S_z} \mathbf{D}_z \widetilde{\boldsymbol{p}}(s) =   \frac{\boldsymbol{\rho}}{S_z} \boldsymbol{f}_w  \\
%  -&\frac{1}{S_z}\mathbf{H}_z^{-1}\left(\mathbf{B}_z(-1)  \left(\frac{1-r_z}{2}\boldsymbol{Z} \widetilde{\boldsymbol{w}}( s) + \frac{1+r_z}{2}\widetilde{\boldsymbol{p}}( s)  \right) 
%  +
%   \mathbf{B}_z(1) \left( \frac{1-r_z}{2}\boldsymbol{Z} \widetilde{\boldsymbol{w}}( s) - \frac{1+r_z}{2}\widetilde{\boldsymbol{p}}( s)  \right) \right)\\
%\end{split}
%\end{equation}
%%%%%%
%%%%%%
%%%%%%
}

We will now eliminate the velocity fields, we have

 { 
 \small
\begin{equation}\label{eq:disc_elemental_pde_secondorder_pml_appendix}
\begin{split}
 &
  ss^*
 %%%%%%%
 \begin{pmatrix} 
\mathbf{H} {\boldsymbol{\kappa}}^{-1}& \mathbf{0} \\
 \mathbf{0}  & \mathbf{H} {\boldsymbol{\kappa}}^{-1}\\
 \end{pmatrix}  
 %%%%%%%
 \begin{bmatrix} 
s\widetilde{\boldsymbol{p}}^{1} \\
%%%%%%%
s\widetilde{\boldsymbol{p}}^{2}
 \end{bmatrix} 
%%%%%%%
+ 
\frac{(s^*S_x^*)}{S_x} 
%%%%%%%
\widetilde{\mathbf{D}}_x^{\dagger}  
\begin{pmatrix} 
\boldsymbol{\rho}^{-1} & \mathbf{0} \\
 \mathbf{0}  & \boldsymbol{\rho}^{-1} \\
 \end{pmatrix} 
 %%%%%%%
  \widetilde{\mathbf{H}}_x
 \widetilde{\mathbf{D}}_x
 \begin{bmatrix} 
\widetilde{\boldsymbol{p}}^{1} \\
%%%%%%%
\widetilde{\boldsymbol{p}}^{2}
 \end{bmatrix} 
%%%%%%%
%%%%%%%
+ 
\frac{(sS_y)^*}{S_y} 
%%%%%%%
\widetilde{\mathbf{D}}_y^{\dagger}  
\begin{pmatrix} 
\boldsymbol{\rho}^{-1} & \mathbf{0} \\
 \mathbf{0}  & \boldsymbol{\rho}^{-1} \\
 \end{pmatrix} 
 %%%%%%%
  \widetilde{\mathbf{H}}_y
 \widetilde{\mathbf{D}}_y
 \begin{bmatrix} 
\widetilde{\boldsymbol{p}}^{1} \\
%%%%%%%
\widetilde{\boldsymbol{p}}^{2}
 \end{bmatrix} 
%%%%%%%
%%%%%%%
\\
&+ 
\frac{\left(sS_z\right)^*}{S_z} 
%%%%%%%
\widetilde{\mathbf{D}}_z^{\dagger}  
\begin{pmatrix} 
\boldsymbol{\rho}^{-1} & \mathbf{0} \\
 \mathbf{0}  & \boldsymbol{\rho}^{-1} \\
 \end{pmatrix} 
 %%%%%%%
  \widetilde{\mathbf{H}}_z
 \widetilde{\mathbf{D}}_z
 \begin{bmatrix} 
\widetilde{\boldsymbol{p}}^{1} \\
%%%%%%%
\widetilde{\boldsymbol{p}}^{2}
 \end{bmatrix} 
 +
 \frac{1}{2S_xZ}
%%%%%
\begin{pmatrix} 
\mathbf{H}_y\mathbf{H}_z & \mathbf{0} \\
 \mathbf{0}  & \mathbf{H}_y\mathbf{H}_z\\
 \end{pmatrix} 
 %%%%%%%
 \begin{pmatrix} 
\mathbf{B}_x{(1,1)}  & -\mathbf{B}_x{(1,-1)}  \\
-\mathbf{B}_x^T{(1,-1)}   & \mathbf{B}_x{(-1,-1)}\\
 \end{pmatrix} 
 %%%%%%%
  %%%%%%%
 %%%%%%%
 \begin{bmatrix} 
\widetilde{\boldsymbol{p}}^{1} \\
%%%%%%%
\widetilde{\boldsymbol{p}}^{2}
 \end{bmatrix} 
%%%%%%%
  \\
%%%%%%%
&
\frac{1+r_x}{2ZS_x}
%%%%%
\begin{pmatrix} 
\mathbf{H}_y\mathbf{H}_z & \mathbf{0} \\
 \mathbf{0}  & \mathbf{H}_y\mathbf{H}_z\\
 \end{pmatrix}  
 %%%%%%%
 \begin{pmatrix} 
\mathbf{B}_z{(-1,-1)}  & \mathbf{0} \\
 \mathbf{0}  & \mathbf{B}_z{(1,1)}\\
 \end{pmatrix} 
 %%%%%%%
  %%%%%%%
 %%%%%%%
 \begin{bmatrix} 
\widetilde{\boldsymbol{p}}^{1} \\
%%%%%%%
\widetilde{\boldsymbol{p}}^{2}
 \end{bmatrix} 
%%%%%%%
+
\frac{1+r_y}{2ZS_y}
%%%%%
\begin{pmatrix} 
\mathbf{H}_x\mathbf{H}_z & \mathbf{0} \\
 \mathbf{0}  & \mathbf{H}_x\mathbf{H}_z\\
 \end{pmatrix} 
 %%%%%%%
 \begin{pmatrix} 
\mathbf{B}_y{(-1,-1)} & \mathbf{0} \\
 \mathbf{0}  & \mathbf{B}_y{(1,1)}\\
 \end{pmatrix} 
 %%%%%%%
  %%%%%%%
 %%%%%%%
 \begin{bmatrix} 
\widetilde{\boldsymbol{p}}^{1} \\
%%%%%%%
\widetilde{\boldsymbol{p}}^{2}
 \end{bmatrix} 
%%%%%%%
\\
&
+
\frac{1+r_z}{2ZS_z}
%%%%%
\begin{pmatrix} 
\mathbf{H}_x\mathbf{H}_y & \mathbf{0} \\
 \mathbf{0}  & \mathbf{H}_x \mathbf{H}_y\\
 \end{pmatrix} 
 %%%%%%%
 \begin{pmatrix} 
\mathbf{B}_z{(-1,-1)}  & \mathbf{0} \\
 \mathbf{0}  & \mathbf{B}_z{(1,1)}\\
 \end{pmatrix} 
 %%%%%%%
  %%%%%%%
 %%%%%%%
 \begin{bmatrix} 
\widetilde{\boldsymbol{p}}^{1} \\
%%%%%%%
\widetilde{\boldsymbol{p}}^{2}
 \end{bmatrix} 
%%%%%%%
\\
&
=   
%%%%%%%
 \begin{pmatrix} 
  %%%%%%%
\mathbf{H} {\boldsymbol{\kappa}}^{-1}& \mathbf{0} \\
 \mathbf{0}  & \mathbf{H} {\boldsymbol{\kappa}}^{-1}\\
 %%%%%%%
 \end{pmatrix} 
 %%%%%%%
   \begin{bmatrix}
 \left( \frac{1}{S_x}{\boldsymbol{f}}_p- \frac{{\boldsymbol{\kappa}}}{sS_yS_x}{\boldsymbol{f}}_{\sigma} -  \frac{{\boldsymbol{\kappa}}}{sS_zS_x}{\boldsymbol{f}}_{\psi}  \right) \\
 %%%%%%
  \left( \frac{1}{S_x}{\boldsymbol{f}}_p- \frac{{\boldsymbol{\kappa}}}{sS_yS_x}{\boldsymbol{f}}_{\sigma} -  \frac{{\boldsymbol{\kappa}}}{sS_zS_x}{\boldsymbol{f}}_{\psi}  \right)
  \end{bmatrix} 
  +
\frac{(s^*S_x^*)}{S_x} 
%%%%%%%
\widetilde{\mathbf{D}}_x^{\dagger}  
\begin{pmatrix} 
\boldsymbol{\rho}^{-1} & \mathbf{0} \\
 \mathbf{0}  & \boldsymbol{\rho}^{-1} \\
 \end{pmatrix} 
 %%%%%%%
  \widetilde{\mathbf{H}}_x
 \begin{bmatrix} 
\frac{\boldsymbol{\rho}}{S_x} \boldsymbol{f}_u \\
%%%%%%%
\frac{\boldsymbol{\rho}}{S_x} \boldsymbol{f}_u
 \end{bmatrix} 
%%%%%%%
\\
&
+ 
\frac{(s^*S_y^*)}{S_y} 
%%%%%%%
\widetilde{\mathbf{D}}_y^{\dagger}  
\begin{pmatrix} 
\boldsymbol{\rho}^{-1} & \mathbf{0} \\
 \mathbf{0}  & \boldsymbol{\rho}^{-1} \\
 \end{pmatrix} 
 %%%%%%%
  \widetilde{\mathbf{H}}_y
 \begin{bmatrix} 
\frac{\boldsymbol{\rho}}{S_y} \boldsymbol{f}_v \\
%%%%%%%
\frac{\boldsymbol{\rho}}{S_y} \boldsymbol{f}_v
 \end{bmatrix} 
%%%%%%%
+
\frac{(s^*S_z^*)}{S_z} 
%%%%%%%
\widetilde{\mathbf{D}}_z^{\dagger} 
\begin{pmatrix} 
\boldsymbol{\rho}^{-1} & \mathbf{0} \\
 \mathbf{0}  & \boldsymbol{\rho}^{-1} \\
 \end{pmatrix} 
 %%%%%%%
  \widetilde{\mathbf{H}}_x
 \begin{bmatrix} 
\frac{\boldsymbol{\rho}}{S_z} \boldsymbol{f}_w\\
%%%%%%%
\frac{\boldsymbol{\rho}}{S_z} \boldsymbol{f}_w
 \end{bmatrix} 
%%%%%%%
\end{split}
\end{equation},
}
where
{
\small
\begin{align*}
 &\widetilde{\mathbf{D}}_x =  \frac{1}{S_x} \left( \begin{pmatrix} 
\mathbf{D}_z & \mathbf{0} \\
 \mathbf{0}  & \mathbf{D}_x\\
 \end{pmatrix}
 + 
%%%%%
\begin{pmatrix} 
\mathbf{H}_x^{-1} & \mathbf{0} \\
 \mathbf{0}  & \mathbf{H}_x^{-1}\\
 \end{pmatrix} 
 %%%%%%%
 \left(
 \frac{1+r_x}{2}
 \begin{pmatrix} 
\mathbf{B}_x{(-1,-1)}  & \mathbf{0} \\
 \mathbf{0}  & -\mathbf{B}_x{(1,1)}\\
 \end{pmatrix} 
 +
 \frac{1}{2}
%%%%%
 %%%%%
 \begin{pmatrix} 
-\mathbf{B}_x{(1,1)} &  \mathbf{B}_x{(1,-1)} \\
-\mathbf{B}_x^T{(1,-1)}   & \mathbf{B}_x{(-1,-1)} \\
 \end{pmatrix} 
 \right)
 \right)
 \\
  &\widetilde{\mathbf{D}}_y =  \frac{1}{S_y} \left( \begin{pmatrix} 
\mathbf{D}_y & \mathbf{0} \\
 \mathbf{0}  & \mathbf{D}_y\\
 \end{pmatrix}
 + 
%%%%%
\begin{pmatrix} 
\mathbf{H}_y^{-1} & \mathbf{0} \\
 \mathbf{0}  & \mathbf{H}_y^{-1}\\
 \end{pmatrix} 
 %%%%%%%
 \left(
 \frac{1+r_y}{2}
 \begin{pmatrix} 
\mathbf{B}_y{(-1,-1)}  & \mathbf{0} \\
 \mathbf{0}  & -\mathbf{B}_y{(1,1)}\\
 \end{pmatrix} 
 \right)
 \right)
 \\
  &\widetilde{\mathbf{D}}_z =  \frac{1}{S_z} \left( \begin{pmatrix} 
\mathbf{D}_z & \mathbf{0} \\
 \mathbf{0}  & \mathbf{D}_z\\
 \end{pmatrix}
 + 
%%%%%
\begin{pmatrix} 
\mathbf{H}_z^{-1} & \mathbf{0} \\
 \mathbf{0}  & \mathbf{H}_z^{-1}\\
 \end{pmatrix} 
 %%%%%%%
 \left(
 \frac{1+r_z}{2}
 \begin{pmatrix} 
\mathbf{B}_z{(-1,-1)}  & \mathbf{0} \\
 \mathbf{0}  & -\mathbf{B}_z{(1,1)}\\
 \end{pmatrix} 
 \right)
 \right)
\end{align*}
}

{
\small
\begin{align*}
 &\widetilde{\mathbf{H}}_x(s, d_x)  =   \begin{pmatrix} 
\mathbf{H} & \mathbf{0} \\
 \mathbf{0}  & \mathbf{H}\\
 \end{pmatrix}
\left( \mathbf{I} + 
%%%%%
\begin{pmatrix} 
\mathbf{H}_x^{-1} & \mathbf{0} \\
 \mathbf{0}  & \mathbf{H}_x^{-1}\\
 \end{pmatrix} 
 %%%%%%%
 \left(
 \frac{(1- r_x)c}{2 sS_x}
 \begin{pmatrix} 
\mathbf{B}_x{(-1,-1)}  & \mathbf{0} \\
 \mathbf{0}  & \mathbf{B}_x{(1,1)} \\
 \end{pmatrix} 
 +
 \frac{c}{2sS_x}
%%%%%
 %%%%%
 \begin{pmatrix} 
\mathbf{B}_x{(1,1)}  & - \mathbf{B}_x{(1,-1)} \\
-\mathbf{B}_x^T{(1,-1)}   & \mathbf{B}_x{(-1,-1)} \\
 \end{pmatrix} 
 \right)
 \right)^{-1}
 \\
 &\widetilde{\mathbf{H}}_y(s, d_y)  =   \begin{pmatrix} 
\mathbf{H} & \mathbf{0} \\
 \mathbf{0}  & \mathbf{H}\\
 \end{pmatrix}
\left( \mathbf{I} + 
%%%%%
\begin{pmatrix} 
\mathbf{H}_y^{-1} & \mathbf{0} \\
 \mathbf{0}  & \mathbf{H}_y^{-1}\\
 \end{pmatrix} 
 %%%%%%%
 \left(
 \frac{(1- r_y)c}{2 sS_y}
 \begin{pmatrix} 
\mathbf{B}_y{(-1,-1)}  & \mathbf{0} \\
 \mathbf{0}  & \mathbf{B}_y{(1,1)} \\
 \end{pmatrix} 
 \right)
 \right)^{-1}
 \\
 &\widetilde{\mathbf{H}}_z(s, d_z) =   \begin{pmatrix} 
\mathbf{H} & \mathbf{0} \\
 \mathbf{0}  & \mathbf{H}\\
 \end{pmatrix}
\left( \mathbf{I} + 
%%%%%
\begin{pmatrix} 
\mathbf{H}_z^{-1} & \mathbf{0} \\
 \mathbf{0}  & \mathbf{H}_z^{-1}\\
 \end{pmatrix} 
 %%%%%%%
 \left(
 \frac{(1- r_z)c}{2 sS_z}
 \begin{pmatrix} 
\mathbf{B}_z{(-1,-1)}  & \mathbf{0} \\
 \mathbf{0}  & \mathbf{B}_z{(1,1)} \\
 \end{pmatrix} 
 \right)
 \right)^{-1}
 \end{align*}
}

Multiply equation \eqref{eq:disc_elemental_pde_secondorder_pml_appendix} by 
$
 \begin{bmatrix} 
s\widetilde{\boldsymbol{p}}^{1} \\
%%%%%%%
s\widetilde{\boldsymbol{p}}^{2}
 \end{bmatrix}^{\boldsymbol{\dagger}}$ 
from the left, and add complex conjugate transpose. Introducing the scalar product \eqref{eq:scalar_product_discrete_laplace} and  the corresponding norm \eqref{eq:norm_discrete_laplace}, we have
%%%%%%%%%
%%%%%%%%%
{\small
\begin{equation}\label{eq:energy_estimate_pml_laplace_corner_appendix}
  \widetilde{\mathcal{E}}^2\left(\widetilde{\mathbf{U}}\left(s\right) \right) +\widetilde{\mathrm{BT}}(s) \le   \widetilde{\mathcal{E}}\left(\widetilde{\mathbf{U}}\left(s\right) \right)\widetilde{\mathcal{E}}\left(\widetilde{\mathbf{F}}\left(s\right) \right), \quad \widetilde{\mathrm{BT}}(s)   = {\Re\left(\frac{1}{S_x}\right)\widetilde{\mathbf{BT}}^{(x)} + \Re\left(\frac{1}{S_y}\right) \widetilde{\mathbf{BT}}^{(y)}  + \Re\left(\frac{1}{S_z}\right)\widetilde{\mathbf{BT}}^{(z)}  + \Re\left(\frac{1}{S_x}\right) \widetilde{\mathbf{IT}}^{(x)}} \ge 0.
\end{equation}
}

\end{document}